\documentclass [11pt]{article}
\usepackage {amsfonts}
\usepackage {mathrsfs}
\usepackage {amsmath}
\usepackage {latexsym}
\usepackage {amssymb}
\usepackage {amsthm}
\usepackage {amscd}
\date{}
\renewcommand{\baselinestretch}{1.2}
\title{\bf Classification of quasi-affine Generalized Dynkin Diagrams with Rank $> 5$ }
\author{ \small Zhengtang Tan $^{a}$,  Shouchuan Zhang $^{b}$   \\
\small $a$.School of Engineering and Design, Hunan Normal University\\ Changsha  410081,   P.R. China \\
\small $b$. Department  of Mathematics,   Hunan University\\
Changsha  410082,  P.R. China\\
\small {\tt Emails:  z9491@sina.cn (SZ);   1843186255@qq.com (ZTT)} }
\date{}

\begin{document}
\newtheorem{Proposition}{Proposition}[section]
\newtheorem{Theorem}[Proposition]{Theorem}
\newtheorem{Definition}[Proposition]{Definition}
\newtheorem{Corollary}[Proposition]{Corollary}
\newtheorem{Lemma}[Proposition]{Lemma}
\newtheorem{Example}[Proposition]{Example}
\newtheorem{Remark}[Proposition]{Remark}

\maketitle

%\addtocounter{section}{-1}

\begin {abstract}  All quasi-affine  connected Generalized Dynkin Diagrams with rank $> 5$ are found.  All quasi-affine  Nichols (Lie braided) algebras with rank $> 5$  are also found.
\vskip.2in
\noindent {\em 2010 Mathematics Subject Classification}: 16W30,  16G10  \\
{\em Keywords}:  Quasi-affine,   Nichols  algebra,   Generalized Dynkin Diagram, Arithmetic GDD.
\emph{}
\end {abstract}
 %\tableofcontents
\section {Introduction and Preliminaries}\label {s0}

Nichols algebras play a fundamental role in the classification of finite-dimensional complex pointed Hopf algebras by means of the lifting method developed by Andruskiewitsch and Schneider \cite {AS02, AS10, AHS08}.
Heckenberger \cite {He06a, He05} classified arithmetic root systems.   Heckenberger \cite {He06b} proved a GDD is  arithmetic  if and only if corresponding matrix is a finite Cartan matrix for GDDs of  Cartan types. W. Wu,   S. Zhang and   Y.-Z. Zhang \cite {WZZ15b} proved that a Nichols Lie braided algebra is a finite dimensional  if and only if its GDD, which   fixed
parameter is of finite order, is  arithmetic.

%\section*{ Preliminaries}\label {s1}

We now recall some basic concepts of the graph theory (see \cite {Ha}).
Let $\Gamma _1$ be a non-empty set and $\Gamma _2 \subseteq \{  \{ u, v\} \mid u, v \in \Gamma _1, \hbox { with } u \not= v \} \subseteq 2 ^{\Gamma_1}.$ Then $\Gamma = (\Gamma _1, \Gamma_2)$ is called a graph;  $\Gamma_1$ is called the vertex set of  $\Gamma$; $\Gamma_2$ is called the edge set of  $\Gamma$; Element $\{u, v\} \in \Gamma_2$ is called an edge. Let $F$ be an algebraically closed field of characteristic zero and $F^* := \{x \mid x\in F, x \not= 0 \}$.  If   $\{x_1,   \cdots,   x_n\}$ is  a basis of  vector space  $V$ and
$C(x_i\otimes  x_j) = q_{ij} x_j\otimes x_i$ with $q_{ij} \in F^*$,
then $V$  is called a braided vector space of diagonal type, $\{x_1,   \cdots,   x_n\}$
is called  canonical basis and $(q_{ij})_{n\times n}$  is called braided matrix. Let
$\widetilde{q}_{i j}:= q_{ij}q_{ji}$ for  $i, j \in \{1, 2, \cdots, n \}$ with $i \not=j$. If let $\Gamma _1= \{1, 2, \cdots, n \}$ and $\Gamma _2 = \{  \{ u, v\} \mid \widetilde{q}_{u v} \not= 1,   u, v \in \Gamma _1, \hbox { with } u \not= v \}$, then $\Gamma = (\Gamma _1, \Gamma_2)$ is a graph. Set $q_{ii}$ over vertex $i$ and $\widetilde{q}_{i j}$ over edge $\{i, j\}$ for $i, j \in \Gamma _1$ with $i \not=j$.
Then $\Gamma = (\Gamma _1, \Gamma_2)$ is called a Generalized Dynkin Diagram of braided vector space $V,$ written as GDD in short(see \cite [Def. 1.2.1] {He05}).
If $\Delta (\mathfrak{B}(V))$ is an arithmetic root system.   then we call its GDD  an arithmetic GDD  for convenience.

Let $(q_{ij})_{n \times n }$ be a braided  matrix. If  $q_{ij} q_{ji}  \left \{  \begin{array}{ll}
 \not= 1,  & \mbox {when }  \mid j-i \mid = 1\\
  =1,   & \mbox {when}  \mid j-i \mid \not= 1    \\
\end{array}\right. $ for any $1\le i\not= j \le n,$ then  $(q_{ij})_{n \times n }$ is called a  chain or labelled  chain.
If $(q_{ij})_{n \times n }$ is  a  chain and \begin {eqnarray} \label {ppe1}(q_{11} q_{1, 2} q_{2, 1} -1)(q_{11} +1)=0;
 (q_{n, n } q_{n, n-1} q_{n-1, n} -1)(q_{n,n} +1)=0;\end {eqnarray} i.e.
  \begin {eqnarray} \label {ppe2}
   q_{ii} +1= q_{i, i-1} q_{i-1, i}q_{i, i+1} q_{i+1, i}-1=0
   \end  {eqnarray}
    \begin {eqnarray} \label {ppe3}
   \mbox { or  }
  q_{ii}q_{i, i-1} q_{i-1, i}=q_{ii} q_{i, i+1} q_{i+1, i}=1,\end  {eqnarray}  $1<i < n$,
   then  the braided matrix's GDD is called a simple chain (see \cite [Def.1]{He06a}).  Conditions (\ref {ppe1}), (\ref {ppe2}) and (\ref {ppe3}) are called simple chain conditions.  Let
$q:= q_{n, n} ^2 q_{n -1, n}q_{n, n-1}$.
%%If $(q_{ij})_{n \times n }$ is a simple chain and  $q_{i, i-1} q_{i-1, i} %%=q^{-1}$, $1< i \le n$, then  $(q_{ij})_{n \times n }$ is called a pure %%simple chain, otherwise, it is called a mixed simple chain.

 For $0\le j \le n$ and $0<i_1 < i_2 < \cdots < i_j \le n,$ let $C_{n,q, i_1, i_2, \cdots, i_j }$ denote  a simple chain which satisfies the condition: $\widetilde{q}_{i, i-1}=q$ if and only if  $i \in \{i_1, i_2, \cdots, i_j\}$,
where
$\widetilde{q} _{1, 0}:= \frac {1}{q_{11}^2\widetilde{q}_{12}}$. For example,
 $n \ge 2$,

$\begin{picture}(100,      15)
\put(27,      1){\makebox(0,     0)[t]{$\bullet$}}
\put(60,      1){\makebox(0,      0)[t]{$\bullet$}}
\put(93,     1){\makebox(0,     0)[t]{$\bullet$}}
\put(159,      1){\makebox(0,      0)[t]{$\bullet$}}
\put(192,     1){\makebox(0,      0)[t]{$\bullet$}}
\put(28,      -1){\line(1,      0){30}}
\put(61,      -1){\line(1,      0){30}}
\put(130,     1){\makebox(0,     0)[t]{$\cdots\cdots\cdots\cdots$}}
\put(160,     -1){\line(1,      0){30}}
\put(22,     -15){1}
\put(58,      -15){2}
\put(91,      -15){3}
\put(157,      -15){$n-1$}
\put(191,      -15){$n$}
\put(22,     10){$q$}
\put(58,      10){$q$}
\put(91,      10){$q$}
\put(157,      10){$q$}
\put(191,      10){$q$}
\put(40,      5){$q^{-1}$}
\put(73,      5){$q^{-1}$}
\put(172,     5){$q^{-1}$}
\put(210,        -1)  {$, q \in F^{*}/\{1, -1\}$. }
\end{picture}$\\

is  $C_{n,q, i_1, i_2, \cdots, i_j }$ with $j=0$.\\ \\

$\begin{picture}(100,      15)
\put(27,      1){\makebox(0,     0)[t]{$\bullet$}}
\put(60,      1){\makebox(0,      0)[t]{$\bullet$}}
\put(93,     1){\makebox(0,     0)[t]{$\bullet$}}
\put(159,      1){\makebox(0,      0)[t]{$\bullet$}}
\put(192,     1){\makebox(0,      0)[t]{$\bullet$}}
\put(28,      -1){\line(1,      0){30}}
\put(61,      -1){\line(1,      0){30}}
\put(130,     1){\makebox(0,     0)[t]{$\cdots\cdots\cdots\cdots$}}
\put(160,     -1){\line(1,      0){30}}
\put(22,     -15){1}
\put(58,      -15){2}
\put(91,      -15){3}
\put(157,      -15){$n-1$}
\put(191,      -15){$n$}
\put(22,     10){$q^{-1}$}
\put(58,      10){$q^{-1}$}
\put(91,      10){$q^{-1}$}
\put(157,      10){$q^{-1}$}
\put(191,      10){$-1$}
\put(40,      5){$q^{}$}
\put(73,      5){$q^{}$}
\put(172,     5){${q}$}
\put(210,        -1)  {$, q \in F^{*}/\{1\}$. }
\end{picture}$\\

is  $C_{n,q, i_1, i_2, \cdots, i_j }$ with $j=n$.

For the convenience, we let
 $C_{1,q, i_1, i_2, \cdots, i_j }$  denote the GDD with length $1$ satisfied the following conditions:   $q= q_{11}$ when $q_{11} \not=-1$; $q$ can be any non-zero number   when $q_{11} =-1$.

Every connected subGDD of
every   arithmetic GDD in Row 1-10 in Table C is  called a classical GDD.  Every connected     GDD which is not classical is called an exception.

If omitting every vertex in a connected GDD with rank $n>1$  is an arithmetic GDD and this GDD is not an arithmetic  GDD, then this GDD is called  a quasi-affine  GDD over  connected subGDD of this GDD with rank $n-1$, quasi-affine  GDD in short. In this case, Nichols algebra $\mathfrak B(V)$
and Nichols Lie braded algebra $\mathfrak L(V)$ are said to be quasi-affine.
In other word, if a GDD is quasi-affine  of a braided vector space $V$ which   fixed
parameter is of finite order, then Nichols algebra and Nichols Lie braded algebra
of every proper subGDD are finite dimensional with $\dim   \mathfrak B { (V)} = \infty$
and $\dim   \mathfrak L { (V)} = \infty$.

In this paper,  using Table $A1$, $A2$, $B$ and $C$ in \cite {He06a, He05}, we find all quasi-affine  connected Generalized Dynkin Diagrams with rank $> 5$. We also  find all quasi-affine   Nichols  algebras and quasi-affine  Nichols Lie braided algebras with rank $> 5$.

In this paper we use following notations without special announcement. every braided vector  $V$  is of diagonal type with dimension $n > 1$; is connected with fixed parameter $q \in R_3$. Let  $F$ be an  algebraic closed base field of characteristic zero. $\mathbb Z =: \{x \mid x \hbox { is an integer }\}.$
$\mathbb N_0 =: \{x \mid x \in \mathbb Z,    x\ge 0\}. $

\section {Properties about arithmetic GDDs  } \label {s2}

\begin {Lemma} \label {mainlemma} A GDD is a classical GDD if and only if it is one of classical Type 1-7, Here
classical types are  listed as follows:\\

 {\ }\ \ \ \ \ \ \ \ \ \ \ \ \ \ \ \ \ \ \ \ \ \ \ \   \ \ \ \ \ \  $\begin{picture}(100,      15)

\put(-125,      -1){ {\rm   Type   1},  $2\le n$.}
\put(80,      1){\makebox(0,      0)[t]{$\bullet$}}

\put(48,      -1){\line(1,      0){33}}
\put(10,      1){\makebox(0,     0)[t]{$C_{n-1, q, i_1, i_2, \cdots, i_j }$}}

\put(-18,     10){$$}
\put(0,      5){$$}
\put(22,     10){$$}
\put(50,      5){$q^{-2}$}

\put(68,      10){$q^2$}

 \ \ \ \ \ \ \ \ \ \ \ \ \ \ \ \ \ \ \
  \ \ \ \ \ \ \ \ \ \ \ \ \ \ \ \ \ \ \ { }$q \in F^{*}\setminus \{1, -1\}$, $0\le j\le n-1.$

 \put(220,      -1) {}
\end{picture}$\\ \\

 {\ }\ \ \ \ \ \ \ \ \ \ \ \ \ \ \ \ \ \ \ \ \ \ \ \   \ \ \ \ \ \  $\begin{picture}(100,      15)

\put(-125,      -1){ {\rm   Type   2}. $2\le n$.}

\put(90,      1){\makebox(0,      0)[t]{$\bullet$}}

\put(58,      -1){\line(1,      0){33}}
\put(27,      1){\makebox(0,     0)[t]{$C_{n-1, q^2, i_1, i_2, \cdots, i_j }$}}

\put(-8,     10){$$}
\put(0,      5){$$}
\put(22,     10){$$}
\put(70,      5){$q^{-2}$}

\put(88,      10){$q$}

  \ \ \ \ \ \ \ \ \ \ \ \ \ \ \ \ \ \ \
  \ \ \ \ \ \ \ \ \ \ \ \ \ \ \ \ \ \ \ $q \in F^{*}\setminus \{1, -1\}$.  $0\le j\le n-1.$
\end{picture}$\\

 {\ }\ \ \ \ \ \ \ \ \ \ \ \ \ \ \ \ \ \ \ \ \ \ \ \   \ \ \ \ \ \  $\begin{picture}(100,      15)

\put(-125,      -1){ {\rm   Type   3}, $2\le n$.}

\put(80,      1){\makebox(0,      0)[t]{$\bullet$}}

\put(48,      -1){\line(1,      0){33}}
\put(27,      1){\makebox(0,     0)[t]{$C_{n-1, q^{-2}, i_1, i_2, \cdots, i_j }$}}

\put(-18,     10){$$}
\put(0,      5){$$}
\put(22,     10){$$}
\put(60,      5){$q^{2}$}

\put(78,      10){$-q^{-1}$}

  \ \ \ \ \ \ \ \ \ \ \ \ \ \ \ \ \ \ \
  \ \ \ \ \ \ \ \ \ \ \ \ \ \ \ \ \ \ \ $q \in F^{*}\setminus \{1, -1\}$, $0\le j\le n-1.$

\end{picture}$\\ \\

 {\ }\ \ \ \ \ \ \ \ \ \ \ \ \ \ \ \ \ \ \ \ \ \ \ \   \ \ \ \ \ \  $\begin{picture}(100,      15)

\put(-125,      -1){{\rm   Type   4}. $2\le n$. }

\put(80,      1){\makebox(0,      0)[t]{$\bullet$}}

\put(48,      -1){\line(1,      0){33}}
\put(27,      1){\makebox(0,     0)[t]{$C_{n-1, -q^{-1}, i_1, i_2, \cdots, i_j }$}}

\put(-18,     10){$$}
\put(0,      5){$$}
\put(22,     10){$$}
\put(60,      5){$-q$}

\put(78,      10){$q$}

  \ \ \ \ \ \ \ \ \ \ \ \ \ \ \ \ \ \ \
  \ \ \ \ \ \ \ \ \ \ \ \ \ \ \ \ \ \ \ $q^3 =1$. $0\le j\le n-1.$

\end{picture}$\\
\\

\ \ \ \     \ \ \ \ \ \  $\begin{picture}(100,      15)

\put(-45,      -1){ {\rm   Type   5}. $3\le n$.}
\put(124,      1){\makebox(0,      0)[t]{$C_{n-2, q^{}, i_1, i_2, \cdots, i_j }$}}
\put(190,     -11){\makebox(0,     0)[t]{$\bullet$}}
\put(190,    15){\makebox(0,     0)[t]{$\bullet$}}
\put(162,    -1){\line(2,     1){27}}
\put(190,      -14){\line(-2,     1){27}}

\put(120,      10){$$}

\put(135,      5){$$}

\put(155,     10){$$}

\put(160,     -20){$q^{-1}$}
\put(165,      15){$q^{-1}$}

\put(193,      -12){$q$}
\put(193,      18){$q$}

 \ \ \ \ \ \ \ \ \ \ \ \ \ \ \ \ \ \ \ \ \ \ \ \ \ \ \ \ \ \ \ \ \ \ \ \ \ \
  \ \ \ \ \ \ \ \ \ \ \ \ \ \ \ \ \ \ \ $q\not= 1,$    $0\le j\le n-2.$

\put(215,        -1)  { }
\end{picture}$\\ \\

\ \ \ \      \ \ \ \ \ \  $\begin{picture}(100,      15)

\put(-45,      -1){{\rm   Type   6}. $3\le n$. }

\put(104,      1){\makebox(0,      0)[t]{$C_{n-2, q^{}, i_1, i_2, \cdots, i_j }$}}
\put(170,     -11){\makebox(0,     0)[t]{$\bullet$}}
\put(170,    15){\makebox(0,     0)[t]{$\bullet$}}

\put(142,    -1){\line(2,     1){27}}
\put(170,      -14){\line(-2,     1){27}}

\put(170,      -14){\line(0,     1){27}}

\put(100,      10){$$}
\put(120,      5){$$}

\put(127,     10){$$}

\put(140,     -20){$q^{-1}$}
\put(145,      15){$q^{-1}$}

\put(178,      -20){$-1$}

\put(178,      0){$q^{2}$}

\put(175,      10){$-1$}

 \ \ \ \ \ \ \ \ \ \ \ \ \ \ \ \ \ \ \ \ \ \ \ \ \ \ \ \ \ \ \ \ \ \ \ \ \ \
  \ \ \ \ \ \ \ \ \ \ \ \ \ \ \ \ \ \ \  $q^2 \not= 1$. $0\le j\le n-2.      $
\put(195,        -1)  {,  }
\end{picture}$\\

{\ }\!\!\!\!\!\!
{\rm   Type   7}. $1\le n$.
$C_{n, q^{-1}, i_1, i_2, \cdots, i_j }.$
$q \not= 1$,  $0\le j\le n$.

\end {Lemma}
{\bf Proof.} We only consider the case $n>4$ since every  GDD classical  type   with $n<5$ is subGDD of a GDD classical type with  $n> 4.$

Necessity is obviously. Now we show the sufficiency.

 Type   1.   Type   1 is classical  by  GDD $1$ of Row $10$ in Table C when  $2 \le j \le n-1$.   Type   1 is classical  by  GDD $1$ of Row $9$ in Table C when  $j=0$.  Type   1 is classical  by GDD $2$ of Row $9$  in Table C and by GDD $1$ of Row $7$  in Table C
 when  $j=1$.

Type 2.  Type 2 is classical  by  GDD $1$ of Row $4$ in Table C when $q^4\not=1$ and   $1 \le j \le n-1$. Type 2 is classical  by  GDD $1$ of Row $3$ in Table C when $q^2\not=1$ and $j=0$. Type 2 is classical  by  GDD $1$ of Row $3$ in Table C when $1 \le j \le n-1.$

Type 3.  Type 3 is classical  by  GDD $2$ of Row $4$ in Table C when $q^4\not=1$ and  $1 \le j \le n-1$.

Type 4.  Type 4 is classical  by  GDD $1$, $2$ of Row $6$ in Table C when  $1 \le j \le n-1$.  Type 4 is classical  by  GDD $1$ of Row $5$ in Table C when  $j=0$.

Type 5.  Type 5 is classical  by  GDD $3$ of Row $10$ in Table C when  $1 \le j \le n-2$ with $q \not=1$.  Type 5 is classical  by  GDD $1$ of Row $8$  in Table C when  $j=0$ with $q^2\not=1$.

Type 6.  Type 6 is classical  by  GDD $2$ of Row $10$ in Table C when  $1 \le j \le n-2$.  Type 1 is classical  by  GDD $3$ of Row $9$ in Table C when  $j=0$. $\Box$

The subGDD of the right hand of  Type i in above Lemma is called head type i, written as h-type i with $i=1,2,\cdots, 6.$ h-type 7 has two forms. Right end is $-1$:    $\begin{picture}(100,      15)  \put(0,      -1){ }

\put(60,      1){\makebox(0,      0)[t]{$\bullet$}}

\put(28,      -1){\line(1,      0){33}}
%\put(27,      1){\makebox(0,     0)[t]{$\bullet$}}

%\put(-14,      1){\makebox(0,     0)[t]{$\bullet$}}

%\put(-14,     -1){\line(1,      0){50}}

\put(22,     10){}
\put(38,      5){$q^{-1}$}

\put(60,      10){$-1^{}$}

  \ \ \ \ \ \ \ \ \ \ \ \ \ \ \ \ \ \ \ {  }
\end{picture}$\\
 written $T6$; right end is $-1$:    $\begin{picture}(100,      15)  \put(-20,      -1){ }

\put(60,      1){\makebox(0,      0)[t]{$\bullet$}}

\put(28,      -1){\line(1,      0){33}}
%\put(27,      1){\makebox(0,     0)[t]{$\bullet$}}

%\put(-14,      1){\makebox(0,     0)[t]{$\bullet$}}

%\put(-14,     -1){\line(1,      0){50}}

\put(22,     10){}
\put(38,      5){$q^{-1}$}

\put(60,      10){$q^{}$}

  \ \ \ \ \ \ \ \ \ \ \ \ \ \ \ \ \ \ \ { written  T5. }
\end{picture}$\\ The left end of every GDD  in above  Lemma is called the tail of the GDD.

{Remark:}
If $q\in R_3$, then  Type   1 turn into the following:\\

  \ \ \ \ \ \  $\begin{picture}(100,      15)

\put(80,      1){\makebox(0,      0)[t]{$\bullet$}}

\put(48,      -1){\line(1,      0){33}}
\put(27,      1){\makebox(0,     0)[t]{$C_{n-1, q^{}, i_1, i_2, \cdots, i_j }$}}

\put(-18,     10){$$}
\put(0,      5){$$}
\put(22,     10){$$}
\put(60,      5){$q$}

\put(78,      10){$q^{-1}$}

 \ \ \ \ \ \ \ \  \ \ \ \ \ \ \ \ \ \ \ \ \ \ \ \ \ \ \ {$,q^3 =1$.  }

\end{picture}$\\

If $q\in R_3$, then Type   3 turn into the following:\\

  \ \ \ \ \ \  $\begin{picture}(100,      15)

\put(80,      1){\makebox(0,      0)[t]{$\bullet$}}

\put(48,      -1){\line(1,      0){33}}
\put(27,      1){\makebox(0,     0)[t]{$C_{n-1, q^{-1}, i_1, i_2, \cdots, i_j }$}}

\put(-18,     10){$$}
\put(0,      5){$$}
\put(22,     10){$$}
\put(60,      5){$q$}

\put(78,      10){$q$}

 \ \ \ \ \ \ \ \  \ \ \ \ \ \ \ \ \ \ \ \ \ \ \ \ \ \ \ {$,q^3 =1$.  }

\end{picture}$\\

If set   $\xi = -q^{-1}$,
then  Type   3 turn into Type 2: \\

   \ \ \ \ \ \  $\begin{picture}(100,      15)

\put(90,      1){\makebox(0,      0)[t]{$\bullet$}}

\put(58,      -1){\line(1,      0){33}}
\put(27,      1){\makebox(0,     0)[t]{$C_{n-1, \xi ^2, i_1, i_2, \cdots, i_j }$}}

\put(-8,     10){$$}
\put(0,      5){$$}
\put(22,     10){$$}
\put(70,      5){$\xi ^{-2}$}

\put(88,      10){$\xi $}

  \ \ \ \ \ \ \ \ \ \ \ \ \ \ \ \ \ \ \   \ \ \ \ \ \ \ \ \ \ \ \ \ \ \ \ \ \ \ {}$\xi  \in F^{*}\setminus \{1, -1\}$.  $0\le j\le n-1.$
\end{picture}$\\ \\
That is, Type 2 and Type 3 are the same.

If $q\in R_3$, then  Type   3 turn into the following: \\ \\
  \ \ \ \ \ \  $\begin{picture}(100,      15)

\put(80,      1){\makebox(0,      0)[t]{$\bullet$}}

\put(48,      -1){\line(1,      0){33}}
\put(27,      1){\makebox(0,     0)[t]{$C_{n-1, q^{-1}, i_1, i_2, \cdots, i_j }$}}

\put(-18,     10){$$}
\put(0,      5){$$}
\put(22,     10){$$}
\put(60,      5){$q$}

\put(78,      10){$-q$}

 \ \ \ \ \ \ \ \  \ \ \ \ \ \ \ \ \ \ \ \ \ \ \ \ \ \ \ {$,q^3 =1$.  $\Box$ }

\end{picture}$ \\

A GDD is called a simple cycle  if  it is a cycle and omitting every vertex in the  GDD
is a simple chain. \\

{\ } \ \ \ \ \ \
  \ \ \ \ \ \
  \ \ \ \ \ \
 $\begin{picture}(100,      15)

\put(-98,     10){GDD:}
\put(-28,     10){$\hbox {h-Type s} + $}
\put(0,      5){$$}
\put(22,     10){ $C_{u, \xi^{}, i_1, i_2, \cdots, i_j }$}
\put(40,      5){$$}

\put(98,      10){$+ \hbox {h-Type t}$}

\put(98,      10){}

\end{picture}$\\
 is called a bi-classical GDD, written as Type s-Type t,  where fixed parameters are  matched as  in
Lemma \ref {mainlemma} with $1\le s, t \le 6.$

A vertex is called an end vertex if its degree  is 1, i.e. there exists only one vertex which connects it.

\begin {Definition} \label {6.1.-2}
{ {\rm (i)}} A connected arithmetic chain   is called a quasi-classical GDD  if omitting an end vertex  is a classical GDD which tail is other end vertex. Further more, the other end vertex is called the tail of the quasi-classical GDD.

{\rm (ii)} A connected  arithmetic  non-chain   is called a quasi-classical GDD if the following conditions hold:  omitting a vertex is a connected classical GDD; omitting an other vertex is also a connected classical GDD; the two tails of the two classical GDDs  are the same. Further more, the two tails of the two classical GDDs is called the tail of the quasi-classical GDD.

{ {\rm (iii)}} A connected arithmetic chain   is called a semi-classical GDD  if omitting an end vertex  is a simple  chain.

{\rm (iv)} A connected  arithmetic  non-chain   is called a semi-classical GDD if the following conditions hold:  omitting a vertex is a connected simple chain; omitting an other vertex is also a connected simple chain.

{\rm (v)} If adding T5 on tail of a  quasi-classical  GDD is an arithmetic GDD   then the  quasi-classical GDD is called continual on the tail via T5.  Otherwise  adding T5 on tail of a  quasi-classical  GDD is called a quasi-affine  continual GDD on the tail via T5.  For T6  we can similarly define these.

%%{\rm (vi)} If adding a diagram on a vertex  of an  arithmetic  GDD is not an arithmetic GDD   then the  new GDD is called non continue on this %%vertex.

%%{\rm (vii)} If omitting  a tail of semi-classical  GDD is connected,    then the new   GDD is called its  pre-GDD.
\end {Definition}

%Remark: In {\rm (ii)} in above definition. two classical GDDs which is obtained  by omitting two vertexes respectively have the same tails, %which is the tail of the quasi-classical GDD.

Remark: Every semi-classical GDD is  a quasi-classical GDD.
 Every  classical GDD is  a semi-classical GDD.

A quasi-classical GDD which is not  semi-classical is called a strict quasi-classical GDD. A semi-classical GDD which is not classical is called a strict semi-classical GDD.

Omitting the longest  simple chain which contains the tail of a quasi-affine GDD is called the head of the GDD.

\begin {Definition} \label {6.1.-3}
%{ {\rm (i)}}   T5 adding   on vertex i of  GDD  is called quasi-affine  continue  of GDD on i
% if GDD is not continue via T5  on i; T6 adding   on vertex i of  GDD  is called quasi-affine  continue of GDD on i
%if GDD is not continue via T6  on i.

{\rm (i)} h-Type i  adding on tail of  semi-classical GDD is called classical + semi-classical, $1\le i \le 4.$
h-Type 5  adding on tail of  semi-classical GDD is called classical + semi-classical  when  semi-classical GDD is continual via T5 ;  h-Type 6  adding on tail of  semi-classical GDD is called classical + semi-classical when  semi-classical GDD is continual via T6.

{\rm (ii)}
The tail of  a semi-classical GDD adding on tail of a semi-classical GDD is called bi-semi-classical, if omitting every vertex of the new  GDD is arithmetic.

 \end {Definition}

 Remark:  Classical + semi-classical can be obtained as follows.
  Adding a vertex on vertex connected tail  such that its tail contained in h-Type 5  if the tail of a semi-classical GDD is not $-1$.
   Adding a vertex on vertex connected tail  such that its tail contained in h- Type 6  if the tail of a semi-classical GDD is  $-1$.

\begin {Lemma} \label {2.2}Assume rank $n >4.$

{\rm (i)} A  GDD is not  an arithmetic GDD  if the GDD contains subGDDs:\\ \\

 $\ \ \ \ \  \ \   
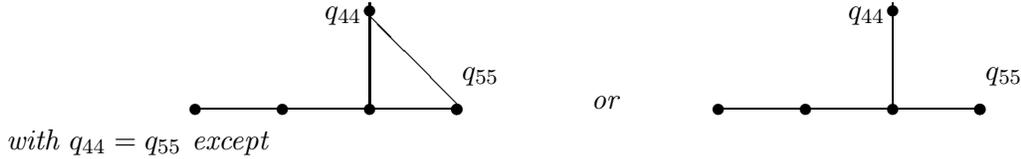
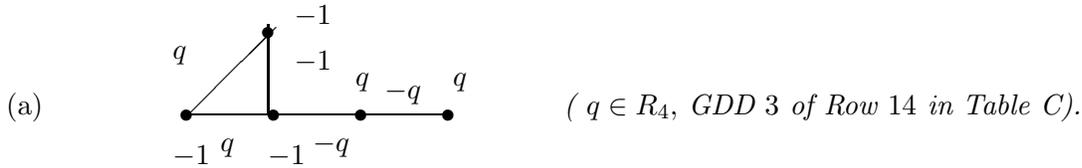
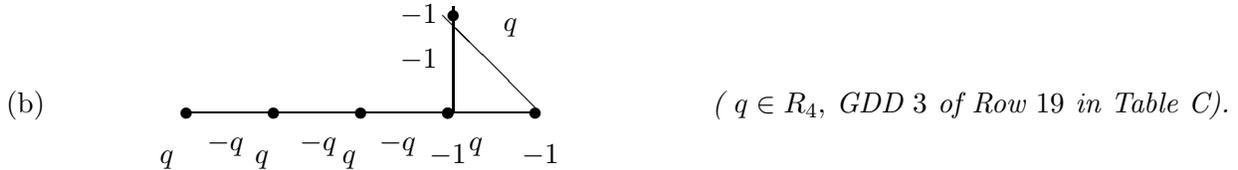
\begin{picture}(100,      15)  \put(-45,      -1){ }
\put(27,      1){\makebox(0,     0)[t]{$\bullet$}}
\put(60,      1){\makebox(0,      0)[t]{$\bullet$}}
\put(93,     1){\makebox(0,     0)[t]{$\bullet$}}

\put(126,      1){\makebox(0,    0)[t]{$\bullet$}}

\put(126,      1){\makebox(0,    0)[t]{$\bullet$}}
\put(28,      -1){\line(1,      0){33}}
\put(61,      -1){\line(1,      0){30}}
\put(94,     -1){\line(1,      0){30}}

\put(22,     10){$$}
\put(40,      5){$$}
\put(58,      10){$$}
\put(70,      5){$$}

\put(91,      10){$$}

\put(102,     5){$$}

\put(128,      10){$q_{55}$}

\put(93,    38){\makebox(0,     0)[t]{$\bullet$}}

\put(93,      -1){\line(0,     1){40}}

\put(76,      33){$q_{44}$}

\put(73,      16){$$}

\put(126,      1){\line(-1,     1){35}}

\put(114,      30){$$}

\put(166,        -1)  {  }
\end{picture}$\ \ \ \ \ \ \ \ \ \ \ \ \ \ \ \ \ \
 $\ \ \ \ \  \ \   \begin{picture}(100,      15)  \put(-45,      -1){ }

\put(-25,      -1){ or}
\put(27,      1){\makebox(0,     0)[t]{$\bullet$}}
\put(60,      1){\makebox(0,      0)[t]{$\bullet$}}
\put(93,     1){\makebox(0,     0)[t]{$\bullet$}}

\put(126,      1){\makebox(0,    0)[t]{$\bullet$}}

\put(126,      1){\makebox(0,    0)[t]{$\bullet$}}
\put(28,      -1){\line(1,      0){33}}
\put(61,      -1){\line(1,      0){30}}
\put(94,     -1){\line(1,      0){30}}

\put(22,     10){$$}
\put(40,      5){$$}
\put(58,      10){$$}
\put(70,      5){$$}

\put(91,      10){$$}

\put(102,     5){$$}

\put(128,      10){$q_{55}$}

\put(93,    38){\makebox(0,     0)[t]{$\bullet$}}

\put(93,      -1){\line(0,     1){40}}

\put(76,      33){$q_{44}$}

\put(73,      16){$$}

\put(166,        -1)  { }
\end{picture}$\\
with  $q_{44} =q_{55}$
except\\ \\ \\

 $\ \ \ \ \  \ \   \begin{picture}(100,      15)  \put(-45,      -1){ }

\put(-45,      -1){ {\rm (a)}}
\put(27,      1){\makebox(0,     0)[t]{$\bullet$}}
\put(60,      1){\makebox(0,      0)[t]{$\bullet$}}
\put(93,     1){\makebox(0,     0)[t]{$\bullet$}}

\put(126,      1){\makebox(0,    0)[t]{$\bullet$}}

\put(28,      -1){\line(1,      0){33}}
\put(61,      -1){\line(1,      0){30}}
\put(94,     -1){\line(1,      0){30}}

\put(22,     -20){$-1$}
\put(40,      -15){$q$}
\put(58,      -20){$-1 ^{}$}
\put(75,      -15){$-q ^{}$}

\put(91,      10){$q ^{}$}

\put(102,     5){$-q ^{}$}

\put(128,      10){$q ^{}$}

\put(58,    32){\makebox(0,     0)[t]{$\bullet$}}

\put(58,      -1){\line(0,     1){34}}

\put(22,    20){$q$}

\put(68,      33){$-1$}

\put(68,      16){$-1 ^{}$}

\put(28,      -1){\line(1,      1){33}}

\put(166,        -1)  { ( $q \in R_4,$ GDD $3$ of Row $14$  in Table C). }
\end{picture}$\\ \\ \\

 $\ \ \ \ \  \ \   \begin{picture}(100,      15)  \put(-45,      -1){ } \put(-45,      -1){ {\rm (b)}}
\put(27,      1){\makebox(0,     0)[t]{$\bullet$}}
\put(60,      1){\makebox(0,      0)[t]{$\bullet$}}
\put(93,     1){\makebox(0,     0)[t]{$\bullet$}}

\put(126,      1){\makebox(0,    0)[t]{$\bullet$}}

\put(159,      1){\makebox(0,    0)[t]{$\bullet$}}

\put(126,      1){\makebox(0,    0)[t]{$\bullet$}}
\put(28,      -1){\line(1,      0){33}}
\put(61,      -1){\line(1,      0){30}}
\put(94,     -1){\line(1,      0){30}}
\put(126,      -1){\line(1,      0){33}}

\put(17,    - 20){$q$}
\put(35,     - 15){$-q ^{}$}
\put(53,     - 20){$q$}
\put(70,     - 15){$-q$}

\put(86,      -20){$q$}
\put(100,      -15){$-q$}
\put(119,      -20){$-1$}
\put(134,    - 15){$q ^{}$}
\put(154,    - 20){${-1}$}

%\put(75,    18){\makebox(0,     0)[t]{$\bullet$}}

\put(60,     -1){\line(1,      0){30}}

%\put(91,      0){\line(-1,     1){17}}
%\put(60,     -1){\line(1,      1){17}}

\put(128,    38){\makebox(0,     0)[t]{$\bullet$}}

\put(128,      -1){\line(0,     1){40}}

\put(108,      33){$-1$}

\put(108,      16){${-1}$}

\put(159,      1){\line(-1,     1){35}}

\put(147,      30){$q ^{}$}

\put(222,        -1)  { ( $q \in R_4,$ GDD $3$ of Row $19$  in Table C). }
\end{picture}$\\ \\

{\rm (c)} Classical GDDs.

{\rm (ii)}  If a GDD contains subGDD
  \\ \\

{\ }

{\ }

 $\ \ \ \ \  \ \   
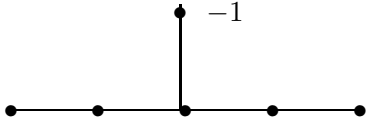
\begin{picture}(100,      15)  \put(-45,      -1){}

\put(-6,      1){\makebox(0,     0)[t]{$\bullet$}}

\put(27,      1){\makebox(0,     0)[t]{$\bullet$}}
\put(60,      1){\makebox(0,      0)[t]{$\bullet$}}
\put(93,     1){\makebox(0,     0)[t]{$\bullet$}}

\put(126,      1){\makebox(0,    0)[t]{$\bullet$}}

%\put(159,      1){\makebox(0,    0)[t]{$\bullet$}}

%\put(126,      1){\makebox(0,    0)[t]{$\bullet$}}

\put(-6,      -1){\line(1,      0){33}}
\put(28,      -1){\line(1,      0){33}}
\put(61,      -1){\line(1,      0){30}}
\put(94,     -1){\line(1,      0){30}}
%\put(126,      -1){\line(1,      0){33}}

\put(22,     10){$$}
\put(40,      5){$$}
\put(58,      -20){$$}
\put(75,      -15){$$}

\put(91,      10){$$}

\put(102,     5){$$}

\put(128,      10){$$}

\put(135,     5){$$}

\put(161,      10){$$}

\put(58,    38){\makebox(0,     0)[t]{$\bullet$}}

\put(58,      -1){\line(0,     1){40}}

\put(68,      33){${-1}$}

\put(68,      16){$$}

\put(166,        -1)  {  }
\end{picture}$\\ \\
then the GDD   is not an arithmetic GDD  when the  fix  parameter $q \not= -1.$

\end {Lemma}

\begin {Lemma} \label {2.6} Assume rank $n >4.$
Chain  GDD is not an arithmetic  if there exist two places where simple chain conditions do not hold except\\

 $\ \ \ \ \  \ \   \begin{picture}(100,      15)  \put(-45,      -1){\rm (a)}
\put(27,      1){\makebox(0,     0)[t]{$\bullet$}}
\put(60,      1){\makebox(0,      0)[t]{$\bullet$}}
\put(93,     1){\makebox(0,     0)[t]{$\bullet$}}

\put(126,      1){\makebox(0,    0)[t]{$\bullet$}}

\put(159,      1){\makebox(0,    0)[t]{$\bullet$}}

\put(126,      1){\makebox(0,    0)[t]{$\bullet$}}
\put(28,      -1){\line(1,      0){33}}
\put(61,      -1){\line(1,      0){30}}
\put(94,     -1){\line(1,      0){30}}
\put(126,      -1){\line(1,      0){33}}

\put(22,     10){$q^{}$}
\put(40,      5){$q^{-1}$}
\put(58,      10){$q^{}$}
\put(70,      5){$q^{-1}$}

\put(88,      10){$-1$}

\put(102,     5){${-1}$}

\put(120,      10){$-1^{}$}

\put(135,     5){$q^{-1}$}

\put(161,      10){${q}$}

\put(186,        -1)  { $q \in R_4,$ GDD $1$ of Row $14$  in Table C. }
\end{picture}$\\

 $\ \ \ \ \  \ \   \begin{picture}(100,      15)  \put(-45,      -1){\rm (b) }
\put(27,      1){\makebox(0,     0)[t]{$\bullet$}}
\put(60,      1){\makebox(0,      0)[t]{$\bullet$}}
\put(93,     1){\makebox(0,     0)[t]{$\bullet$}}

\put(126,      1){\makebox(0,    0)[t]{$\bullet$}}

\put(159,      1){\makebox(0,    0)[t]{$\bullet$}}

\put(126,      1){\makebox(0,    0)[t]{$\bullet$}}
\put(28,      -1){\line(1,      0){33}}
\put(61,      -1){\line(1,      0){30}}
\put(94,     -1){\line(1,      0){30}}
\put(126,      -1){\line(1,      0){33}}

\put(22,     10){$q^{2}$}
\put(40,      5){$q^{-2}$}
\put(58,      10){$q^{2}$}
\put(70,      5){$q^{-2}$}

\put(87,      10){$-1$}

\put(102,     5){$q^{-1}$}

\put(128,      10){$q^{}$}

\put(135,     5){$q^{-2}$}

\put(161,      10){$q^{2}$}

\put(186,        -1)  {  $q \in R_5,$ GDD $1$ of Row $15$  in Table C. }
\end{picture}$\\

 $\ \ \ \ \  \ \   \begin{picture}(100,      15)  \put(-45,      -1){\rm (c) }
\put(27,      1){\makebox(0,     0)[t]{$\bullet$}}
\put(60,      1){\makebox(0,      0)[t]{$\bullet$}}
\put(93,     1){\makebox(0,     0)[t]{$\bullet$}}

\put(126,      1){\makebox(0,    0)[t]{$\bullet$}}

\put(159,      1){\makebox(0,    0)[t]{$\bullet$}}
\put(192,      1){\makebox(0,    0)[t]{$\bullet$}}

\put(126,      1){\makebox(0,    0)[t]{$\bullet$}}
\put(28,      -1){\line(1,      0){33}}
\put(61,      -1){\line(1,      0){30}}
\put(94,     -1){\line(1,      0){30}}
\put(126,      -1){\line(1,      0){33}}
\put(159,      -1){\line(1,      0){33}}

\put(22,     10){$q$}
\put(36,      5){$q ^{-1}$}
\put(55,     10){$q$}

\put(75,      5){$q ^{-1}$}

\put(91,      10){$q ^{}$}

\put(102,     5){$q ^{-1}$}

\put(120,      10){$-1 ^{}$}
\put(135,     5){${-1}$}
\put(150,      10){${-1}$}

\put(168,     5){$q ^{-1}$}
\put(192,      10){$q ^{}$}

\put(210,        -1)  {  $q \in R_4,$ GDD $1$ of Row $19\emph{}$  in Table  C. }
\end{picture}$\\

 $\ \ \ \ \  \ \   \begin{picture}(100,      15)  \put(-45,      -1){\rm (d) }

\put(-45,      -1){}
\put(27,      1){\makebox(0,     0)[t]{$\bullet$}}
\put(60,      1){\makebox(0,      0)[t]{$\bullet$}}
\put(93,     1){\makebox(0,     0)[t]{$\bullet$}}

\put(126,      1){\makebox(0,    0)[t]{$\bullet$}}

\put(159,      1){\makebox(0,    0)[t]{$\bullet$}}
\put(28,      -1){\line(1,      0){33}}
\put(61,      -1){\line(1,      0){30}}
\put(94,     -1){\line(1,      0){30}}
\put(126,      -1){\line(1,      0){33}}

\put(22,     10){$q ^{}$}
\put(40,      5){$q ^{-1}$}
\put(58,      10){$q$}
\put(70,      5){$q ^{-1}$}

\put(86,      10){${-1}$}

\put(102,     5){$q ^{}$}

\put(112,      10){$-q ^{}$}
\put(130,     5){$q ^{}$}
\put(153,      10){${-1}$}

\put(200,        -1)  {   $q \in R_3,$ GDD $7$ of Row $13$  in Table C. }
\end{picture}$

\end {Lemma}

\begin {Lemma} \label {2.79}
Every quasi-affine GDD with rank $n >5$ over classical GDD is one of following.
{\rm (i)} Bi-classical GDDs. {\rm (ii)}  Simple cycles. {\rm (iii)} quasi-affine GDDs over non classical GDDs. {\rm (iv)}
Adding h-Type 7 on two end vertexes of \\ \\

  $\ \ \ 
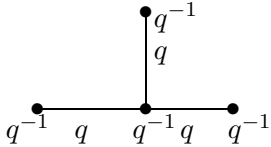
\begin{picture}(100,      15)
\put(-45,      -1){ }
\put(60,      1){\makebox(0,      0)[t]{$\bullet$}}

\put(28,      -1){\line(1,      0){33}}
\put(27,      1){\makebox(0,     0)[t]{$\bullet$}}
\put(-14,      1){\makebox(0,     0)[t]{$\bullet$}}

\put(-14,     -1){\line(1,      0){50}}

\put(-26,     -12){$q^{-1}$}
\put(0,      -12){$q^{}$}
\put(22,     -12){$q^{-1}$}
\put(40,      -12){$q^{}$}
\put(58,      -12){$q^{-1}$}

\put(27,    38){\makebox(0,     0)[t]{$\bullet$}}

\put(27,      0){\line(0,     1){35}}

\put(30,      30){$q^{-1}$}

\put(30,      18){$q^{}$}

\put(130,        -1)  {.   }

\end{picture}$

\end {Lemma}

{\bf Proof.} If it does not hold, then  there exists an quasi-affine   GDD over classical Type i which   is  not  bi-classical, simple cycle and over non classical GDD.

{\rm (a)} If the GDD  is a  cycle, then we can obtain a chain such that h-Type i is the middle of the chain by omitting a vertex in cycle. This chain is not classical, which is contradiction.

{\rm (b)} If the GDD is a chain,  then the GDD  is bi- classical, which is contradiction.

{\rm (c)} If the GDD is not cycle and containing a cycle , then
there is an vertex which does not contain in a cycle. Omitting the vertex in whole is Type 6.
Consequently the GDD is a bi-classical since omitting one vertex in head of the Type 6 is classical. It is a contradiction.

{\rm (d)} If the GDD is not chain and does not  contains any cycles, then it is an quasi-affine over  Type 5.  It has to be a GDD in Lemma. $\Box$

\begin {Lemma} \label {2.77} Assume rank $n >4.$
strict   semi-classical chains are listed:\\

 $\ \ \ \ \  \ \   \begin{picture}(100,      15)  \put(-45,      -1){\rm (a)}
 \put(27,      1){\makebox(0,     0)[t]{$\bullet$}}
\put(60,      1){\makebox(0,      0)[t]{$\bullet$}}
\put(93,     1){\makebox(0,     0)[t]{$\bullet$}}

\put(126,      1){\makebox(0,    0)[t]{$\bullet$}}

\put(159,      1){\makebox(0,    0)[t]{$\bullet$}}

\put(126,      1){\makebox(0,    0)[t]{$\bullet$}}
\put(28,      -1){\line(1,      0){33}}
\put(61,      -1){\line(1,      0){30}}
\put(94,     -1){\line(1,      0){30}}
\put(126,      -1){\line(1,      0){33}}

\put(22,     10){$-1$}
\put(40,      5){$q ^{-1}$}
\put(58,      10){$-1$}
\put(70,      5){$q ^{-1}$}

\put(91,      10){$q$}

\put(102,     5){$q ^{-1}$}

\put(128,      10){$q$}
\put(135,     5){$q ^{-1}$}
\put(159,      10){$q$}

\put(200,        -1)  {  $q \in  R_{3}$.   GDD $2$ of Row $12$  in Table C. }
\end{picture}$\\

 $\ \ \ \ \  \ \   \begin{picture}(100,      15)  \put(-45,      -1){\rm (b)}
 \put(27,      1){\makebox(0,     0)[t]{$\bullet$}}
\put(60,      1){\makebox(0,      0)[t]{$\bullet$}}
\put(93,     1){\makebox(0,     0)[t]{$\bullet$}}

\put(126,      1){\makebox(0,    0)[t]{$\bullet$}}

\put(159,      1){\makebox(0,    0)[t]{$\bullet$}}

\put(126,      1){\makebox(0,    0)[t]{$\bullet$}}
\put(28,      -1){\line(1,      0){33}}
\put(61,      -1){\line(1,      0){30}}
\put(94,     -1){\line(1,      0){30}}
\put(126,      -1){\line(1,      0){33}}

\put(22,     10){$-1$}
\put(40,      5){$q ^{-1}$}
\put(58,      10){$-1$}
\put(70,      5){$q ^{-1}$}

\put(91,      10){$q$}

\put(102,     5){$q ^{-1}$}

\put(120,      10){$-1$}
\put(135,     5){$q ^{}$}
\put(159,      10){$-1$}

\put(200,        -1)  {  $q \in  R_{3}$.   GDD $5$ of Row $13$  in Table C. }
\end{picture}$\\

 $\ \ \ \ \  \ \   \begin{picture}(100,      15)  \put(-45,      -1){\rm (c)}{ }
\put(27,      1){\makebox(0,     0)[t]{$\bullet$}}
\put(60,      1){\makebox(0,      0)[t]{$\bullet$}}
\put(93,     1){\makebox(0,     0)[t]{$\bullet$}}

\put(126,      1){\makebox(0,    0)[t]{$\bullet$}}

\put(159,      1){\makebox(0,    0)[t]{$\bullet$}}

\put(126,      1){\makebox(0,    0)[t]{$\bullet$}}
\put(28,      -1){\line(1,      0){33}}
\put(61,      -1){\line(1,      0){30}}
\put(94,     -1){\line(1,      0){30}}
\put(126,      -1){\line(1,      0){33}}

\put(22,     10){$-1$}
\put(40,      5){$q ^{-1}$}
\put(58,      10){$-1$}
\put(70,      5){$q ^{-1}$}

\put(88,      10){$-1$}

\put(102,     5){$q ^{}$}

\put(120,      10){$-1$}
\put(135,     5){$q ^{-1}$}
\put(159,      10){$q$}

\put(200,        -1)  {  $q \in  R_{3}$.   GDD $8$ of Row $13$  in Table C. }
\end{picture}$\\

 $\ \ \ \ \  \ \   \begin{picture}(100,      15)  \put(-45,      -1){\rm (d)}{}
\put(27,      1){\makebox(0,     0)[t]{$\bullet$}}
\put(60,      1){\makebox(0,      0)[t]{$\bullet$}}
\put(93,     1){\makebox(0,     0)[t]{$\bullet$}}

\put(126,      1){\makebox(0,    0)[t]{$\bullet$}}

\put(159,      1){\makebox(0,    0)[t]{$\bullet$}}

\put(126,      1){\makebox(0,    0)[t]{$\bullet$}}
\put(28,      -1){\line(1,      0){33}}
\put(61,      -1){\line(1,      0){30}}
\put(94,     -1){\line(1,      0){30}}
\put(126,      -1){\line(1,      0){33}}

\put(22,     10){$-1$}
\put(40,      5){$q ^{-1}$}
\put(58,      10){$q ^{-1}$}
\put(70,      5){$q ^{}$}

\put(87,      10){$-1$}

\put(102,     5){$q ^{-1}$}

\put(128,      10){$q$}
\put(135,     5){$q ^{-1}$}
\put(159,      10){$q$}

\put(200,        -1)  {   $q \in  R_{3}$. GDD $11$ of Row $13$  in Table C. }
\end{picture}$\\

 $\ \ \ \ \  \ \   \begin{picture}(100,      15)  \put(-45,      -1){\rm (e)}{}
\put(27,      1){\makebox(0,     0)[t]{$\bullet$}}
\put(60,      1){\makebox(0,      0)[t]{$\bullet$}}
\put(93,     1){\makebox(0,     0)[t]{$\bullet$}}

\put(126,      1){\makebox(0,    0)[t]{$\bullet$}}

\put(159,      1){\makebox(0,    0)[t]{$\bullet$}}

\put(126,      1){\makebox(0,    0)[t]{$\bullet$}}
\put(28,      -1){\line(1,      0){33}}
\put(61,      -1){\line(1,      0){30}}
\put(94,     -1){\line(1,      0){30}}
\put(126,      -1){\line(1,      0){33}}

\put(22,     10){$q$}
\put(40,      5){$-q ^{}$}
\put(58,      10){$-1$}
\put(70,      5){$-q ^{}$}

\put(91,      10){$q$}

\put(102,     5){$-q ^{}$}

\put(128,      10){$q$}
\put(135,     5){$-q ^{}$}
\put(159,      10){$q$}

\put(200,        -1)  {   $q \in R_4,$ GDD $4$ of Row $14$  in Table C. }
\end{picture}$\\

 $\ \ \ \ \  \ \   \begin{picture}(100,      15)  \put(-45,      -1){\rm (f)}{}
\put(27,      1){\makebox(0,     0)[t]{$\bullet$}}
\put(60,      1){\makebox(0,      0)[t]{$\bullet$}}
\put(93,     1){\makebox(0,     0)[t]{$\bullet$}}

\put(126,      1){\makebox(0,    0)[t]{$\bullet$}}

\put(159,      1){\makebox(0,    0)[t]{$\bullet$}}

\put(126,      1){\makebox(0,    0)[t]{$\bullet$}}
\put(28,      -1){\line(1,      0){33}}
\put(61,      -1){\line(1,      0){30}}
\put(94,     -1){\line(1,      0){30}}
\put(126,      -1){\line(1,      0){33}}

\put(22,     10){$-1$}
\put(40,      5){$q ^{}$}
\put(50,      10){$-1$}
\put(70,      5){$q ^{-2}$}

\put(91,      10){$q ^{2}$}

\put(102,     5){$q ^{-2}$}

\put(128,      10){$q ^{2}$}
\put(135,     5){$q ^{-2}$}
\put(159,      10){$q ^{2}$}

\put(200,        -1)  {  $q \in R_5,$ GDD $6$ of Row $15$  in Table C. }
\end{picture}$\\

 $\ \ \ \ \  \ \   \begin{picture}(100,      15)  \put(-45,      -1){\rm (g)}{}
\put(27,      1){\makebox(0,     0)[t]{$\bullet$}}
\put(60,      1){\makebox(0,      0)[t]{$\bullet$}}
\put(93,     1){\makebox(0,     0)[t]{$\bullet$}}

\put(126,      1){\makebox(0,    0)[t]{$\bullet$}}

\put(159,      1){\makebox(0,    0)[t]{$\bullet$}}

\put(126,      1){\makebox(0,    0)[t]{$\bullet$}}
\put(28,      -1){\line(1,      0){33}}
\put(61,      -1){\line(1,      0){30}}
\put(94,     -1){\line(1,      0){30}}
\put(126,      -1){\line(1,      0){33}}

\put(22,     10){$-1$}
\put(40,      5){$q ^{-1}$}
\put(58,      10){$-1$}
\put(70,      5){$q ^{-1}$}

\put(91,      10){$q ^{}$}

\put(102,     5){$q ^{-1}$}

\put(128,      10){$q ^{}$}
\put(135,     5){$q ^{-1}$}
\put(159,      10){$-1 ^{}$}

\put(200,        -1)  { $q \in R_3,$  GDD $3$ of Row $13$  in Table C. }
\end{picture}$\\

 $\ \ \ \ \  \ \   \begin{picture}(100,      15)  \put(-45,      -1){\rm (h)}{}
\put(27,      1){\makebox(0,     0)[t]{$\bullet$}}
\put(60,      1){\makebox(0,      0)[t]{$\bullet$}}
\put(93,     1){\makebox(0,     0)[t]{$\bullet$}}

\put(126,      1){\makebox(0,    0)[t]{$\bullet$}}

\put(159,      1){\makebox(0,    0)[t]{$\bullet$}}
\put(192,      1){\makebox(0,    0)[t]{$\bullet$}}

\put(126,      1){\makebox(0,    0)[t]{$\bullet$}}
\put(28,      -1){\line(1,      0){33}}
\put(61,      -1){\line(1,      0){30}}
\put(94,     -1){\line(1,      0){30}}
\put(126,      -1){\line(1,      0){33}}
\put(159,      -1){\line(1,      0){33}}

\put(22,     10){$-1$}
\put(40,      5){$q ^{-1}$}
\put(58,      10){$-1$}
\put(74,      5){$q ^{-1}$}

\put(91,      10){$q$}

\put(102,     5){$q ^{-1}$}

\put(124,      10){$q$}
\put(135,     5){$q ^{-1}$}
\put(159,      10){$q$}

\put(168,     5){$q ^{-1}$}
\put(192,      10){$q$}

\put(168,     5){$q ^{-1}$}
\put(192,      10){$q$}

\put(200,        -1)  { $q \in R_3,$  GDD $2$ of Row $18$  in Table B. }
\end{picture}$\\

 $\ \ \ \ \  \ \   \begin{picture}(100,      15)  \put(-45,      -1){\rm (i)}{}
\put(27,      1){\makebox(0,     0)[t]{$\bullet$}}
\put(60,      1){\makebox(0,      0)[t]{$\bullet$}}
\put(93,     1){\makebox(0,     0)[t]{$\bullet$}}

\put(126,      1){\makebox(0,    0)[t]{$\bullet$}}

\put(159,      1){\makebox(0,    0)[t]{$\bullet$}}
\put(192,      1){\makebox(0,    0)[t]{$\bullet$}}

\put(126,      1){\makebox(0,    0)[t]{$\bullet$}}
\put(28,      -1){\line(1,      0){33}}
\put(61,      -1){\line(1,      0){30}}
\put(94,     -1){\line(1,      0){30}}
\put(126,      -1){\line(1,      0){33}}
\put(159,      -1){\line(1,      0){33}}

\put(22,     10){$q$}
\put(40,      5){$q ^{-1}$}
\put(58,      10){$-1$}
\put(74,      5){$q ^{-1}$}

\put(91,      10){$q$}

\put(102,     5){$q ^{-1}$}

\put(124,      10){$q$}
\put(135,     5){$q ^{-1}$}
\put(159,      10){$q$}

\put(168,     5){$q ^{-1}$}
\put(192,      10){$q$}

\put(168,     5){$q ^{-1}$}
\put(192,      10){$q$}

\put(200,        -1)  {  $q \in R_4,$ GDD $4$ of Row $19$  in Table C. }
\end{picture}$\\

\end {Lemma}

\begin {Lemma}\label {2.13}  Assume rank $n >5.$ If GDD is an arithmetic non chain,  then omitting some vertex is connected simple chain.

\end {Lemma}

\begin {Lemma}\label {2.63} Assume rank $n >4.$

{\rm (i)}  All  arithmetic GDDs are quasi-classical.

{\rm (ii)}  All  classical GDDs are semi-classical.

\end {Lemma}

\begin {Lemma} \label {2.93} Assume rank $n >5.$
 {\rm (i)}  Bi-classical GDDs are quasi-affine.

 {\rm (ii)} Classical  + semi-classical GDDs are quasi-affine.

 {\rm (iii)} Bi-semi-classical GDDs are quasi-affine.

\end {Lemma}
{\bf Proof.} We only need prove that  these GDDs are not arithmetic.  If  a non chain GDD satisfied conditions is  an arithmetic GDD, then we obtain a contradiction by Lemma \ref {2.13}. If  a chain GDD satisfied conditions is  an arithmetic GDD, then we obtain a contradiction by Lemma \ref {2.6}.  \hfill $\Box$

\begin {Lemma} \label {2.93'}
 Adding a vertex on tail of a strict quasi-classical GDD in Table C is not quasi-affine except adding h-type 7.

\end {Lemma}
{\bf Proof.}  It follows from Lemma \ref {2.6} and Lemma \ref {2.13}.
  \hfill $\Box$

\begin {Lemma} \label {2.93''}
 Adding a vertex on tail of a strict semi-classical GDD in Table C is not quasi-affine except\\

 $\ \ \ \ \  \ \
$\\ \\
and   adding h-type i for $ i= 1, 2, 4, 5, 6, 7.$

\end {Lemma}
{\bf Proof.}  It follows from Lemma \ref {2.77} and Lemma \ref {2.13}.
  \hfill $\Box$

\section {  Main result} \label {s4}

\begin {Lemma} \label {3.1}
All classical + semi-classical  GDDs with rank $n >5$ are listed.
\\ \\

$\ \ \ \ \  \ \   %
$

\end {Lemma}

{\ }

\begin {Lemma} \label {3.2} All quasi-affine continual GDDs with rank $n >5$ are listed.
\\

$\ \ \ \ \  \ \   %
$\\

\end {Lemma}

\begin {Theorem} \label {Main} All quasi-affine  connected GDDs with rank $n >5$ are listed.

 {\rm (i)}  Bi-classical GDDs.
{\rm (ii)} Simple cycles.
{\rm (iii)}  Quasi-affine continual GDDs.
{\rm (iv)} Classical + semi-classical GDDs.
{\rm (v)} Other cases
.\\ \\ \\ \\

 $\ \ \ \ \  \ \   %
$\\

\end {Theorem}

{\bf Proof.}  Sufficiency. Considering Lemma \ref {2.93} we  obtain
that  {\rm (i)}, {\rm (iv)} and {\rm (v)} are quasi-affine. It is clear that {\rm (ii)} and {\rm (iii)} are quasi-affine. cl 1 and cl 2 in   {\rm (vi)} are quasi-affine by  Lemma \ref {2.79}. The others in {\rm (vi)} are quasi-affine by subsection \ref {ss1}-
subsection \ref {ss99}.

Necessity. Set $A:=\{x \mid x \hbox{ is a quasi-affine GDD in subsection } \ref {ss1}-
subsection \ref {ss99} \}$. $B:=\{x \mid x \hbox{ is a  GDD in  } {\rm (iii)},  {\rm (iv)}, {\rm (v)} \hbox{ and } {\rm (vi)}; x \hbox{ is not } \hbox { cl.}1, \hbox { cl.}2 \}$. It is clear that $A = B.$
If GDD is quasi-affine over a non classical GDD, then it is  in subsection \ref {ss1}-
subsection \ref {ss99}. If GDD is quasi-affine over a  classical GDD, then it is  in subsection \ref {ss1}-
subsection \ref {ss99} or in  {\rm (i)}, {\rm (ii)} and {\rm (vi)} by Lemma \ref {2.79}.

 We write the proof according  to the following method.

(1) If a GDD is  quasi-affine over two GDDs, we write it in below one.

 (2) If a GDD is  quasi-affine non chain , we write it in one over non chain.

(3) If a non cycle  GDD is  quasi-affine with cycle, we write it in ones over GDD with cycle.

(4) If a non cycle  GDD is  quasi-affine over classical GDD, we write it in ones over non classical GDD or bi-classical GDD.

(5)
We do not consider tails by Lemma \ref {2.93'} and Lemma \ref {2.93''}.
\subsection  {Quasi-affine  GDDs  over
GDD $1$ of Row $11$ in Table C} \label {ss1}
   {\rm (i)}   Omitting  Vertex 1 and adding on  Vertex 5. We find all quasi-affine GDDs  adding a vertex on vertex 5 of GDD $1$ of Row $11$. We have to consider all  GDDs adding a vertex on vertex 5 of GDD $1$ of Row $11$, in which omitting Vertex 1 is an arithmetic GDD.\\

   $\ \ \ \ \  \ \
$\\
is GDD $19$ of Row $18.$
{\rm (c)} There are not any other cases by Lemma \ref {2.6}.

 {\rm (iii)}
Cycle.

Every  quasi-affine cycle over GDD $1$ of Row $11$  consists of a GDD in Part {\rm (i)} and a GDD in Part {\rm (ii)}.  Of course, order of Vertex 6 of GDD in Part {\rm (i)} is the same as  order of Vertex 1 of GDD  Part {\rm (ii)}.
\\ \\ \\

 $   %
$\\
Omitting  Vertex 2,  Vertex 3,  Vertex 4 are  arithmetic.  It is quasi-affine.
\\ \\ \\

{\rm (b)} to {\rm (b)}
 $\ \ \ \ \  \ \   %
$
   \subsection  {Quasi-affine  GDDs  over
GDD $2$ of Row $11$ in Table C} \label {sub3.2}

 {\rm (i)}   Omitting Vertex 1 and adding  on  Vertex 5.

 {\rm (ii)}    Adding on  Vertex 2, Vertex  4. It is not quasi-affine since it contains a proper subGDD which is not an arithmetic GDD.
 {\rm (iv)} Omitting Vertex 1 and adding  on  Vertex 3.\\ \\ \\

 $\ \ \ \ \  \ \   %
$\\ \\
Omitting   Vertex 4  in whole is an arithmetic GDD   by  GDD $2$ of Row $11$.
It is quasi-affine.\\ \\ \\

 $\ \ \ \ \  \ \   %
$\\ \\
It is repeated.
{\rm (c)} There are not any other cases by Lemma \ref {2.13}.
  \subsection  {Quasi-affine  GDDs  over
GDD $3$ of Row $11$ in Table C} \label {sub3.2}
{\rm (i)} Omitting  Vertex  1 and adding   on Vertex 4.\\ \\

 $\ \ \ \ \  \ \   %
$\\ \\ Omitting 4  in whole is  GDD $1$ of Row $11$.
 It is quasi-affine    by  Lemma \ref {2.2} {\rm (ii)}.\\ \\

 $\ \ \ \ \  \ \   %
$\\ \\
 It is repeated.

{\rm (c)} There are not any other cases by Lemma \ref {2.13}.

 {\rm (ii)}   Adding  on  Vertex 3. It is not quasi-affine since it contains a proper subGDD which is not an arithmetic GDD.

 {\rm (iii)} Omitting Vertex 5 and adding on  Vertex 2.

 There are not any  cases by Lemma \ref {2.13}.

 {\rm (iv)} Omitting Vertex  5 and adding  on Vertex 1.
     \subsection  {Quasi-affine  GDDs  over
GDD $4$ of Row $11$  in Table C }\label {sub3.2}
 {\rm (i)}  Omitting  Vertex  1 and adding  on  Vertex 4.\\ \\

 $\ \ \ \ \  \ \   %
$\\ \\
 It is quasi-affine    by  Lemma \ref {2.2}  {\rm (ii)}.\\ \\

 $\ \ \ \ \  \ \   %
$\\ \\
 It is repeated.
{\rm (c)} There are not any other cases by Lemma \ref {2.13}.
 {\rm (ii)} Omitting Vertex 1 and adding  on  Vertex 5.

 {\rm (iii)}  Adding  on  Vertex 3. It is not quasi-affine since it contains a proper subGDD which is not an arithmetic GDD.

 {\rm (iv)} Adding on  Vertex 2.\\ \\

 $\ \ \ \ \  \ \   %
$\\ \\
Omitting Vertex 1 in whole  is GDD $10$ of Row $12$  in Table C. It is repeated.
{\rm (b)} There are not any other cases by Lemma \ref {2.13}.
\\ \\

$\ \ \ \ \  \ \   %
$\\ \\ \\ Omitting 6 Vertex in whole  is  not arithmetic by Lemma \ref {2.2} {\rm (i)}.
\\ \\

$\ \ \ \ \  \ \   %
$\\ \\ \\ Omitting Vertex  1 in whole  is  GDD $10$ of Row $12$. It is repeated.
   \subsection  {Quasi-affine  GDDs  over
GDD $1$ of Row $12$ in Table C }\label {sub3.2}
 {\rm (i)}  Omitting Vertex 1 and adding  on Vertex  5.\\

 $\ \ \ \ \  \ \   %
$
\\  It is repeated.\\

 $\ \ \ \ \  \ \   %
$\\
It is repeated.
{\rm (c)} There are not any other cases by Lemma \ref {2.63}.
    \subsection  {Quasi-affine  GDDs  over
GDD $2$ of Row $12$ in Table C}\label {sub3.2}
 {\rm (i)}  Omitting Vertex   1 and adding  on  Vertex 5.\\

 $\ \ \ \ \  \ \   %
$\\
 It is repeated.\\

 $\ \ \ \ \  \ \   %
$\\
It is repeated.

    \subsection  {Quasi-affine  GDDs  over
GDD $3$ of Row $12$ in Table C}\label {sub3.2}
 {\rm (i)}  Omitting  Vertex  1 and adding  on  Vertex 4.
It    is empty   by Lemma \ref {2.63}.
   {\rm (ii)}  Omitting  Vertex  1 and adding  on  Vertex 5.\\ \\

 $\ \ \ \ \  \ \   %
$\\ \\
 It is repeated.\\

 $\ \ \ \ \  \ \   %
$\\ \\
 It is repeated.
{\rm (c)} There are not any other cases by Lemma \ref {2.13}.
 {\rm (iii)}  Adding  on  Vertex 3. It is not quasi-affine since it contains a proper subGDD which is not an arithmetic GDD.
 {\rm (iv)}  Adding  on Vertex 2.\\ \\

 $\ \ \ \ \  \ \   %
$\\ \\
 It is repeated.
\\ \\

 $\ \ \ \ \  \ \   %
$\\ \\
Onitting Vertex  5 in whole is not arithmetic   by  Lemma \ref {2.2}
 {\rm (i)}.\\ \\

 $\ \ \ \ \  \ \   %
$\\ \\
 Onitting  Vertex 6 in whole  is not an arithmetic GDD   by Row 11-13.
{\rm (e)} There are not any other cases by Lemma \ref {2.13}.
   \subsection  {Quasi-affine  GDDs  over
GDD $4$ of Row $12$ in Table C}\label {sub3.2}
  {\rm (i)}  Omitting  Vertex  1 and adding  on  Vertex 4.

 It    is empty  by  Lemma \ref {2.2}
 {\rm (i)}.

 {\rm (ii)} Omitting Vertex 1 and adding  on Vertex  5.\\ \\

 $\ \ \ \ \  \ \   %
$\\ \\
 It is repeated.\\ \\

 $\ \ \ \ \  \ \   %
$\\ \\
 It is repeated.
{\rm (c)} There are not any other cases by Lemma \ref {2.13}.
 {\rm (iii)}    Adding  on  Vertex 3. It is not quasi-affine since it contains a proper subGDD which is not an arithmetic GDD.
 {\rm (iv)} Adding on  Vertex 2.\\ \\

 $\ \ \ \ \  \ \   %
$\\ \\
 It is repeated.
 \\ \\

$\ \ \ \ \  \ \   %
$\\ \\
 Omitting  Vertex 6 in whole is not  arithmetic by Row 11-13.
{\rm (c)} There are not any other cases by Lemma \ref {2.13}.
    \subsection  {Quasi-affine  GDDs  over
GDD $5$ of Row $12$ in Table C}\label {sub3.2}
  {\rm (i)}  Omitting  Vertex  1 and adding  on  Vertex 4.
  It    is empty  by  Lemma \ref {2.2}
 {\rm (i)}.

   {\rm (ii)}  Omitting  Vertex  1 and adding  on  Vertex 5.\\ \\

 $\ \ \ \ \  \ \   %
$\\ \\
 It is quasi-affine    by Lemma \ref {2.63}.\\ \\

 $\ \ \ \ \  \ \   %
$\\ \\
It is repeated.
{\rm (c)} There are not any other cases by Lemma \ref {2.13}.

 {\rm (iii)}    Adding  on  Vertex 3. It is not quasi-affine since it contains a proper subGDD which is not an arithmetic GDD.

 {\rm (iv)} Omitting Vertex 1 and adding  on Vertex 2.
 Omitting  Vertex 6 in whole
 is not  arithmetic  by Lemma \ref {2.13}.
   \subsection  {Quasi-affine  GDDs  over
GDD $6$ of Row $12$ in Table C}\label {sub3.2}
{\rm (i)}  Omitting Vertex 1 and adding   on Vertex 4.
It is empty
by  Lemma \ref {2.2}
 {\rm (i)}.
   {\rm (ii)}   Omitting  Vertex  1 and adding  on  Vertex 5.\\ \\

 $\ \ \ \ \  \ \   %
$\\ \\
 It is quasi-affine    by Lemma \ref {2.63}.\\ \\

 $\ \ \ \ \  \ \   %
$\\ \\
It is repeated.
{\rm (c)} There are not any other cases by Lemma \ref {2.13}.

{\rm (iii)}   Adding   on Vertex 3. It is not quasi-affine since it contains a proper subGDD which is not an arithmetic GDD.

 {\rm (iv)} Omitting  Vertex 5 and adding  on Vertex 2. It    is empty   by Lemma \ref {2.13}.

    \subsection  {Quasi-affine  GDDs  over
GDD $7$ of Row $12$ in Table C}\label {sub3.2}
 {\rm (i)} Omitting  Vertex  1 and adding  on  Vertex 4.\\ \\

 $\ \ \ \ \  \ \   %
$\\ \\
Omitting  Vertex 4 in whole  is   GDD $13$ of Row $12$.
 It is quasi-affine
by  Lemma \ref {2.2} {\rm (ii)}.\\ \\

 $\ \ \ \ \  \ \   %
$\\ \\
Omitting  Vertex 4 in whole  is not  arithmetic  by Row 11- 13.
{\rm (c)} There are not any other cases by Lemma \ref {2.13}.
 {\rm (ii)} Omitting Vertex 1 and adding   on Vertex 5.

 {\rm (iii)}   Adding  on Vertex  3. It is not quasi-affine since it contains a proper subGDD which is not an arithmetic GDD.

 {\rm (iv)} Omitting Vertex 1 and adding  on Vertex 2.
 Omitting  Vertex 6  in whole is not an arithmetic GDD   by Lemma \ref {2.13}.
    \subsection  {Quasi-affine  GDDs  over
GDD $8$ of Row $12$ in Table C}\label {sub3.2}
   {\rm (i)}  Omitting  Vertex  1 and adding on  Vertex  4.
   It    is empty  by  Lemma \ref {2.2}
 {\rm (i)}.

    {\rm (ii)}  Omitting  Vertex  1 and adding  on  Vertex 5.\\ \\

 $\ \ \ \ \  \ \
$\\ \\
It is repeated.\\ \\

 $\ \ \ \ \  \ \   %
$\\ \\
 It is quasi-affine    by Row 17-Row 18 or by Lemma \ref {2.63}.
{\rm (c)} There are not any other cases by Lemma \ref {2.13}.

 {\rm (iii)}   Adding  on  Vertex 3. It is not quasi-affine since it contains a proper subGDD which is not an arithmetic GDD.

 {\rm (iv)} Adding  on Vertex 2.\\ \\

 $\ \ \ \ \  \ \   %
$\\  \\ Omitting Vertex 1 in whole is  GDD $14$ of Row $13$  in Table C. It is repeated.\\ \\

 $\ \ \ \ \  \ \   %
$\\ \\
 Omitting  Vertex 5  in whole is not an arithmetic GDD
by  Lemma \ref {2.2}
 {\rm (i)}.\\ \\

 $\ \ \ \ \  \ \   %
$\\ \\
 Omitting  Vertex 6  in whole is GDD $3$ of Row $12.$
Omitting  Vertex 5  in whole is GDD $6$ of Row $13.$  It is quasi-affine.
{\rm (d)} There are not any other cases by Lemma \ref {2.13}.
   \subsection  {Quasi-affine  GDDs  over
GDD $9$ of Row $12$ in Table C}\label {sub3.2}
 {\rm (i)}  Omitting Vertex 1 and adding  on  Vertex 5.\\

 $\ \ \ \ \  \ \   %
$\\ \\
 It is quasi-affine     by Lemma \ref {2.63}.\\

 $\ \ \ \ \  \ \   %
$\\ \\
 It is repeated.
{\rm (c)} There are not any other cases by Lemma \ref {2.13}.
 {\rm (ii)} Omitting Vertex 1 and adding  on  Vertex 3.
It is empty  by Lemma \ref {2.63}.

 {\rm (iii)}    Adding on  Vertex 2, Vertex 4. It is not quasi-affine since it contains a proper subGDD which is not an arithmetic GDD.
    \subsection  {Quasi-affine  GDDs  over
GDD $10$ of Row $12$ in Table C}\label {sub3.2}
  {\rm (i)} Omitting  Vertex  1 and adding  on  Vertex 4.\\ \\

 $\ \ \ \ \  \ \   %
$\\ \\
 It is quasi-affine    Lemma \ref {2.63}\\ \\

 $\ \ \ \ \  \ \   %
$\\ \\
 It is quasi-affine    by  Lemma \ref {2.2} {\rm (ii)}.
{\rm (c)} There are not any other cases by Lemma \ref {2.13}.

 {\rm (ii)} Omitting Vertex 1 and adding  on  Vertex 5.
{\rm (a)} There are not any other cases by Lemma \ref {2.13}.

{\rm (iii)}   Adding  on  Vertex 3. It is not quasi-affine since it contains a proper subGDD which is not an arithmetic GDD.

  {\rm (iv)} Adding on  Vertex 2.\\ \\

 $\ \ \ \ \  \ \   %
$\\ \\
 Omitting Vertex  5  in whole is   GDD $10$ of Row $12$.
Omitting  Vertex 6  in whole is   GDD $4$ of Row $11$.
 It is quasi-affine.\\ \\

 $\ \ \ \ \  \ \   %
$\\ \\
Omitting  Vertex 5  in whole is not an arithmetic GDD   by Lemma \ref {2.2}
 {\rm (i)}.\\ \\

$\ \ \ \ \  \ \   %
$\\ \\ Omitting  Vertex 5 in whole  is not  arithmetic  by Row 11-Row 13.
{\rm (d)} There are not any other cases by Lemma \ref {2.13}.
   \subsection  {Quasi-affine  GDDs  over
GDD $11$ of Row $12$ in Table C}\label {sub3.2}
 {\rm (i)}   Adding  on Vertex 2, Vertex 4. It is not quasi-affine since it contains a proper subGDD which is not an arithmetic GDD.

 {\rm (ii)} Omitting Vertex 1 and adding on Vertex 5.

 {\rm (iii)} Omitting Vertex 1 and adding  on  Vertex 3.\\ \\ \\

$\ \ \ \ \  \ \   %
$\\ \\
 Omitting  Vertex 6 in whole  is  GDD $11$ of Row $12.$
     It is quasi-affine. \\ \\ \\

$\ \ \ \ \  \ \   %
$\\ \\
 Omitting Vertex  6 in whole  is  not  arithmetic   by Row 11-13. {\rm (c)} There are not any other cases by Lemma \ref {2.13}.
 \subsection  {Quasi-affine  GDDs  over
GDD $12$ of Row $12$ in Table C}\label {sub3.2}
 {\rm (i)} Omitting   Vertex 1 and adding  on Vertex  5.\\

 $\ \ \ \ \  \ \   %
$\\
  It is quasi-affine    by Lemma \ref {2.63}.\\

 $\ \ \ \ \  \ \   %
$\\
 It is repeated.
{\rm (c)} There are not any other cases by Lemma \ref {2.6}.

 {\rm (iii)} Cycle. \\ \\

{\ }

{\ }

{\ }

 $\ \ \ \ \  \ \   %
$\\
Omitting  Vertex 2,  Vertex 3,  Vertex 4  in whole are  arithmetic.  It is quasi-affine.
  \subsection  {Quasi-affine  GDDs  over
GDD $13$ of Row $12$ in Table C}\emph{\label {sub3.2}}
{\rm (i)}  Omitting Vertex 1 and adding  on  Vertex 5.

  {\rm (iii)} Cycle.
\\ \\ \\

{\rm (b)} to {\rm (b)} $  %
$\\
 Omitting  Vertex 2, Vertex  3,  Vertex 4  in whole are  arithmetic.  It is quasi-affine.
   \subsection  {Quasi-affine  GDDs  over
GDD $14$ of Row $12$ in Table C}\label {sub3.2}
 {\rm (i)}  Omitting  Vertex  5 and adding  on  Vertex 3.\\ \\

 $\ \ \ \ \  \ \   %
$\\ \\
Omitting Vertex 4  in whole is GDD $1$ of Row $11$.
  It is quasi-affine    by Lemma \ref {2.13}.  \\ \\

 $\ \ \ \ \  \ \   %
$\\ \\
{\rm (c)} There are not any other cases by Lemma \ref {2.13}.

 {\rm (i)}   Adding on  Vertex 2. It is not quasi-affine since it contains a proper subGDD which is not an arithmetic GDD.

 {\rm (ii)} Omitting  Vertex 1 and adding  on  Vertex 4. It is not arithmetic  by Lemma \ref {2.13}.

 {\rm (iii)} Omitting Vertex 1 and adding  on  Vertex 5.
    \subsection  {Quasi-affine  GDDs  over
GDD $15$ of Row $12$ in Table C}\label {sub3.2}
  {\rm (i)}  Omitting  Vertex  1 and adding  on  Vertex 5.\\

 $\ \ \ \ \  \ \   %
$\\
 It is quasi-affine     by Row 17-18.\\

 $\ \ \ \ \  \ \   %
$\\
 It is repeated.
{\rm (c)} There are not any other cases by Lemma \ref {2.6}.

  {\rm (iii)} Cycle.\\ \\ \\

{\rm (a)} to {\rm (a)}
 $  %
$\\
Omitting  Vertex 2,  Vertex 3, Vertex  4  in whole are  arithmetic.  It is quasi-affine.
     \subsection  {Quasi-affine  GDDs  over
GDD $1$ of Row $13$ in Table C}\label {sub3.2}
  {\rm (i)}  Omitting   Vertex 1 and adding  on  Vertex 5.
\\

 $\ \ \ \ \  \ \   %
$\\
  It is quasi-affine    by Lemma \ref {2.63}.\\

 $\ \ \ \ \  \ \   %
$\\
 It is repeated in GDD $20$ of Row $13$.
{\rm (c)} There are not any other cases by Lemma \ref {2.6}.
   \subsection  {Quasi-affine  GDDs  over
GDD $2$ of Row $13$ in Table C}\label {sub3.2}
{\rm (i)} Omitting  Vertex  1 and adding on   Vertex 5.
 It    is empty  by Row 11-13 and Lemma \ref {2.6}.
      \subsection  {Quasi-affine  GDDs  over
GDD $3$ of Row $13$ in Table C}\label {sub3.2}
  {\rm (i)} Omitting   Vertex 1 and adding on   Vertex 5.

  $\ \ \ \ \  \ \   %
$\\
 It is repeated.\\

 $\ \ \ \ \  \ \   %
$\\
 It is quasi-affine    Lemma \ref {2.63}.
{\rm (c)} There are not any other cases by Lemma \ref {2.6}.\\

 $\ \ \ \ \  \ \   %
$\\  It is quasi-affine     by Lemma \ref {2.63}.

 {\rm (iii)}  Cycle.
 \\ \\ \\

  $  %
$\\
Omitting  Vertex 2  in whole is not   arithmetic  by Lemma \ref {2.6}.
 \\ \\ \\ \\

  $  %
$\\
Omitting  Vertex 2, Vertex 3, Vertex 4  in whole are  arithmetic.
     \subsection  {Quasi-affine  GDDs  over
GDD $4$ of Row $13$ in Table C}\label {sub3.2}
 {\rm (i)} Omitting Vertex 1 and adding on   Vertex 5.
  It is empty by Row 11-13.
    \subsection  {Quasi-affine  GDDs  over
GDD $5$ of Row $13$ in Table C}\label {sub3.2}
 {\rm (i)} Omitting  Vertex  1 and adding on   Vertex 5.
It    is empty  by Row 11-13 and  Lemma \ref {2.6}.
    \subsection  {Quasi-affine  GDDs  over
GDD $6$ of Row $13$ in Table C}\label {sub3.2}
 {\rm (i)} Omitting  Vertex  1 and adding on  Vertex  4.
 It    is empty  by Lemma \ref {2.63}.

   {\rm (ii)}  Omitting  Vertex  1 and adding on  Vertex  5.\\

 $\ \ \ \ \  \ \   
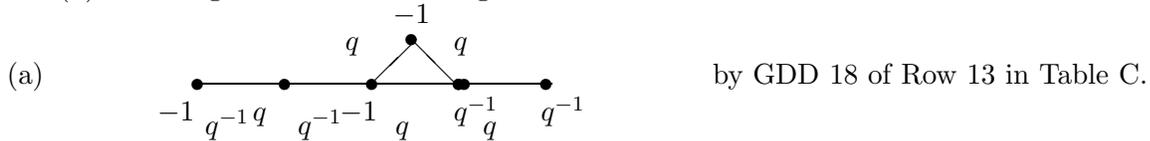
\begin{picture}(100,      15)  \put(-45,      -1){} \put(-45,      -1){{\rm (a)} }
\put(27,      1){\makebox(0,     0)[t]{$\bullet$}}
\put(60,      1){\makebox(0,      0)[t]{$\bullet$}}
\put(93,     1){\makebox(0,     0)[t]{$\bullet$}}

\put(126,      1){\makebox(0,    0)[t]{$\bullet$}}

\put(128,      1){\makebox(0,    0)[t]{$\bullet$}}
\put(159,      1){\makebox(0,    0)[t]{$\bullet$}}

\put(126,      1){\makebox(0,    0)[t]{$\bullet$}}
\put(28,      -1){\line(1,      0){33}}
\put(61,      -1){\line(1,      0){30}}
\put(94,     -1){\line(1,      0){30}}

\put(128,      -1){\line(1,      0){33}}

\put(12,     -15){$-1$}
\put(30,     -20){$q ^{-1}$}
\put(48,      -15){$q$}
\put(65,      -20){$q ^{-1}$}

\put(81,      -15){${-1}$}

\put(102,     -20){$q ^{}$}

\put(124,      -15){$q ^{-1}$}
\put(135,     -20){$q ^{}$}

\put(157,      -15){$q ^{-1}$}

\put(108,    18){\makebox(0,     0)[t]{$\bullet$}}

\put(93,     -1){\line(1,      0){30}}

\put(124,      0){\line(-1,     1){17}}
\put(93,     -1){\line(1,      1){17}}

\put(83,    12){$q$}

\put(101,    22){$-1$}
\put(124,    12){$q$}

\put(222,        -1)  {by GDD $18$ of Row $13$  in Table C. }
\end{picture}$\\ \\
 It is repeated.\\

 $\ \ \ \ \  \ \   
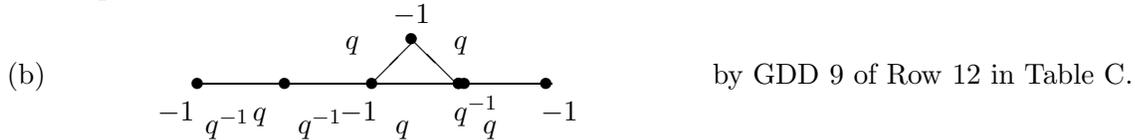
\begin{picture}(100,      15)  \put(-45,      -1){} \put(-45,      -1){{\rm (b)} }
\put(27,      1){\makebox(0,     0)[t]{$\bullet$}}
\put(60,      1){\makebox(0,      0)[t]{$\bullet$}}
\put(93,     1){\makebox(0,     0)[t]{$\bullet$}}

\put(126,      1){\makebox(0,    0)[t]{$\bullet$}}

\put(128,      1){\makebox(0,    0)[t]{$\bullet$}}
\put(159,      1){\makebox(0,    0)[t]{$\bullet$}}

\put(126,      1){\makebox(0,    0)[t]{$\bullet$}}
\put(28,      -1){\line(1,      0){33}}
\put(61,      -1){\line(1,      0){30}}
\put(94,     -1){\line(1,      0){30}}

\put(128,      -1){\line(1,      0){33}}

\put(12,     -15){$-1$}
\put(30,     -20){$q ^{-1}$}
\put(48,      -15){$q$}
\put(65,      -20){$q ^{-1}$}

\put(81,      -15){${-1}$}

\put(102,     -20){$q ^{}$}

\put(124,      -15){$q ^{-1}$}
\put(135,     -20){$q ^{}$}

\put(157,      -15){${-1}$}

\put(108,    18){\makebox(0,     0)[t]{$\bullet$}}

\put(93,     -1){\line(1,      0){30}}

\put(124,      0){\line(-1,     1){17}}
\put(93,     -1){\line(1,      1){17}}

\put(83,    12){$q$}

\put(101,    22){$-1$}
\put(124,    12){$q$}

\put(222,        -1)  {by GDD $9$ of Row $12$  in Table C. }
\end{picture}$\\ \\
 It is quasi-affine    by Lemma \ref {2.63}.
{\rm (c)} There are not any other cases by Lemma \ref {2.63}.

{\rm (iii)} Adding on   Vertex 3. It is not quasi-affine since it contains a proper subGDD which is not an arithmetic GDD.

 {\rm (iv)} Adding on  Vertex 2.\\ \\

 $\ \ \ \ \  \ \   
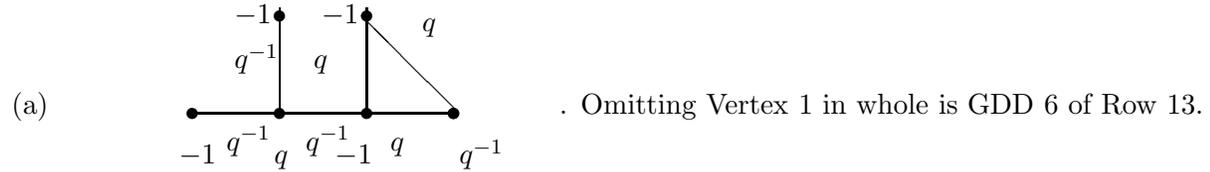
\begin{picture}(100,      15)  \put(-45,      -1){} \put(-45,      -1){ {\rm (a)}}
\put(27,      1){\makebox(0,     0)[t]{$\bullet$}}
\put(60,      1){\makebox(0,      0)[t]{$\bullet$}}
\put(93,     1){\makebox(0,     0)[t]{$\bullet$}}

\put(126,      1){\makebox(0,    0)[t]{$\bullet$}}

\put(126,      1){\makebox(0,    0)[t]{$\bullet$}}
\put(28,      -1){\line(1,      0){33}}
\put(61,      -1){\line(1,      0){30}}
\put(94,     -1){\line(1,      0){30}}

\put(22,    -20){${-1}$}
\put(40,     -15){$q ^{-1}$}
\put(58,      -20){$q ^{}$}
\put(70,      -15){$q ^{-1}$}

\put(81,     -20){${-1}$}

\put(102,    -15){$q ^{}$}

\put(128,     -20){$q^{-1}$}

\put(60,    38){\makebox(0,     0)[t]{$\bullet$}}

\put(60,      -1){\line(0,     1){40}}

\put(43,      33){$-1$}

\put(43,      16){$q ^{-1}$}

\put(93,    38){\makebox(0,     0)[t]{$\bullet$}}

\put(93,      -1){\line(0,     1){40}}

\put(76,      33){${-1}$}

\put(73,      16){$q ^{}$}

\put(126,      1){\line(-1,     1){35}}

\put(114,      30){$q ^{}$}

\put(166,        -1)  {. Omitting  Vertex 1 in whole  is GDD $6$ of Row $13$. }

\end{picture}$\\ \\
Omitting  Vertex 4  in whole is not an arithmetic GDD   by Lemma \ref {2.2}
 {\rm (i)}.\\ \\

 $\ \ \ \ \  \ \   
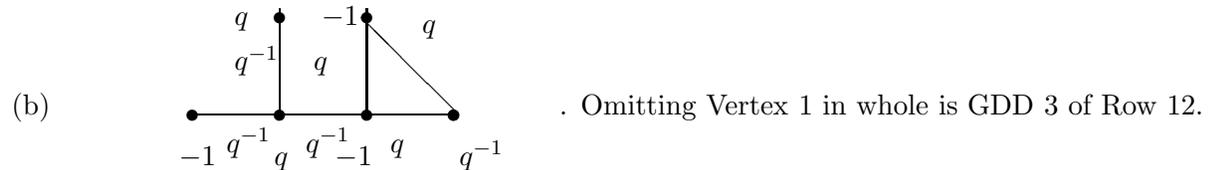
\begin{picture}(100,      15)  \put(-45,      -1){} \put(-45,      -1){ {\rm (b)}}
\put(27,      1){\makebox(0,     0)[t]{$\bullet$}}
\put(60,      1){\makebox(0,      0)[t]{$\bullet$}}
\put(93,     1){\makebox(0,     0)[t]{$\bullet$}}

\put(126,      1){\makebox(0,    0)[t]{$\bullet$}}

\put(126,      1){\makebox(0,    0)[t]{$\bullet$}}
\put(28,      -1){\line(1,      0){33}}
\put(61,      -1){\line(1,      0){30}}
\put(94,     -1){\line(1,      0){30}}

\put(22,    -20){${-1}$}
\put(40,     -15){$q ^{-1}$}
\put(58,      -20){$q ^{}$}
\put(70,      -15){$q ^{-1}$}

\put(81,     -20){${-1}$}

\put(102,    -15){$q ^{}$}

\put(128,     -20){$q^{-1}$}

\put(60,    38){\makebox(0,     0)[t]{$\bullet$}}

\put(60,      -1){\line(0,     1){40}}

\put(43,      33){$q$}

\put(43,      16){$q ^{-1}$}

\put(93,    38){\makebox(0,     0)[t]{$\bullet$}}

\put(93,      -1){\line(0,     1){40}}

\put(76,      33){${-1}$}

\put(73,      16){$q ^{}$}

\put(126,      1){\line(-1,     1){35}}

\put(114,      30){$q ^{}$}

\put(166,        -1)  {.  Omitting  Vertex 1 in whole  is  GDD $3$ of Row $12$. }

\end{picture}$\\ \\
 Omitting  Vertex 5 in whole  is  GDD $14$ of Row $13$  in Table C.
 It is repeated.\\ \\

 $\ \ \ \ \  \ \   
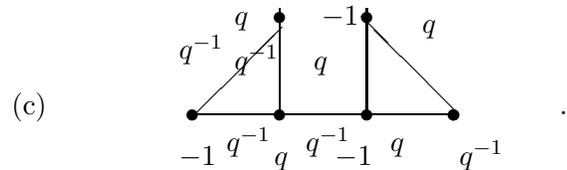
\begin{picture}(100,      15)  \put(-45,      -1){} \put(-45,      -1){ {\rm (c)}}
\put(27,      1){\makebox(0,     0)[t]{$\bullet$}}
\put(60,      1){\makebox(0,      0)[t]{$\bullet$}}
\put(93,     1){\makebox(0,     0)[t]{$\bullet$}}

\put(126,      1){\makebox(0,    0)[t]{$\bullet$}}

\put(126,      1){\makebox(0,    0)[t]{$\bullet$}}
\put(28,      -1){\line(1,      0){33}}
\put(61,      -1){\line(1,      0){30}}
\put(94,     -1){\line(1,      0){30}}

\put(22,    -20){${-1}$}
\put(40,     -15){$q ^{-1}$}
\put(58,      -20){$q ^{}$}
\put(70,      -15){$q ^{-1}$}

\put(81,     -20){${-1}$}

\put(102,    -15){$q ^{}$}

\put(128,     -20){$q ^{-1}$}

\put(60,    38){\makebox(0,     0)[t]{$\bullet$}}

\put(22,    20){$q ^{-1}$}

\put(28,      -1){\line(1,      1){33}}

\put(60,      -1){\line(0,     1){40}}

\put(43,      33){$q$}

\put(43,      16){$q ^{-1}$}

\put(93,    38){\makebox(0,     0)[t]{$\bullet$}}

\put(93,      -1){\line(0,     1){40}}

\put(76,      33){${-1}$}

\put(73,      16){$q ^{}$}

\put(126,      1){\line(-1,     1){35}}

\put(114,      30){$q ^{}$}

\put(166,        -1)  {. }

\end{picture}$\\ \\ Omitting  Vertex 6  in whole is not  arithmetic   by Row 11-13.
{\rm (d)} There are not any other cases by Lemma \ref {2.13}.
   \subsection  {Quasi-affine  GDDs  over
GDD $7$ of Row $13$ in Table C}\label {sub3.2}
 {\rm (i)}  Omitting  Vertex  1 and adding on   Vertex 5.
It    is empty   by Lemma \ref {2.63}.
    \subsection  {Quasi-affine  GDDs  over
GDD $8$ of Row $13$ in Table C}\label {sub3.2}
{\rm (i)}  Omitting  Vertex  1 and adding on   Vertex 5.
  It is empty  by Lemma \ref {2.63}.
    \subsection  {Quasi-affine  GDDs  over
GDD $9$ of Row $13$}\label {sub3.2}
{\rm (i)} Omitting  Vertex  1 and adding on   Vertex 4.
It    is empty   by Lemma \ref {2.63}.

    {\rm (ii)} Omitting Vertex   1 and adding on   Vertex 5.  It    is empty  by Row 11-13.

 {\rm (iii)}   Adding on   Vertex 3. It is not quasi-affine since it contains a proper subGDD which is not an arithmetic GDD.

   {\rm (iv)}  Adding on  Vertex 2.\\ \\

 $\ \ \ \ \  \ \   
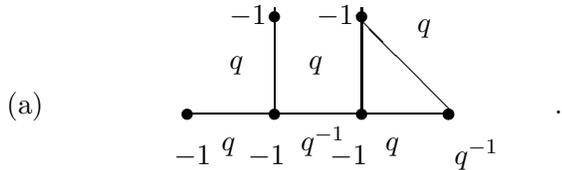
\begin{picture}(100,      15)  \put(-45,      -1){} \put(-45,      -1){ {\rm (a)}}
\put(27,      1){\makebox(0,     0)[t]{$\bullet$}}
\put(60,      1){\makebox(0,      0)[t]{$\bullet$}}
\put(93,     1){\makebox(0,     0)[t]{$\bullet$}}

\put(126,      1){\makebox(0,    0)[t]{$\bullet$}}

\put(126,      1){\makebox(0,    0)[t]{$\bullet$}}
\put(28,      -1){\line(1,      0){33}}
\put(61,      -1){\line(1,      0){30}}
\put(94,     -1){\line(1,      0){30}}

\put(22,    -20){${-1}$}
\put(40,     -15){$q ^{}$}
\put(50,      -20){${-1}$}
\put(70,      -15){$q ^{-1}$}

\put(81,     -20){${-1}$}

\put(102,    -15){$q ^{}$}

\put(128,     -20){$q^{-1}$}

\put(60,    38){\makebox(0,     0)[t]{$\bullet$}}

\put(60,      -1){\line(0,     1){40}}

\put(43,      33){$-1$}

\put(43,      16){$q ^{}$}

\put(93,    38){\makebox(0,     0)[t]{$\bullet$}}

\put(93,      -1){\line(0,     1){40}}

\put(76,      33){$-1$}

\put(73,      16){$q ^{}$}

\put(126,      1){\line(-1,     1){35}}

\put(114,      30){$q ^{}$}

\put(166,        -1)  {.}

\end{picture}$\\ \\
Omitting  Vertex 5  in whole is not an arithmetic GDD   by Lemma \ref {2.2}
 {\rm (i)}.\\ \\

 $\ \ \ \ \  \ \   
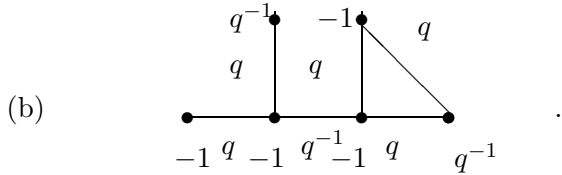
\begin{picture}(100,      15)  \put(-45,      -1){} \put(-45,      -1){ {\rm (b)}}
\put(27,      1){\makebox(0,     0)[t]{$\bullet$}}
\put(60,      1){\makebox(0,      0)[t]{$\bullet$}}
\put(93,     1){\makebox(0,     0)[t]{$\bullet$}}

\put(126,      1){\makebox(0,    0)[t]{$\bullet$}}

\put(126,      1){\makebox(0,    0)[t]{$\bullet$}}
\put(28,      -1){\line(1,      0){33}}
\put(61,      -1){\line(1,      0){30}}
\put(94,     -1){\line(1,      0){30}}

\put(22,    -20){${-1}$}
\put(40,     -15){$q ^{}$}
\put(50,      -20){${-1}$}
\put(70,      -15){$q ^{-1}$}

\put(81,     -20){${-1}$}

\put(102,    -15){$q ^{}$}

\put(128,     -20){$q^{-1}$}

\put(60,    38){\makebox(0,     0)[t]{$\bullet$}}

\put(60,      -1){\line(0,     1){40}}

\put(43,      33){$q ^{-1}$}

\put(43,      16){$q ^{}$}

\put(93,    38){\makebox(0,     0)[t]{$\bullet$}}

\put(93,      -1){\line(0,     1){40}}

\put(76,      33){$-1$}

\put(73,      16){$q ^{}$}

\put(126,      1){\line(-1,     1){35}}

\put(114,      30){$q ^{}$}

\put(166,        -1)  {. }

\end{picture}$\\ \\  Omitting  Vertex 1  in whole is not  arithmetic by Row 11-13.
\\ \\

 $\ \ \ \ \  \ \   
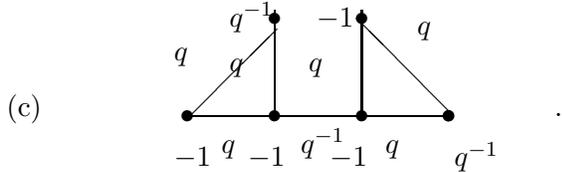
\begin{picture}(100,      15)  \put(-45,      -1){} \put(-45,      -1){ {\rm (c)}}
\put(27,      1){\makebox(0,     0)[t]{$\bullet$}}
\put(60,      1){\makebox(0,      0)[t]{$\bullet$}}
\put(93,     1){\makebox(0,     0)[t]{$\bullet$}}

\put(126,      1){\makebox(0,    0)[t]{$\bullet$}}

\put(126,      1){\makebox(0,    0)[t]{$\bullet$}}
\put(28,      -1){\line(1,      0){33}}
\put(61,      -1){\line(1,      0){30}}
\put(94,     -1){\line(1,      0){30}}

\put(22,    -20){${-1}$}
\put(40,     -15){$q ^{}$}
\put(50,      -20){${-1}$}
\put(70,      -15){$q ^{-1}$}

\put(81,     -20){${-1}$}

\put(102,    -15){$q ^{}$}

\put(128,     -20){$q ^{-1}$}

\put(60,    38){\makebox(0,     0)[t]{$\bullet$}}

\put(22,    20){$q ^{}$}

\put(28,      -1){\line(1,      1){33}}

\put(60,      -1){\line(0,     1){40}}

\put(43,      33){$q ^{-1}$}

\put(43,      16){$q ^{}$}

\put(93,    38){\makebox(0,     0)[t]{$\bullet$}}

\put(93,      -1){\line(0,     1){40}}

\put(76,      33){${-1}$}

\put(73,      16){$q ^{}$}

\put(126,      1){\line(-1,     1){35}}

\put(114,      30){$q ^{}$}

\put(166,        -1)  {.  }

\end{picture}$\\ \\
Omitting  Vertex 5  in whole is not arithmetic by Row 11-13.

{\rm (d)} There are not any other cases by Lemma \ref {2.13}.
    \subsection  {Quasi-affine  GDDs  over
GDD $10$ of Row $13$ in Table C}\label {sub3.2}
{\rm (i)}  Omitting Vertex 1 and adding on  Vertex  4.
 It    is empty  by  Lemma \ref {2.2}
 {\rm (i)}.

     {\rm (ii)} Omitting Vertex 1 and adding on  Vertex  5.\\ \\

 $\ \ \ \ \  \ \   
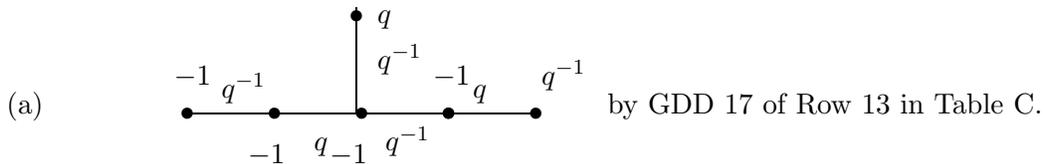
\begin{picture}(100,      15)  \put(-45,      -1){}

\put(-45,      -1){ {\rm (a)}}
\put(27,      1){\makebox(0,     0)[t]{$\bullet$}}
\put(60,      1){\makebox(0,      0)[t]{$\bullet$}}
\put(93,     1){\makebox(0,     0)[t]{$\bullet$}}

\put(126,      1){\makebox(0,    0)[t]{$\bullet$}}

\put(159,      1){\makebox(0,    0)[t]{$\bullet$}}

\put(126,      1){\makebox(0,    0)[t]{$\bullet$}}
\put(28,      -1){\line(1,      0){33}}
\put(61,      -1){\line(1,      0){30}}
\put(94,     -1){\line(1,      0){30}}
\put(126,      -1){\line(1,      0){33}}

\put(22,     10){${-1}$}
\put(40,      5){$q^{-1}$}
\put(50,      -20){${-1}$}
\put(75,      -15){$q ^{}$}

\put(81,      -20){${-1}$}

\put(102,    -15){$q ^{-1}$}

\put(120,      10){${-1}$}

\put(135,     5){$q ^{}$}

\put(161,      10){$q^{-1}$}

\put(91,    38){\makebox(0,     0)[t]{$\bullet$}}

\put(91,      -1){\line(0,     1){40}}

\put(99,      33){$q ^{}$}

\put(99,      16){$q ^{-1}$}

\put(186,        -1)  {by GDD $17$ of Row $13$  in Table C. }
\end{picture}$\\ \\
 It is repeated.\\ \\

 $\ \ \ \ \  \ \   
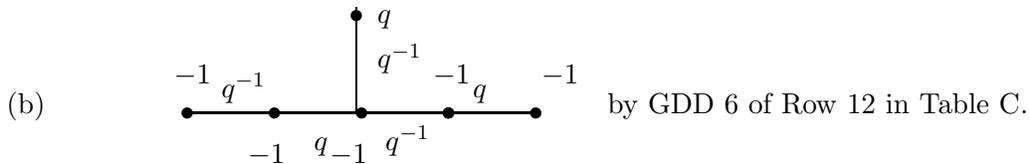
\begin{picture}(100,      15)  \put(-45,      -1){}

\put(-45,      -1){ {\rm (b)}}
\put(27,      1){\makebox(0,     0)[t]{$\bullet$}}
\put(60,      1){\makebox(0,      0)[t]{$\bullet$}}
\put(93,     1){\makebox(0,     0)[t]{$\bullet$}}

\put(126,      1){\makebox(0,    0)[t]{$\bullet$}}

\put(159,      1){\makebox(0,    0)[t]{$\bullet$}}

\put(126,      1){\makebox(0,    0)[t]{$\bullet$}}
\put(28,      -1){\line(1,      0){33}}
\put(61,      -1){\line(1,      0){30}}
\put(94,     -1){\line(1,      0){30}}
\put(126,      -1){\line(1,      0){33}}

\put(22,     10){${-1}$}
\put(40,      5){$q^{-1}$}
\put(50,      -20){${-1}$}
\put(75,      -15){$q ^{}$}

\put(81,      -20){${-1}$}

\put(102,    -15){$q ^{-1}$}

\put(120,      10){${-1}$}

\put(135,     5){$q ^{}$}

\put(161,      10){${-1}$}

\put(91,    38){\makebox(0,     0)[t]{$\bullet$}}

\put(91,      -1){\line(0,     1){40}}

\put(99,      33){$q ^{}$}

\put(99,      16){$q ^{-1}$}

\put(186,        -1)  {by GDD $6$ of Row $12$  in Table C. }
\end{picture}$\\ \\
Omitting  Vertex 4  in whole is a simple chain.   It is quasi-affine.
{\rm (c)} There are not any other cases by Lemma \ref {2.13}.

 {\rm (iii)}   Adding on   Vertex 3. It is not quasi-affine since it contains a proper subGDD which is not an arithmetic GDD.

   {\rm (iv)}  Adding on  Vertex 2.\\ \\

 $\ \ \ \ \  \ \   
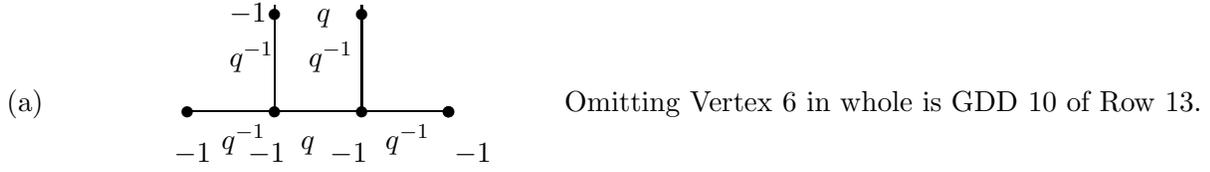
\begin{picture}(100,      15)  \put(-45,      -1){}

\put(-45,      -1){ {\rm (a)}}
\put(27,      1){\makebox(0,     0)[t]{$\bullet$}}
\put(60,      1){\makebox(0,      0)[t]{$\bullet$}}
\put(93,     1){\makebox(0,     0)[t]{$\bullet$}}

\put(126,      1){\makebox(0,    0)[t]{$\bullet$}}

\put(126,      1){\makebox(0,    0)[t]{$\bullet$}}
\put(28,      -1){\line(1,      0){33}}
\put(61,      -1){\line(1,      0){30}}
\put(94,     -1){\line(1,      0){30}}

\put(22,     -20){$-1$}
\put(40,     -15){$q ^{-1}$}
\put(50,      -20){$-1$}
\put(70,      -15){$q ^{}$}

\put(81,      -20){$-1$}

\put(102,     -15){$q ^{-1}$}

\put(128,      -20){${-1}$}

\put(60,    38){\makebox(0,     0)[t]{$\bullet$}}

\put(60,      -1){\line(0,     1){40}}

\put(43,      33){$-1$}

\put(43,      16){$q ^{-1}$}

\put(93,    38){\makebox(0,     0)[t]{$\bullet$}}

\put(93,      -1){\line(0,     1){40}}

\put(76,      33){$q$}

\put(73,      16){$q ^{-1}$}

\put(166,        -1)  { Omitting  Vertex 6 in whole  is GDD $10$ of Row $13$. }
\end{picture}$\\ \\
 Omitting  Vertex 5  in whole is not an arithmetic GDD   by Lemma \ref {2.2}
 {\rm (i)}.\\ \\

 $\ \ \ \ \  \ \   
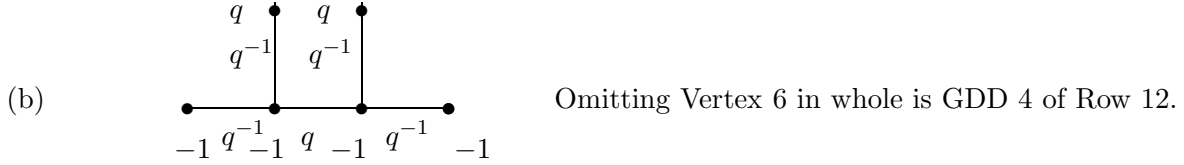
\begin{picture}(100,      15)  \put(-45,      -1){}

\put(-45,      -1){ {\rm (b)}}
\put(27,      1){\makebox(0,     0)[t]{$\bullet$}}
\put(60,      1){\makebox(0,      0)[t]{$\bullet$}}
\put(93,     1){\makebox(0,     0)[t]{$\bullet$}}

\put(126,      1){\makebox(0,    0)[t]{$\bullet$}}

\put(126,      1){\makebox(0,    0)[t]{$\bullet$}}
\put(28,      -1){\line(1,      0){33}}
\put(61,      -1){\line(1,      0){30}}
\put(94,     -1){\line(1,      0){30}}

\put(22,     -20){$-1$}
\put(40,      -15){$q ^{-1}$}
\put(50,      -20){$-1$}
\put(70,      -15){$q ^{}$}

\put(81,      -20){$-1$}

\put(102,     -15){$q ^{-1}$}

\put(128,      -20){${-1}$}

\put(60,    38){\makebox(0,     0)[t]{$\bullet$}}

\put(60,      -1){\line(0,     1){40}}

\put(43,      33){$q$}

\put(43,      16){$q ^{-1}$}

\put(93,    38){\makebox(0,     0)[t]{$\bullet$}}

\put(93,      -1){\line(0,     1){40}}

\put(76,      33){$q$}

\put(73,      16){$q ^{-1}$}

\put(166,        -1)  {Omitting  Vertex 6 in whole  is GDD $4$ of Row $12$. }
\end{picture}$\\ \\
 Omitting Vertex 5  in whole is  GDD $10$ of Row $13$. Omitting  Vertex 6  in whole is  GDD $4$ of Row $12$.
  It is quasi-affine.\\ \\

$\ \ \ \ \  \ \   
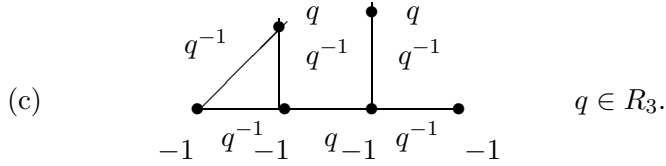
\begin{picture}(130,      15)

\put(-45,      -1){{\rm (c)} }
\put(27,      1){\makebox(0,     0)[t]{$\bullet$}}
\put(60,      1){\makebox(0,      0)[t]{$\bullet$}}
\put(93,     1){\makebox(0,     0)[t]{$\bullet$}}

\put(126,      1){\makebox(0,    0)[t]{$\bullet$}}

\put(28,      -1){\line(1,      0){33}}
\put(61,      -1){\line(1,      0){30}}
\put(94,     -1){\line(1,      0){30}}

\put(12,     -20){${-1}$}
\put(36,      -15){$q ^{-1}$}
\put(48,      -20){$-1 ^{}$}
\put(75,      -15){$q ^{}$}

\put(81,      -20){$-1 ^{}$}

\put(102,     -15){$q ^{-1}$}

\put(128,      -20){${-1}$}

\put(58,    32){\makebox(0,     0)[t]{$\bullet$}}

\put(58,      -1){\line(0,     1){34}}

\put(22,    20){$q ^{-1}$}

\put(68,      33){$q ^{}$}

\put(68,      16){$q ^{-1}$}

\put(28,      -1){\line(1,      1){33}}

\put(93,    38){\makebox(0,     0)[t]{$\bullet$}}

\put(93,      -1){\line(0,     1){40}}

\put(106,      33){$q$}

\put(103,      16){$q ^{-1}$}

\put(166,        -1)  {  $q \in R_3$. }
\end{picture}$\\ \\ Omitting  Vertex 5 in whole  is GDD $9$ of Row $13$. Omitting  Vertex 6 in whole  is GDD $3$ of Row $12$. It is quasi-affine.
{\rm (d)} There are not any other cases by Lemma \ref {2.13}.

    \subsection  {Quasi-affine  GDDs  over
GDD $11$ of Row $13$ in Table C}\label {sub3.2}
{\rm (i)} Omitting Vertex   1 and adding on  Vertex  5.
 It    is empty   by Lemma \ref {2.63}.
    \subsection  {Quasi-affine  GDDs  over
GDD $12$ of Row $13$ in Table C}\label {sub3.2}
 {\rm (i)} Omitting  Vertex  1 and adding on   Vertex 4.
 It    is empty  by  Lemma \ref {2.63}.

   {\rm (ii)}  Omitting Vertex 1 and adding on  Vertex  5.   It    is empty  by  Lemma \ref {2.63}.

 {\rm (iii)}   Adding on   Vertex 3. It is not quasi-affine since it contains a proper subGDD which is not an arithmetic GDD.

   {\rm (iv)} Adding on  Vertex 2.\\ \\

 $\ \ \ \ \  \ \   
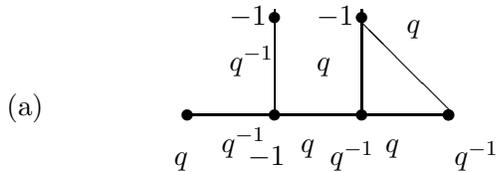
\begin{picture}(100,      15)  \put(-45,      -1){} \put(-45,      -1){ {\rm (a)}}
\put(27,      1){\makebox(0,     0)[t]{$\bullet$}}
\put(60,      1){\makebox(0,      0)[t]{$\bullet$}}
\put(93,     1){\makebox(0,     0)[t]{$\bullet$}}

\put(126,      1){\makebox(0,    0)[t]{$\bullet$}}

\put(126,      1){\makebox(0,    0)[t]{$\bullet$}}
\put(28,      -1){\line(1,      0){33}}
\put(61,      -1){\line(1,      0){30}}
\put(94,     -1){\line(1,      0){30}}

\put(22,    -20){$q ^{}$}
\put(40,     -15){$q ^{-1}$}
\put(50,      -20){${-1}$}
\put(70,      -15){$q ^{}$}

\put(81,     -20){$q ^{-1}$}

\put(102,    -15){$q ^{}$}

\put(128,     -20){$q^{-1}$}

\put(60,    38){\makebox(0,     0)[t]{$\bullet$}}

\put(60,      -1){\line(0,     1){40}}

\put(43,      33){$-1$}

\put(43,      16){$q ^{-1}$}

\put(93,    38){\makebox(0,     0)[t]{$\bullet$}}

\put(93,      -1){\line(0,     1){40}}

\put(76,      33){${-1}$}

\put(76,      16){$q ^{}$}

\put(126,      1){\line(-1,     1){35}}

\put(110,      30){$q ^{}$}

\put(166,        -1)  {
.}

\end{picture}$\\ \\ Omitting  Vertex 1 in whole  is not  arithmetic  by Row 11-Row 13.\\ \\

 $\ \ \ \ \  \ \   
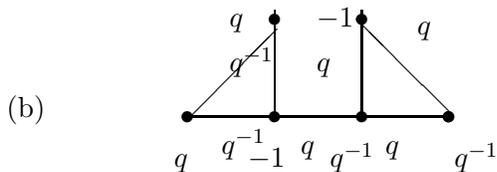
\begin{picture}(100,      15)  \put(-45,      -1){} \put(-45,      -1){ {\rm (b)}}
\put(27,      1){\makebox(0,     0)[t]{$\bullet$}}
\put(60,      1){\makebox(0,      0)[t]{$\bullet$}}
\put(93,     1){\makebox(0,     0)[t]{$\bullet$}}

\put(126,      1){\makebox(0,    0)[t]{$\bullet$}}

\put(126,      1){\makebox(0,    0)[t]{$\bullet$}}
\put(28,      -1){\line(1,      0){33}}
\put(61,      -1){\line(1,      0){30}}
\put(94,     -1){\line(1,      0){30}}

\put(22,    -20){$q ^{}$}
\put(40,     -15){$q ^{-1}$}
\put(50,      -20){${-1}$}
\put(70,      -15){$q ^{}$}

\put(81,     -20){$q ^{-1}$}

\put(102,    -15){$q ^{}$}

\put(128,     -20){$q ^{-1}$}

\put(60,    38){\makebox(0,     0)[t]{$\bullet$}}

\put(22,    20){$$}

\put(28,      -1){\line(1,      1){33}}

\put(60,      -1){\line(0,     1){40}}

\put(43,      33){$q$}

\put(43,      16){$q ^{-1}$}

\put(93,    38){\makebox(0,     0)[t]{$\bullet$}}

\put(93,      -1){\line(0,     1){40}}

\put(76,      33){${-1}$}

\put(76,      16){$q ^{}$}

\put(126,      1){\line(-1,     1){35}}

\put(114,      30){$q ^{}$}

\put(166,        -1)  {.  }

\end{picture}$\\ \\
Omitting  Vertex  Vertex 4  in whole is not  arithmetic  by Lemma \ref {2.2}
 {\rm (i)}.
 {\rm (c)} There are not any other cases by  Lemma \ref {2.13}.
    \subsection  {Quasi-affine  GDDs  over
GDD $13$ of Row $13$ in Table C}\label {sub3.2}
 {\rm (i)} Omitting Vertex 1 and adding on   Vertex 4.
 It    is empty  by Lemma \ref {2.63}.

 {\rm (ii)} Omitting Vertex 1 and adding on   Vertex 5.\\ \\

 $\ \ \ \ \  \ \   
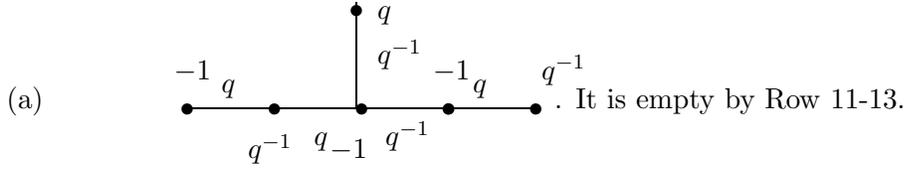
\begin{picture}(100,      15)  \put(-45,      -1){}

\put(-45,      -1){ {\rm (a)}}
\put(27,      1){\makebox(0,     0)[t]{$\bullet$}}
\put(60,      1){\makebox(0,      0)[t]{$\bullet$}}
\put(93,     1){\makebox(0,     0)[t]{$\bullet$}}

\put(126,      1){\makebox(0,    0)[t]{$\bullet$}}

\put(159,      1){\makebox(0,    0)[t]{$\bullet$}}

\put(126,      1){\makebox(0,    0)[t]{$\bullet$}}
\put(28,      -1){\line(1,      0){33}}
\put(61,      -1){\line(1,      0){30}}
\put(94,     -1){\line(1,      0){30}}
\put(126,      -1){\line(1,      0){33}}

\put(22,     10){${-1}$}
\put(40,      5){$q^{}$}
\put(50,      -20){$q ^{-1}$}
\put(75,      -15){$q ^{}$}

\put(81,      -20){${-1}$}

\put(102,    -15){$q ^{-1}$}

\put(120,      10){${-1}$}

\put(135,     5){$q ^{}$}

\put(161,      10){$q ^{-1}$}

\put(91,    38){\makebox(0,     0)[t]{$\bullet$}}

\put(91,      -1){\line(0,     1){40}}

\put(99,      33){$q ^{}$}

\put(99,      16){$q ^{-1}$}

\put(166,        -1)  {. It is empty   by Row 11-13.}
\end{picture}$\\ \\ \\

 $\ \ \ \ \  \ \   
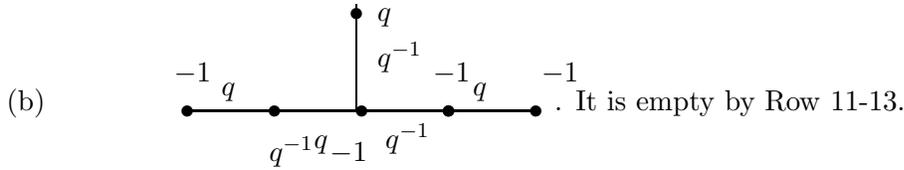
\begin{picture}(100,      15)  \put(-45,      -1){}

\put(-45,      -1){ {\rm (b)}}
\put(27,      1){\makebox(0,     0)[t]{$\bullet$}}
\put(60,      1){\makebox(0,      0)[t]{$\bullet$}}
\put(93,     1){\makebox(0,     0)[t]{$\bullet$}}

\put(126,      1){\makebox(0,    0)[t]{$\bullet$}}

\put(159,      1){\makebox(0,    0)[t]{$\bullet$}}

\put(126,      1){\makebox(0,    0)[t]{$\bullet$}}
\put(28,      -1){\line(1,      0){33}}
\put(61,      -1){\line(1,      0){30}}
\put(94,     -1){\line(1,      0){30}}
\put(126,      -1){\line(1,      0){33}}

\put(22,     10){${-1}$}
\put(40,      5){$q^{}$}
\put(58,      -20){$q ^{-1}$}
\put(75,      -15){$q ^{}$}

\put(81,      -20){${-1}$}

\put(102,    -15){$q ^{-1}$}

\put(120,      10){${-1}$}

\put(135,     5){$q ^{}$}

\put(161,      10){${-1}$}

\put(91,    38){\makebox(0,     0)[t]{$\bullet$}}

\put(91,      -1){\line(0,     1){40}}

\put(99,      33){$q ^{}$}

\put(99,      16){$q ^{-1}$}

\put(166,        -1)  {. It is empty   by Row 11-13.}

\end{picture}$\\ \\
{\rm (c)} There are not any other cases by Lemma \ref {2.13}.

 {\rm (iii)}    Adding on   Vertex 3. It is not quasi-affine since it contains a proper subGDD which is not arithmetic.

 {\rm (iv)} Adding on  Vertex 2.\\ \\

 $\ \ \ \ \  \ \   
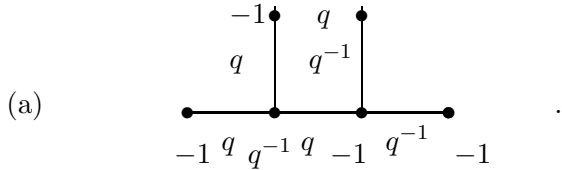
\begin{picture}(100,      15)  \put(-45,      -1){}

\put(-45,      -1){ {\rm (a)}}
\put(27,      1){\makebox(0,     0)[t]{$\bullet$}}
\put(60,      1){\makebox(0,      0)[t]{$\bullet$}}
\put(93,     1){\makebox(0,     0)[t]{$\bullet$}}

\put(126,      1){\makebox(0,    0)[t]{$\bullet$}}

\put(126,      1){\makebox(0,    0)[t]{$\bullet$}}
\put(28,      -1){\line(1,      0){33}}
\put(61,      -1){\line(1,      0){30}}
\put(94,     -1){\line(1,      0){30}}

\put(22,    -20){$-1$}
\put(40,      -15){$q ^{}$}
\put(50,      -20){$q ^{-1}$}
\put(70,      -15){$q ^{}$}

\put(81,      -20){$-1$}

\put(102,     -15){$q ^{-1}$}

\put(128,      -20){${-1}$}

\put(60,    38){\makebox(0,     0)[t]{$\bullet$}}

\put(60,      -1){\line(0,     1){40}}

\put(43,      33){$-1$}

\put(43,      16){$q ^{}$}

\put(93,    38){\makebox(0,     0)[t]{$\bullet$}}

\put(93,      -1){\line(0,     1){40}}

\put(76,      33){$q$}

\put(73,      16){$q ^{-1}$}

\put(166,        -1)  {. }
\end{picture}$\\ \\ Omitting   Vertex 1  in whole is   GDD $13$ of Row $13$  in Table C.
Omitting   Vertex 5  in whole is not arithmetic by  Lemma \ref {2.2}
 {\rm (i)}.\\ \\

 $\ \ \ \ \  \ \   
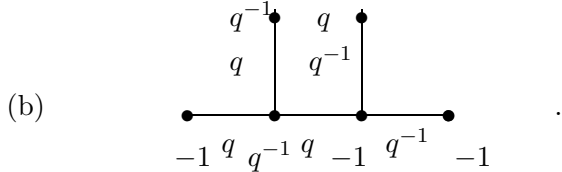
\begin{picture}(100,      15)  \put(-45,      -1){}

\put(-45,      -1){ {\rm (b)}}
\put(27,      1){\makebox(0,     0)[t]{$\bullet$}}
\put(60,      1){\makebox(0,      0)[t]{$\bullet$}}
\put(93,     1){\makebox(0,     0)[t]{$\bullet$}}

\put(126,      1){\makebox(0,    0)[t]{$\bullet$}}

\put(126,      1){\makebox(0,    0)[t]{$\bullet$}}
\put(28,      -1){\line(1,      0){33}}
\put(61,      -1){\line(1,      0){30}}
\put(94,     -1){\line(1,      0){30}}

\put(22,     -20){$-1$}
\put(40,      -15){$q ^{}$}
\put(50,      -20){$q ^{-1}$}
\put(70,      -15){$q ^{}$}

\put(81,      -20){$-1$}

\put(102,     -15){$q ^{-1}$}

\put(128,      -20){${-1}$}

\put(60,    38){\makebox(0,     0)[t]{$\bullet$}}

\put(60,      -1){\line(0,     1){40}}

\put(43,      33){$q ^{-1}$}

\put(43,      16){$q ^{}$}

\put(93,    38){\makebox(0,     0)[t]{$\bullet$}}

\put(93,      -1){\line(0,     1){40}}

\put(76,      33){$q$}

\put(73,      16){$q ^{-1}$}

\put(166,        -1)  {. }
\end{picture}$\\ \\ Omitting   Vertex  Vertex 1  in whole is  not  arithmetic   by Row 11-13.\\ \\

$\ \ \ \ \  \ \   
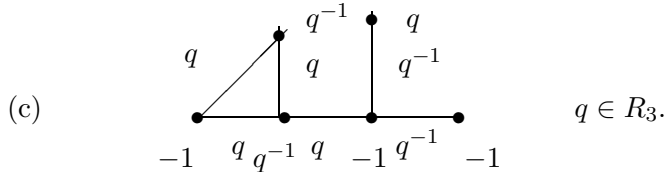
\begin{picture}(130,      15)

\put(-45,      -1){{\rm (c)} }
\put(27,      1){\makebox(0,     0)[t]{$\bullet$}}
\put(60,      1){\makebox(0,      0)[t]{$\bullet$}}
\put(93,     1){\makebox(0,     0)[t]{$\bullet$}}

\put(126,      1){\makebox(0,    0)[t]{$\bullet$}}

\put(28,      -1){\line(1,      0){33}}
\put(61,      -1){\line(1,      0){30}}
\put(94,     -1){\line(1,      0){30}}

\put(12,     -20){${-1}$}
\put(40,      -15){$q ^{}$}
\put(48,      -20){$ q^{-1}$}
\put(70,      -15){$q ^{}$}

\put(85,      -20){$-1 ^{}$}

\put(102,     -15){$q ^{-1}$}

\put(128,      -20){${-1}$}

\put(58,    32){\makebox(0,     0)[t]{$\bullet$}}

\put(58,      -1){\line(0,     1){34}}

\put(22,    20){$q ^{}$}

\put(68,      33){$q ^{-1}$}

\put(68,      16){$q ^{}$}

\put(28,      -1){\line(1,      1){33}}

\put(93,    38){\makebox(0,     0)[t]{$\bullet$}}

\put(93,      -1){\line(0,     1){40}}

\put(106,      33){$q$}

\put(103,      16){$q ^{-1}$}

\put(166,        -1)  {  $q \in R_3$. }
\end{picture}$\\ \\ Omitting  Vertex 5 in whole  is not  arithmetic  by Row 11-13.

{\rm (d)} There are not any other cases by Lemma \ref {2.13}.
    \subsection  {Quasi-affine  GDDs  over
GDD $14$ of Row $13$ in Table C}\label {sub3.2}
{\rm (i)} Omitting  Vertex  1 and adding on   Vertex 4.
It    is empty  by  Lemma \ref {2.2}
 {\rm (i)}.

    {\rm (ii)}  Omitting Vertex 1 and adding on   Vertex 5.\\ \\

    $\ \ \ \ \  \ \   \begin{picture}(100,      15)  \put(-45,      -1){}

\put(-45,      -1){ {\rm (a)}}
\put(27,      1){\makebox(0,     0)[t]{$\bullet$}}
\put(60,      1){\makebox(0,      0)[t]{$\bullet$}}
\put(93,     1){\makebox(0,     0)[t]{$\bullet$}}

\put(126,      1){\makebox(0,    0)[t]{$\bullet$}}

\put(159,      1){\makebox(0,    0)[t]{$\bullet$}}

\put(126,      1){\makebox(0,    0)[t]{$\bullet$}}
\put(28,      -1){\line(1,      0){33}}
\put(61,      -1){\line(1,      0){30}}
\put(94,     -1){\line(1,      0){30}}
\put(126,      -1){\line(1,      0){33}}

\put(22,     10){$q ^{-1}$}
\put(40,      5){$q^{}$}
\put(50,      -20){${-1}$}
\put(70,      -15){$q ^{-1}$}

\put(91,      -20){$q ^{}$}

\put(102,    -15){$q ^{-1}$}

\put(120,      10){${-1}$}

\put(135,     5){$q ^{}$}

\put(161,      10){$q ^{-1}$}

\put(91,    38){\makebox(0,     0)[t]{$\bullet$}}

\put(91,      -1){\line(0,     1){40}}

\put(99,      33){$q ^{}$}

\put(99,      16){$q ^{-1}$}

\put(186,        -1)  {by GDD $14$ of Row $13$  in Table C. }
\end{picture}$\\ \\
 It is quasi-affine    by Lemma \ref {2.63}.\\ \\

    $\ \ \ \ \  \ \   \begin{picture}(100,      15)  \put(-45,      -1){}

\put(-45,      -1){ {\rm (b)}}
\put(27,      1){\makebox(0,     0)[t]{$\bullet$}}
\put(60,      1){\makebox(0,      0)[t]{$\bullet$}}
\put(93,     1){\makebox(0,     0)[t]{$\bullet$}}

\put(126,      1){\makebox(0,    0)[t]{$\bullet$}}

\put(159,      1){\makebox(0,    0)[t]{$\bullet$}}

\put(126,      1){\makebox(0,    0)[t]{$\bullet$}}
\put(28,      -1){\line(1,      0){33}}
\put(61,      -1){\line(1,      0){30}}
\put(94,     -1){\line(1,      0){30}}
\put(126,      -1){\line(1,      0){33}}

\put(22,     10){$q ^{-1}$}
\put(40,      5){$q^{}$}
\put(50,      -20){${-1}$}
\put(70,      -15){$q ^{-1}$}

\put(91,      -20){$q ^{}$}

\put(102,    -15){$q ^{-1}$}

\put(120,      10){${-1}$}

\put(135,     5){$q ^{}$}

\put(161,      10){${-1}$}

\put(91,    38){\makebox(0,     0)[t]{$\bullet$}}

\put(91,      -1){\line(0,     1){40}}

\put(99,      33){$q ^{}$}

\put(99,      16){$q ^{-1}$}

\put(186,        -1)  {by GDD $8$ of Row $12$  in Table C. }
\end{picture}$\\ \\
 It is quasi-affine    by Lemma \ref {2.63}.
{\rm (c)} There are not any other cases by Lemma \ref {2.13}.

 {\rm (iii)}   Adding on   Vertex 3. It is not quasi-affine since it contains a proper subGDD which is not arithmetic.

 {\rm (iv)} Adding on  Vertex 2.\\ \\

 $\ \ \ \ \  \ \   
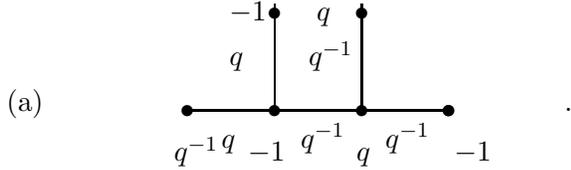
\begin{picture}(100,      15)  \put(-45,      -1){}

\put(-45,      -1){ {\rm (a)}}
\put(27,      1){\makebox(0,     0)[t]{$\bullet$}}
\put(60,      1){\makebox(0,      0)[t]{$\bullet$}}
\put(93,     1){\makebox(0,     0)[t]{$\bullet$}}

\put(126,      1){\makebox(0,    0)[t]{$\bullet$}}

\put(126,      1){\makebox(0,    0)[t]{$\bullet$}}
\put(28,      -1){\line(1,      0){33}}
\put(61,      -1){\line(1,      0){30}}
\put(94,     -1){\line(1,      0){30}}

\put(22,     -20){$q ^{-1}$}
\put(40,      -15){$q ^{}$}
\put(50,      -20){$-1$}
\put(70,      -15){$q ^{-1}$}

\put(91,      -20){$q$}

\put(102,     -15){$q ^{-1}$}

\put(128,      -20){${-1}$}

\put(60,    38){\makebox(0,     0)[t]{$\bullet$}}

\put(60,      -1){\line(0,     1){40}}

\put(43,      33){${-1}$}

\put(43,      16){$q ^{}$}

\put(93,    38){\makebox(0,     0)[t]{$\bullet$}}

\put(93,      -1){\line(0,     1){40}}

\put(76,      33){$q$}

\put(73,      16){$q ^{-1}$}

\put(166,        -1)  { . }
\end{picture}$\\ \\ Omitting  Vertex 1  in whole is GDD $8$ of Row $12$  in Table C.
Omitting  Vertex 6 in whole  is not  arithmetic  by Row 11- 13.\\ \\

$\ \ \ \ \  \ \   
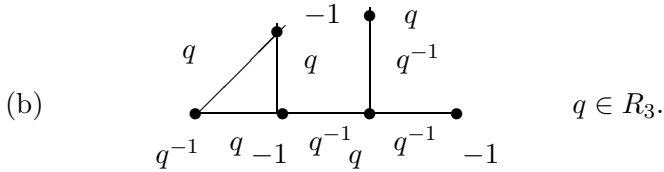
\begin{picture}(130,      15)

\put(-45,      -1){{\rm (b)} }
\put(27,      1){\makebox(0,     0)[t]{$\bullet$}}
\put(60,      1){\makebox(0,      0)[t]{$\bullet$}}
\put(93,     1){\makebox(0,     0)[t]{$\bullet$}}

\put(126,      1){\makebox(0,    0)[t]{$\bullet$}}

\put(28,      -1){\line(1,      0){33}}
\put(61,      -1){\line(1,      0){30}}
\put(94,     -1){\line(1,      0){30}}

\put(12,     -20){$q ^{-1}$}
\put(40,      -15){$q ^{}$}
\put(48,      -20){$-1 ^{}$}
\put(70,      -15){$q ^{-1}$}

\put(85,      -20){$q ^{}$}

\put(102,     -15){$q ^{-1}$}

\put(128,      -20){${-1}$}

\put(58,    32){\makebox(0,     0)[t]{$\bullet$}}

\put(58,      -1){\line(0,     1){34}}

\put(22,    20){$q ^{}$}

\put(68,      33){$-1 ^{}$}

\put(68,      16){$q ^{}$}

\put(28,      -1){\line(1,      1){33}}

\put(93,    38){\makebox(0,     0)[t]{$\bullet$}}

\put(93,      -1){\line(0,     1){40}}

\put(106,      33){$q$}

\put(103,      16){$q ^{-1}$}

\put(166,        -1)  {  $q \in R_3$. }
\end{picture}$\\ \\ Omitting  Vertex 5 in whole  is GDD $6$ of Row $13$. Omitting  Vertex 6 in whole  is GDD $3$ of Row $12$. It is quasi-affine.

 {\rm (c)} There are not any other cases by Lemma \ref {2.13}.
    \subsection  {Quasi-affine  GDDs  over
GDD $15$ of Row $13$}\label {sub3.2}
 {\rm (i)}  Omitting  Vertex  1 and adding on   Vertex 4.
 It    is empty  by Lemma \ref {2.2}
 {\rm (i)}.

    {\rm (ii)} Omitting   Vertex 1 and adding on   Vertex 5.\\ \\

 $\ \ \ \ \  \ \   \begin{picture}(100,      15)  \put(-45,      -1){}

\put(-45,      -1){ {\rm (a)}}
\put(27,      1){\makebox(0,     0)[t]{$\bullet$}}
\put(60,      1){\makebox(0,      0)[t]{$\bullet$}}
\put(93,     1){\makebox(0,     0)[t]{$\bullet$}}

\put(126,      1){\makebox(0,    0)[t]{$\bullet$}}

\put(159,      1){\makebox(0,    0)[t]{$\bullet$}}

\put(126,      1){\makebox(0,    0)[t]{$\bullet$}}
\put(28,      -1){\line(1,      0){33}}
\put(61,      -1){\line(1,      0){30}}
\put(94,     -1){\line(1,      0){30}}
\put(126,      -1){\line(1,      0){33}}

\put(22,     10){${-1}$}
\put(40,      5){$q^{-1}$}
\put(58,      -20){${-1}$}
\put(75,      -15){$q ^{}$}

\put(91,      -20){$q ^{-1}$}

\put(112,    -15){$q ^{}$}

\put(120,      10){${-1}$}

\put(135,     5){$q ^{-1}$}

\put(161,      10){$q ^{}$}

\put(91,    38){\makebox(0,     0)[t]{$\bullet$}}

\put(91,      -1){\line(0,     1){40}}

\put(99,      33){$q ^{}$}

\put(99,      16){$q ^{-1}$}

\put(176,        -1)  {by GDD $5$ of Row $12$  in Table C. }
\end{picture}$\\ \\
 It is quasi-affine    by Lemma \ref {2.63}.\\ \\

 $\ \ \ \ \  \ \   \begin{picture}(100,      15)  \put(-45,      -1){}

\put(-45,      -1){ {\rm (b)}}
\put(27,      1){\makebox(0,     0)[t]{$\bullet$}}
\put(60,      1){\makebox(0,      0)[t]{$\bullet$}}
\put(93,     1){\makebox(0,     0)[t]{$\bullet$}}

\put(126,      1){\makebox(0,    0)[t]{$\bullet$}}

\put(159,      1){\makebox(0,    0)[t]{$\bullet$}}

\put(126,      1){\makebox(0,    0)[t]{$\bullet$}}
\put(28,      -1){\line(1,      0){33}}
\put(61,      -1){\line(1,      0){30}}
\put(94,     -1){\line(1,      0){30}}
\put(126,      -1){\line(1,      0){33}}

\put(22,     10){${-1}$}
\put(40,      5){$q^{-1}$}
\put(58,      -20){${-1}$}
\put(75,      -15){$q ^{}$}

\put(91,      -20){$q ^{-1}$}

\put(112,    -15){$q ^{}$}

\put(120,      10){${-1}$}

\put(135,     5){$q ^{-1}$}

\put(161,      10){$-1 ^{}$}

\put(91,    38){\makebox(0,     0)[t]{$\bullet$}}

\put(91,      -1){\line(0,     1){40}}

\put(99,      33){$q ^{}$}

\put(99,      16){$q ^{-1}$}

\put(176,        -1)  {by GDD $15$ of Row $13$  in Table C. }
\end{picture}$\\ \\
 It is quasi-affine    by Lemma \ref {2.63}.
{\rm (c)} There are not any other cases by Lemma \ref {2.13}.

 {\rm (iii)}   Adding on   Vertex 3. It is not quasi-affine since it contains a proper subGDD which is not arithmetic.

  {\rm (iv)} Adding on  Vertex 2. Omitting  Vertex 6  in whole is  not  arithmetic  by Lemma \ref {2.13}.
    \subsection  {Quasi-affine  GDDs  over
GDD $16$ of Row $13$ in Table C}\label {sub3.2}
{\rm (i)}  Omitting   Vertex 1 and adding on   Vertex 4.
   It    is empty  by Lemma \ref {2.13}.

  {\rm (ii)}  Omitting Vertex 1 and adding on   Vertex 5.
   It    is empty  by Row 11-13.

 {\rm (iii)}  Adding on  Vertex  3. It is not quasi-affine since it contains a proper subGDD which is not arithmetic.

  {\rm (iv)} Adding on  Vertex 2.    Omitting  Vertex  Vertex 6  in whole is  not arithmetic  by Lemma \ref {2.13}.
    \subsection  {Quasi-affine  GDDs  over
GDD $17$ of Row $13$ in Table C}\label {sub3.2}
 {\rm (i)}  Omitting   Vertex 1 and adding on  Vertex  4.
It    is empty   by Lemma \ref {2.2}
 {\rm (i)}.

    {\rm (ii)} Omitting  Vertex  1 and adding on   Vertex 5.\\ \\

   $\ \ \ \ \  \ \   \begin{picture}(100,      15)  \put(-45,      -1){}

\put(-45,      -1){ {\rm (a)}}
\put(27,      1){\makebox(0,     0)[t]{$\bullet$}}
\put(60,      1){\makebox(0,      0)[t]{$\bullet$}}
\put(93,     1){\makebox(0,     0)[t]{$\bullet$}}

\put(126,      1){\makebox(0,    0)[t]{$\bullet$}}

\put(159,      1){\makebox(0,    0)[t]{$\bullet$}}

\put(126,      1){\makebox(0,    0)[t]{$\bullet$}}
\put(28,      -1){\line(1,      0){33}}
\put(61,      -1){\line(1,      0){30}}
\put(94,     -1){\line(1,      0){30}}
\put(126,      -1){\line(1,      0){33}}

\put(22,     10){$q ^{-1}$}
\put(40,      5){$q^{}$}
\put(50,      -20){${-1}$}
\put(70,      -15){$q ^{-1}$}

\put(81,      -20){${-1}$}

\put(102,    -15){$q ^{}$}

\put(118,      10){${-1}$}

\put(135,     5){$q ^{-1}$}

\put(161,      10){$q ^{}$}

\put(91,    38){\makebox(0,     0)[t]{$\bullet$}}

\put(91,      -1){\line(0,     1){40}}

\put(99,      33){$q ^{}$}

\put(99,      16){$q ^{-1}$}

\put(186,        -1)  {by GDD $4$ of Row $12$  in Table C. }
\end{picture}$\\ \\
 It is quasi-affine     by Lemma \ref {2.63}.\\ \\

   $\ \ \ \ \  \ \   \begin{picture}(100,      15)  \put(-45,      -1){}

\put(-45,      -1){ {\rm (b)}}
\put(27,      1){\makebox(0,     0)[t]{$\bullet$}}
\put(60,      1){\makebox(0,      0)[t]{$\bullet$}}
\put(93,     1){\makebox(0,     0)[t]{$\bullet$}}

\put(126,      1){\makebox(0,    0)[t]{$\bullet$}}

\put(159,      1){\makebox(0,    0)[t]{$\bullet$}}

\put(126,      1){\makebox(0,    0)[t]{$\bullet$}}
\put(28,      -1){\line(1,      0){33}}
\put(61,      -1){\line(1,      0){30}}
\put(94,     -1){\line(1,      0){30}}
\put(126,      -1){\line(1,      0){33}}

\put(22,     10){$q ^{-1}$}
\put(40,      5){$q^{}$}
\put(50,      -20){${-1}$}
\put(70,      -15){$q ^{-1}$}

\put(81,      -20){${-1}$}

\put(102,    -15){$q ^{}$}

\put(118,      10){${-1}$}

\put(135,     5){$q ^{-1}$}

\put(161,      10){$-1 ^{}$}

\put(91,    38){\makebox(0,     0)[t]{$\bullet$}}

\put(91,      -1){\line(0,     1){40}}

\put(99,      33){$q ^{}$}

\put(99,      16){$q ^{-1}$}

\put(186,        -1)  {by GDD $10$ of Row $13$. }
\end{picture}$\\ \\
 It is quasi-affine    by Lemma \ref {2.63}.
{\rm (c)} There are not any other cases by Lemma \ref {2.13}.

 {\rm (iii)}   Adding on   Vertex 3. It is not quasi-affine since it contains a proper subGDD which is not arithmetic.

   {\rm (iv)} Adding on  Vertex 2.   Omitting   Vertex  Vertex 6  in whole is not arithmetic by Lemma \ref {2.13}.
     \subsection  {Quasi-affine  GDDs  over
GDD $18$ of Row $13$ in Table C}\label {sub3.2}
 {\rm (i)}  Omitting   Vertex 1 and adding on   Vertex 5.\\

 $\ \ \ \ \  \ \   \begin{picture}(100,      15)  \put(-45,      -1){} \put(-45,      -1){{\rm (a)} }
\put(27,      1){\makebox(0,     0)[t]{$\bullet$}}
\put(60,      1){\makebox(0,      0)[t]{$\bullet$}}
\put(93,     1){\makebox(0,     0)[t]{$\bullet$}}

\put(126,      1){\makebox(0,    0)[t]{$\bullet$}}

%\put(128,      1){\makebox(0,    0)[t]{$\bullet$}}
\put(159,      1){\makebox(0,    0)[t]{$\bullet$}}

\put(28,      -1){\line(1,      0){33}}
\put(61,      -1){\line(1,      0){30}}
\put(94,     -1){\line(1,      0){30}}

\put(128,      -1){\line(1,      0){33}}

\put(22,     -15){$q ^{-1}$}
\put(40,     -20){$q ^{}$}
\put(58,      -15){$q ^{-1}$}
\put(75,      -20){$q ^{}$}

\put(87,      -15){${-1}$}

\put(102,     -20){$q ^{-1}$}

\put(124,      -15){$q ^{}$}
\put(135,     -20){$q ^{-1}$}

\put(157,      -15){$q ^{}$}

\put(75,    18){\makebox(0,     0)[t]{$\bullet$}}

\put(60,     -1){\line(1,      0){30}}

\put(91,      0){\line(-1,     1){17}}
\put(60,     -1){\line(1,      1){17}}

\put(50,    12){$q$}

\put(68,    22){$-1$}
\put(91,    12){$q$}

\put(222,        -0)  {by GDD $3$ of Row $12$  in Table C. }
\end{picture}$\\ \\
 It is quasi-affine     by Lemma \ref {2.63}.\\

 $\ \ \ \ \  \ \   \begin{picture}(100,      15)  \put(-45,      -1){} \put(-45,      -1){{\rm (b)} }
\put(27,      1){\makebox(0,     0)[t]{$\bullet$}}
\put(60,      1){\makebox(0,      0)[t]{$\bullet$}}
\put(93,     1){\makebox(0,     0)[t]{$\bullet$}}

\put(126,      1){\makebox(0,    0)[t]{$\bullet$}}

\put(128,      1){\makebox(0,    0)[t]{$\bullet$}}
\put(159,      1){\makebox(0,    0)[t]{$\bullet$}}

\put(126,      1){\makebox(0,    0)[t]{$\bullet$}}
\put(28,      -1){\line(1,      0){33}}
\put(61,      -1){\line(1,      0){30}}
\put(94,     -1){\line(1,      0){30}}

\put(128,      -1){\line(1,      0){33}}

\put(22,     -15){$q ^{-1}$}
\put(40,     -20){$q ^{}$}
\put(58,      -15){$q ^{-1}$}
\put(75,      -20){$q ^{}$}

\put(87,      -15){${-1}$}

\put(102,     -20){$q ^{-1}$}

\put(124,      -15){$q ^{}$}
\put(135,     -20){$q ^{-1}$}

\put(157,      -15){$-1 ^{}$}

\put(75,    18){\makebox(0,     0)[t]{$\bullet$}}

\put(60,     -1){\line(1,      0){30}}

\put(91,      0){\line(-1,     1){17}}
\put(60,     -1){\line(1,      1){17}}

\put(50,    12){$q$}

\put(68,    22){$-1$}
\put(91,    12){$q$}

\put(222,        -1)  { by GDD $6$ of Row $13$  in Table C. }
\end{picture}$\\ \\
 It is quasi-affine     by Lemma \ref {2.63}.
{\rm (c)} There are not any other cases by Lemma \ref {2.13}.

 {\rm (ii)}   Adding on  Vertex 2 ,4. It is not quasi-affine since it contains a proper subGDD which is not arithmetic.

 {\rm (iii)} Omitting Vertex 1 and adding on   Vertex 3.  It    is empty  by Lemma \ref {2.63}.
    \subsection  {Quasi-affine  GDDs  over
GDD $19$ of Row $13$ in Table C}\label {sub3.2}
  {\rm (i)}  Omitting  Vertex  1 and adding on   Vertex 5.\\

 $\ \ \ \ \  \ \   \begin{picture}(100,      15)  \put(-45,      -1){}

\put(-45,      -1){ {\rm (a)}}
\put(27,      1){\makebox(0,     0)[t]{$\bullet$}}
\put(60,      1){\makebox(0,      0)[t]{$\bullet$}}
\put(93,     1){\makebox(0,     0)[t]{$\bullet$}}

\put(126,      1){\makebox(0,    0)[t]{$\bullet$}}

\put(159,      1){\makebox(0,    0)[t]{$\bullet$}}
\put(192,      1){\makebox(0,    0)[t]{$\bullet$}}

\put(126,      1){\makebox(0,    0)[t]{$\bullet$}}
\put(28,      -1){\line(1,      0){33}}
\put(61,      -1){\line(1,      0){30}}
\put(94,     -1){\line(1,      0){30}}
\put(126,      -1){\line(1,      0){33}}
\put(159,      -1){\line(1,      0){33}}

\put(22,     10){$q ^{-1}$}
\put(40,      5){$q ^{}$}
\put(54,      10){$-1$}
\put(70,      5){$q ^{-1}$}

\put(87,      10){${-1}$}

\put(102,     5){$q ^{-1}$}

\put(124,      10){$q ^{}$}
\put(135,     5){$q ^{-1}$}
\put(159,      10){$q ^{}$}

\put(168,     5){$q ^{-1}$}
\put(192,      10){$q$}

\put(222,        -1)  {by GDD $2$ of Row $12$  in Table C. }
\end{picture}$\\
 It is quasi-affine    by Lemma \ref {2.63}.\\

 $\ \ \ \ \  \ \   \begin{picture}(100,      15)  \put(-45,      -1){}

\put(-45,      -1){ {\rm (b)}}
\put(27,      1){\makebox(0,     0)[t]{$\bullet$}}
\put(60,      1){\makebox(0,      0)[t]{$\bullet$}}
\put(93,     1){\makebox(0,     0)[t]{$\bullet$}}

\put(126,      1){\makebox(0,    0)[t]{$\bullet$}}

\put(159,      1){\makebox(0,    0)[t]{$\bullet$}}
\put(192,      1){\makebox(0,    0)[t]{$\bullet$}}

\put(126,      1){\makebox(0,    0)[t]{$\bullet$}}
\put(28,      -1){\line(1,      0){33}}
\put(61,      -1){\line(1,      0){30}}
\put(94,     -1){\line(1,      0){30}}
\put(126,      -1){\line(1,      0){33}}
\put(159,      -1){\line(1,      0){33}}

\put(22,     10){$q ^{-1}$}
\put(40,      5){$q ^{}$}
\put(54,      10){$-1$}
\put(70,      5){$q ^{-1}$}

\put(87,      10){${-1}$}

\put(102,     5){$q ^{-1}$}

\put(124,      10){$q ^{}$}
\put(135,     5){$q ^{-1}$}
\put(159,      10){$q ^{}$}

\put(168,     5){$q ^{-1}$}
\put(192,      10){$ -1$}

\put(222,        -1)  {by GDD $3$ of Row $13$  in Table C. }
\end{picture}$\\
 It is quasi-affine    by Lemma \ref {2.63}.
{\rm (c)} There are not any other cases by Lemma \ref {2.6}.
 {\rm (ii)} cycle.
 {\rm (a)} to {\rm (a)}  is empty.{\rm (a)} to {\rm (a')}  is empty.{\rm (b)} to {\rm (a)}  is empty.\\ \\ \\ \\

 $ \ \ \ \ \  \ \  
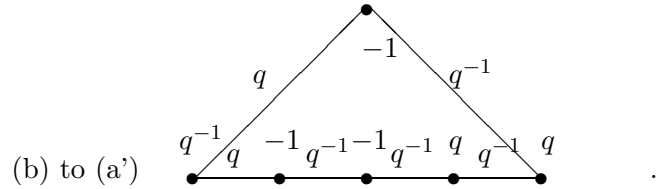
\begin{picture}(100,      15)  \put(-45,      -1){ {\rm (b)} to {\rm (a')}}

\put(27,      1){\makebox(0,     0)[t]{$\bullet$}}
\put(60,      1){\makebox(0,      0)[t]{$\bullet$}}
\put(93,     1){\makebox(0,     0)[t]{$\bullet$}}

\put(126,      1){\makebox(0,    0)[t]{$\bullet$}}

\put(159,      1){\makebox(0,    0)[t]{$\bullet$}}
\put(28,      -1){\line(1,      0){33}}
\put(61,      -1){\line(1,      0){30}}
\put(94,     -1){\line(1,      0){30}}
\put(126,      -1){\line(1,      0){33}}

\put(22,     10){$q ^{-1}$}
\put(40,      5){$q ^{}$}
\put(54,      10){$-1$}
\put(70,      5){$q ^{-1}$}

\put(87,      10){${-1}$}

\put(102,     5){$q ^{-1}$}

\put(124,      10){$q ^{}$}
\put(135,     5){$q ^{-1}$}
\put(159,      10){$q ^{}$}

\put(93,     65){\makebox(0,     0)[t]{$\bullet$}}

\put(28,      -1){\line(1,      1){65}}

\put(159,      -1){\line(-1,      1){65}}

\put(91,      45){$-1 ^{}$}
\put(50,      35){$q$}
\put(124,      35){$q ^{-1}$}

\put(200,  -1)  {. }
\end{picture}$\\
Omitting  Vertex 2,  Vertex 3,  Vertex 4  in whole are  arithmetic,  It is quasi-affine.
    \subsection  {Quasi-affine  GDDs  over
GDD $20$ of Row $13$ in Table C}\label {sub3.2}
 {\rm (i)}   Omitting Vertex 1 and adding on   Vertex 5.\\

 $\ \ \ \ \  \ \   \begin{picture}(100,      15)  \put(-45,      -1){}

\put(-45,      -1){ {\rm (a)}}
\put(27,      1){\makebox(0,     0)[t]{$\bullet$}}
\put(60,      1){\makebox(0,      0)[t]{$\bullet$}}
\put(93,     1){\makebox(0,     0)[t]{$\bullet$}}

\put(126,      1){\makebox(0,    0)[t]{$\bullet$}}

\put(159,      1){\makebox(0,    0)[t]{$\bullet$}}
\put(192,      1){\makebox(0,    0)[t]{$\bullet$}}

\put(126,      1){\makebox(0,    0)[t]{$\bullet$}}
\put(28,      -1){\line(1,      0){33}}
\put(61,      -1){\line(1,      0){30}}
\put(94,     -1){\line(1,      0){30}}
\put(126,      -1){\line(1,      0){33}}
\put(159,      -1){\line(1,      0){33}}

\put(22,     10){${-1}$}
\put(40,      5){$q ^{-1}$}
\put(58,      10){$-1$}
\put(74,      5){$q ^{}$}

\put(91,      10){$q ^{-1}$}

\put(102,     5){$q ^{-1}$}

\put(124,      10){$q ^{}$}
\put(135,     5){$q ^{-1}$}
\put(159,      10){$q ^{}$}

\put(168,     5){$q ^{-1}$}
\put(192,      10){$q$}

\put(222,        -1)  {by GDD $1$ of Row $12$  in Table C. }
\end{picture}$\\
 It is quasi-affine     by Lemma \ref {2.63}.\\

 $\ \ \ \ \  \ \   \begin{picture}(100,      15)  \put(-45,      -1){}

\put(-45,      -1){ {\rm (b)}}
\put(27,      1){\makebox(0,     0)[t]{$\bullet$}}
\put(60,      1){\makebox(0,      0)[t]{$\bullet$}}
\put(93,     1){\makebox(0,     0)[t]{$\bullet$}}

\put(126,      1){\makebox(0,    0)[t]{$\bullet$}}

\put(159,      1){\makebox(0,    0)[t]{$\bullet$}}
\put(192,      1){\makebox(0,    0)[t]{$\bullet$}}

\put(126,      1){\makebox(0,    0)[t]{$\bullet$}}
\put(28,      -1){\line(1,      0){33}}
\put(61,      -1){\line(1,      0){30}}
\put(94,     -1){\line(1,      0){30}}
\put(126,      -1){\line(1,      0){33}}
\put(159,      -1){\line(1,      0){33}}

\put(22,     10){${-1}$}
\put(40,      5){$q ^{-1}$}
\put(58,      10){$-1$}
\put(74,      5){$q ^{}$}

\put(91,      10){$q ^{-1}$}

\put(102,     5){$q ^{-1}$}

\put(124,      10){$q ^{}$}
\put(135,     5){$q ^{-1}$}
\put(159,      10){$q ^{}$}

\put(168,     5){$q ^{-1}$}
\put(192,      10){$-1$}

\put(222,        -1)  {by GDD $1$ of Row $13$  in Table C. }
\end{picture}$\\
 It is quasi-affine     by Lemma \ref {2.63}.
{\rm (c)} There are not any other cases by Lemma \ref {2.6}.
 {\rm (ii)} Cycle. {\rm (a)} to {\rm (a)}  is empty. {\rm (a)} to {\rm (a')}  is empty. {\rm (b)} to {\rm (a)}  is empty.\\ \\ \\ \\

 $\ \ \ \ \  \ \   
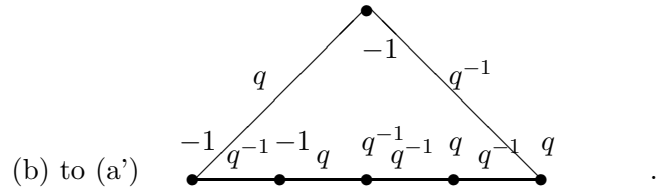
\begin{picture}(100,      15)  \put(-45,      -1){ {\rm (b)} to {\rm (a')}}

\put(-45,      -1){}
\put(27,      1){\makebox(0,     0)[t]{$\bullet$}}
\put(60,      1){\makebox(0,      0)[t]{$\bullet$}}
\put(93,     1){\makebox(0,     0)[t]{$\bullet$}}

\put(126,      1){\makebox(0,    0)[t]{$\bullet$}}

\put(159,      1){\makebox(0,    0)[t]{$\bullet$}}
\put(28,      -1){\line(1,      0){33}}
\put(61,      -1){\line(1,      0){30}}
\put(94,     -1){\line(1,      0){30}}
\put(126,      -1){\line(1,      0){33}}

\put(22,     10){${-1}$}
\put(40,      5){$q ^{-1}$}
\put(58,      10){$-1$}
\put(74,      5){$q ^{}$}

\put(91,      10){$q ^{-1}$}

\put(102,     5){$q ^{-1}$}

\put(124,      10){$q ^{}$}
\put(135,     5){$q ^{-1}$}
\put(159,      10){$q ^{}$}

\put(93,     65){\makebox(0,     0)[t]{$\bullet$}}

\put(28,      -1){\line(1,      1){65}}

\put(159,      -1){\line(-1,      1){65}}

\put(91,      45){${-1}$}
\put(50,      35){$q$}
\put(124,      35){$q ^{-1}$}

\put(200,  -1)  {. }
\end{picture}$\\
Omitting  Vertex 2, Vertex 3, Vertex 4  in whole are  arithmetic.  It is quasi-affine.
   \subsection  {Quasi-affine  GDDs  over
GDD $21$ of Row $13$ in Table C}\label {sub3.2}
{\rm (i)}   Omitting Vertex 1 and adding on   Vertex 5.
   It    is empty  by Row 11-13.
   \subsection  {Quasi-affine  GDDs  over
GDD $1$ of Row $14$ in Table C}\label {sub3.2}
 {\rm (i)}   Omitting Vertex 1 and adding on   Vertex 5.

 $\ \ \ \ \  \ \   \begin{picture}(100,      15)  \put(-45,      -1){}

\put(-45,      -1){ {\rm (a)}}
\put(27,      1){\makebox(0,     0)[t]{$\bullet$}}
\put(60,      1){\makebox(0,      0)[t]{$\bullet$}}
\put(93,     1){\makebox(0,     0)[t]{$\bullet$}}

\put(126,      1){\makebox(0,    0)[t]{$\bullet$}}

\put(159,      1){\makebox(0,    0)[t]{$\bullet$}}
\put(192,      1){\makebox(0,    0)[t]{$\bullet$}}

\put(126,      1){\makebox(0,    0)[t]{$\bullet$}}
\put(28,      -1){\line(1,      0){33}}
\put(61,      -1){\line(1,      0){30}}
\put(94,     -1){\line(1,      0){30}}
\put(126,      -1){\line(1,      0){33}}
\put(159,      -1){\line(1,      0){33}}

\put(22,     10){$q ^{}$}
\put(40,      5){$-q ^{}$}
\put(58,      10){$q$}
\put(74,      5){${-q}$}

\put(91,      10){${-1}$}

\put(102,     5){$-1 ^{}$}

\put(124,      10){$-1 ^{}$}
\put(135,     5){$-q ^{}$}
\put(159,      10){$q ^{}$}

\put(168,     5){$-q ^{}$}
\put(192,      10){$q$}

\put(222,        -1)  {   $q \in R_4,$ by GDD $1$ of Row $14$. }
\end{picture}$
\\
 It is quasi-affine    by Lemma \ref {2.6}.\\

 $\ \ \ \ \  \ \   \begin{picture}(100,      15)  \put(-45,      -1){}

\put(-45,      -1){ {\rm (a')}}
\put(27,      1){\makebox(0,     0)[t]{$\bullet$}}
\put(60,      1){\makebox(0,      0)[t]{$\bullet$}}
\put(93,     1){\makebox(0,     0)[t]{$\bullet$}}

\put(126,      1){\makebox(0,    0)[t]{$\bullet$}}

\put(159,      1){\makebox(0,    0)[t]{$\bullet$}}
\put(192,      1){\makebox(0,    0)[t]{$\bullet$}}

\put(126,      1){\makebox(0,    0)[t]{$\bullet$}}
\put(28,      -1){\line(1,      0){33}}
\put(61,      -1){\line(1,      0){30}}
\put(94,     -1){\line(1,      0){30}}
\put(126,      -1){\line(1,      0){33}}
\put(159,      -1){\line(1,      0){33}}

\put(22,     10){$q ^{}$}
\put(40,      5){$-q ^{}$}
\put(58,      10){$q$}
\put(74,      5){${-q}$}

\put(91,      10){${-1}$}

\put(102,     5){$-1 ^{}$}

\put(124,      10){$-1 ^{}$}
\put(135,     5){$-q ^{}$}
\put(159,      10){$q ^{}$}

\put(168,     5){$-q ^{}$}
\put(192,      10){$-1$}

\put(200,  -1)  { (  $q \in R_4$ ). }
\end{picture}$
\\
 It    is empty  by Lemma \ref {2.6}.
{\rm (b)} There are not any other cases by Row 14 or  by Lemma \ref {2.6}.
 {\rm (ii)} Cycle.
(a) to (a) empty.\\ \\ \\ \\

 $\ \ \ \ \  \ \   
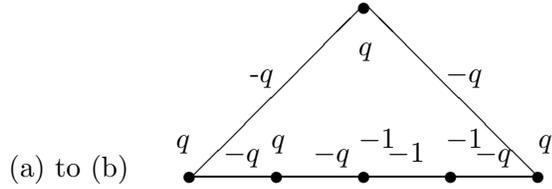
\begin{picture}(100,      15)  \put(-45,      -1){}

\put(-45,      -1){  (a) to (b) }
\put(27,      1){\makebox(0,     0)[t]{$\bullet$}}
\put(60,      1){\makebox(0,      0)[t]{$\bullet$}}
\put(93,     1){\makebox(0,     0)[t]{$\bullet$}}

\put(126,      1){\makebox(0,    0)[t]{$\bullet$}}

\put(159,      1){\makebox(0,    0)[t]{$\bullet$}}
\put(28,      -1){\line(1,      0){33}}
\put(61,      -1){\line(1,      0){30}}
\put(94,     -1){\line(1,      0){30}}
\put(126,      -1){\line(1,      0){33}}

\put(22,     10){$q ^{}$}
\put(40,      5){$-q ^{}$}
\put(58,      10){$q$}
\put(74,      5){$-q ^{}$}

\put(91,      10){${-1}$}

\put(102,     5){${-1}$}

\put(124,      10){${-1}$}
\put(135,     5){$-q ^{}$}
\put(159,      10){$q ^{}$}

\put(93,     65){\makebox(0,     0)[t]{$\bullet$}}

\put(28,      -1){\line(1,      1){65}}

\put(159,      -1){\line(-1,      1){65}}

\put(91,      45){$q ^{}$}
\put(50,      35){-$q ^{}$}
\put(124,      35){$-q ^{}$}

\put(200,  -1)  {. }
\end{picture}$\\ Omitting  Vertex 2,  Vertex 3,  Vertex 4 are arithmetic. It is quasi-affine.
    \subsection  {Quasi-affine  GDDs  over
GDD $2$ of Row $14$ in Table C}\label {sub3.2}
 {\rm (i)}   Omitting Vertex 1 and adding on   Vertex 5.\\

 $\ \ \ \ \  \ \   
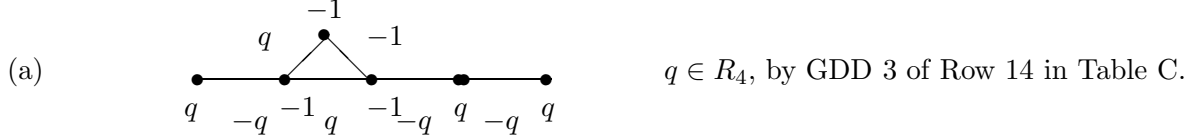
\begin{picture}(100,      15)  \put(-45,      -1){} \put(-45,      -1){{\rm (a)} }
\put(27,      1){\makebox(0,     0)[t]{$\bullet$}}
\put(60,      1){\makebox(0,      0)[t]{$\bullet$}}
\put(93,     1){\makebox(0,     0)[t]{$\bullet$}}

\put(126,      1){\makebox(0,    0)[t]{$\bullet$}}

\put(128,      1){\makebox(0,    0)[t]{$\bullet$}}
\put(159,      1){\makebox(0,    0)[t]{$\bullet$}}

%\put(126,      1){\makebox(0,    0)[t]{$\bullet$}}
\put(28,      -1){\line(1,      0){33}}
\put(61,      -1){\line(1,      0){30}}
\put(94,     -1){\line(1,      0){30}}

\put(128,      -1){\line(1,      0){33}}

\put(22,     -15){$q$}
\put(40,     -20){$-q ^{}$}
\put(58,      -15){$-1$}
\put(75,      -20){$q ^{}$}

\put(91,      -15){${-1}$}

\put(102,     -20){$-q ^{}$}

\put(124,      -15){$q ^{}$}
\put(135,     -20){$-q ^{}$}

\put(157,      -15){$q ^{}$}

\put(75,    18){\makebox(0,     0)[t]{$\bullet$}}

\put(60,     -1){\line(1,      0){30}}

\put(91,      0){\line(-1,     1){17}}
\put(60,     -1){\line(1,      1){17}}

\put(50,    12){$q$}

\put(68,    22){$-1$}
\put(91,    12){$-1$}

\put(200,  -1)  {   $q \in R_4,$ by GDD $3$ of Row $14$  in Table C. }
\end{picture}$\\ \\
It is repeated.
{\rm (b)} There are not any other cases by Row 14 or  by Lemma \ref {2.13}.

 {\rm (ii)}   Adding on  Vertex 2 ,4. It is not quasi-affine since it contains a proper subGDD which is not an arithmetic GDD.

 {\rm (iii)} Omitting Vertex 1 and adding on   Vertex 3.  It    is empty  by Row 14 or  by Lemma \ref {2.63}.
   \subsection  {Quasi-affine  GDDs  over
GDD $3$ of Row $14$ in Table C}\label {sub3.2}
  {\rm (i)}   Omitting Vertex 1 and adding on   Vertex 5.\\

 $\ \ \ \ \  \ \   
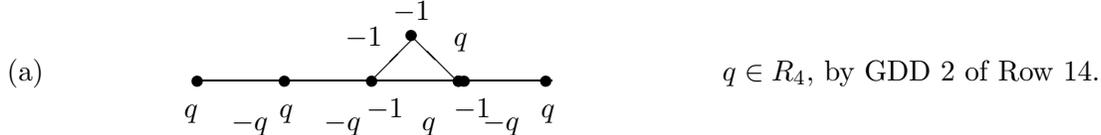
\begin{picture}(100,      15)  \put(-45,      -1){} \put(-45,      -1){{\rm (a)} }
\put(27,      1){\makebox(0,     0)[t]{$\bullet$}}
\put(60,      1){\makebox(0,      0)[t]{$\bullet$}}
\put(93,     1){\makebox(0,     0)[t]{$\bullet$}}

\put(126,      1){\makebox(0,    0)[t]{$\bullet$}}

\put(128,      1){\makebox(0,    0)[t]{$\bullet$}}
\put(159,      1){\makebox(0,    0)[t]{$\bullet$}}

\put(126,      1){\makebox(0,    0)[t]{$\bullet$}}
\put(28,      -1){\line(1,      0){33}}
\put(61,      -1){\line(1,      0){30}}
\put(94,     -1){\line(1,      0){30}}

\put(128,      -1){\line(1,      0){33}}

\put(22,     -15){$q$}
\put(40,     -20){$-q ^{}$}
\put(58,      -15){$q$}
\put(75,      -20){$-q ^{}$}

\put(91,      -15){${-1}$}

\put(112,     -20){$q ^{}$}

\put(124,      -15){$-1 ^{}$}
\put(135,     -20){$-q ^{}$}

\put(157,      -15){$q ^{}$}

\put(108,    18){\makebox(0,     0)[t]{$\bullet$}}

\put(93,     -1){\line(1,      0){30}}

\put(124,      0){\line(-1,     1){17}}
\put(93,     -1){\line(1,      1){17}}

\put(83,    12){$-1$}

\put(101,    22){$-1$}
\put(124,    12){$q$}

\put(222,        -1)  {   $q \in R_4,$ by GDD $2$ of Row $14$. }
\end{picture}$\\ \\
 It is quasi-affine    by Lemma \ref {2.63}.
{\rm (b)} There are not any other cases by Row 14 or by Lemma \ref {2.13}.

 {\rm (ii)} Omitting Vertex 1 and adding on   Vertex 4.  It    is empty   by Row 14.

 {\rm (iii)}   Adding on   Vertex 3. It is not quasi-affine since it contains a proper subGDD which is not an arithmetic GDD   by Row 14.

 {\rm (iv)} Omitting Vertex  5 and adding on  Vertex 2. It    is empty   by Row 14.
    \subsection  {Quasi-affine  GDDs  over
GDD $4$ of Row $14$ in Table C}\label {sub3.2}
 {\rm (i)}   Omitting Vertex 1 and adding on   Vertex 5.
It    is empty  by Row 14.
   \subsection  {Quasi-affine  GDDs  over
GDD $5$ of Row $14$ in Table C}\label {sub3.2}
 {\rm (i)}   Omitting  Vertex 5 and adding on   Vertex 3.\\ \\

 $\ \ \ \ \  \ \   
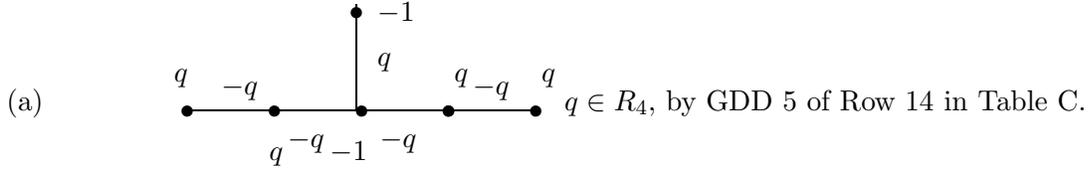
\begin{picture}(100,      15)  \put(-45,      -1){}

\put(-45,      -1){ {\rm (a)}}
\put(27,      1){\makebox(0,     0)[t]{$\bullet$}}
\put(60,      1){\makebox(0,      0)[t]{$\bullet$}}
\put(93,     1){\makebox(0,     0)[t]{$\bullet$}}

\put(126,      1){\makebox(0,    0)[t]{$\bullet$}}

\put(159,      1){\makebox(0,    0)[t]{$\bullet$}}

\put(126,      1){\makebox(0,    0)[t]{$\bullet$}}
\put(28,      -1){\line(1,      0){33}}
\put(61,      -1){\line(1,      0){30}}
\put(94,     -1){\line(1,      0){30}}
\put(126,      -1){\line(1,      0){33}}

\put(22,     10){$q ^{}$}
\put(40,      5){$-q^{}$}
\put(58,      -20){$q ^{}$}
\put(65,      -15){$-q ^{}$}

\put(81,      -20){${-1}$}

\put(100,    -15){$-q ^{}$}

\put(128,      10){$q ^{}$}

\put(135,     5){$-q ^{}$}

\put(161,      10){$q$}

\put(91,    38){\makebox(0,     0)[t]{$\bullet$}}

\put(91,      -1){\line(0,     1){40}}

\put(99,      33){${-1}$}

\put(99,      16){$q ^{}$}

\put(166,        -1)  {   $q \in R_4,$ by GDD $5$ of Row $14$  in Table C. }
\end{picture}$\\ \\
 It is quasi-affine    Row 19 or  by Lemma \ref {2.63}.
{\rm (b)} There are not any other cases by Row 14.

 {\rm (ii)}   Adding on  Vertex 2, 4. It is not quasi-affine since it contains a proper subGDD which is not an arithmetic GDD   by Row 14.

 {\rm (iii)} Omitting Vertex 1 and adding on   Vertex 5.\\ \\

 $\ \ \ \ \  \ \   
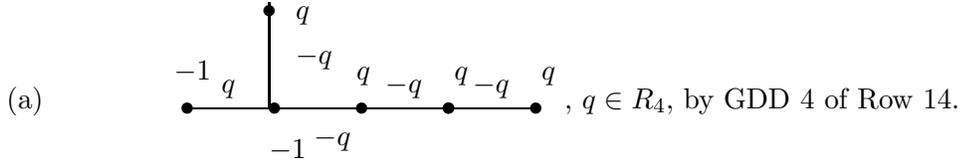
\begin{picture}(100,      15)  \put(-45,      -1){}

\put(-45,      -1){ \rm (a) }
\put(27,      1){\makebox(0,     0)[t]{$\bullet$}}
\put(60,      1){\makebox(0,      0)[t]{$\bullet$}}
\put(93,     1){\makebox(0,     0)[t]{$\bullet$}}

\put(126,      1){\makebox(0,    0)[t]{$\bullet$}}

\put(159,      1){\makebox(0,    0)[t]{$\bullet$}}

\put(126,      1){\makebox(0,    0)[t]{$\bullet$}}
\put(28,      -1){\line(1,      0){33}}
\put(61,      -1){\line(1,      0){30}}
\put(94,     -1){\line(1,      0){30}}
\put(126,      -1){\line(1,      0){33}}

\put(22,     10){${-1}$}
\put(40,      5){$q^{}$}
\put(58,      -20){${-1}$}
\put(75,      -15){$-q ^{}$}

\put(91,      10){$q ^{}$}

\put(102,     5){$-q ^{}$}

\put(128,      10){$q ^{}$}

\put(135,     5){$-q ^{}$}

\put(161,      10){$q$}

\put(58,    38){\makebox(0,     0)[t]{$\bullet$}}

\put(58,      -1){\line(0,     1){40}}

\put(68,      33){$q ^{}$}

\put(68,      16){$-q ^{}$}

\put(166,        -1)  { , $q\in R_4$, by GDD $4$ of Row $14$. }
\end{picture}$\\ \\ \\
 It is quasi-affine     by Row 19.\\ \\

 $\ \ \ \ \  \ \   
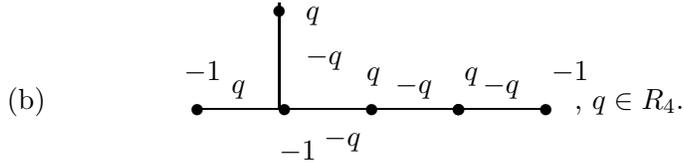
\begin{picture}(100,      15)  \put(-45,      -1){}

\put(-45,      -1){{\rm (b)} }
\put(27,      1){\makebox(0,     0)[t]{$\bullet$}}
\put(60,      1){\makebox(0,      0)[t]{$\bullet$}}
\put(93,     1){\makebox(0,     0)[t]{$\bullet$}}

\put(126,      1){\makebox(0,    0)[t]{$\bullet$}}

\put(159,      1){\makebox(0,    0)[t]{$\bullet$}}

\put(126,      1){\makebox(0,    0)[t]{$\bullet$}}
\put(28,      -1){\line(1,      0){33}}
\put(61,      -1){\line(1,      0){30}}
\put(94,     -1){\line(1,      0){30}}
\put(126,      -1){\line(1,      0){33}}

\put(22,     10){${-1}$}
\put(40,      5){$q^{}$}
\put(58,      -20){${-1}$}
\put(75,      -15){$-q ^{}$}

\put(91,      10){$q ^{}$}

\put(102,     5){$-q ^{}$}

\put(128,      10){$q ^{}$}

\put(135,     5){$-q ^{}$}

\put(161,      10){$-1$}

\put(58,    38){\makebox(0,     0)[t]{$\bullet$}}

\put(58,      -1){\line(0,     1){40}}

\put(68,      33){$q ^{}$}

\put(68,      16){$-q ^{}$}

\put(166,        -1)  { , $q\in R_4$. }
\end{picture}$\\ \\
It  is not arithmetic  by Row 14.
{\rm (c)} There are not any other cases by Row 14 or  by Lemma \ref {2.6}.
   \subsection  {Quasi-affine  GDDs  over
GDD $6$ of Row $14$}\label {sub3.2}
 {\rm (i)} Omitting  Vertex 5 and adding on   Vertex 3.\\ \\

 $\ \ \ \ \  \ \   
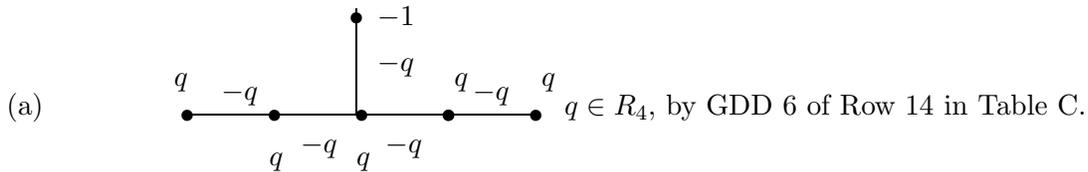
\begin{picture}(100,      15)  \put(-45,      -1){}

\put(-45,      -1){ {\rm (a)}}
\put(27,      1){\makebox(0,     0)[t]{$\bullet$}}
\put(60,      1){\makebox(0,      0)[t]{$\bullet$}}
\put(93,     1){\makebox(0,     0)[t]{$\bullet$}}

\put(126,      1){\makebox(0,    0)[t]{$\bullet$}}

\put(159,      1){\makebox(0,    0)[t]{$\bullet$}}

\put(126,      1){\makebox(0,    0)[t]{$\bullet$}}
\put(28,      -1){\line(1,      0){33}}
\put(61,      -1){\line(1,      0){30}}
\put(94,     -1){\line(1,      0){30}}
\put(126,      -1){\line(1,      0){33}}

\put(22,     10){$q ^{}$}
\put(40,      5){$-q^{}$}
\put(58,      -20){$q ^{}$}
\put(70,      -15){$-q ^{}$}

\put(91,      -20){$q ^{}$}

\put(102,    -15){$-q ^{}$}

\put(128,      10){$q ^{}$}

\put(135,     5){$-q ^{}$}

\put(161,      10){$q$}

\put(91,    38){\makebox(0,     0)[t]{$\bullet$}}

\put(91,      -1){\line(0,     1){40}}

\put(99,      33){${-1}$}

\put(99,      16){$-q ^{}$}

\put(166,        -1)  {   $q \in R_4,$ by GDD $6$ of Row $14$  in Table C. }
\end{picture}$\\ \\
 It is quasi-affine    by Row 19 or  Lemma \ref {2.2}
 {\rm (i)}.\\ \\

 $\ \ \ \ \  \ \   
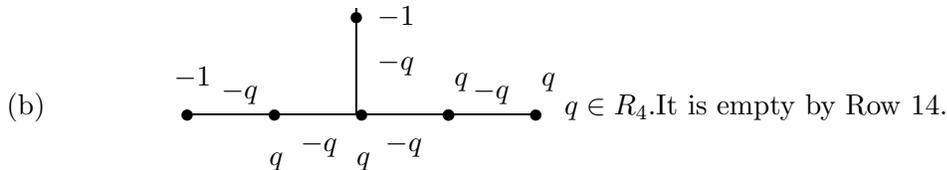
\begin{picture}(100,      15)  \put(-45,      -1){}

\put(-45,      -1){ {\rm (b)}}
\put(27,      1){\makebox(0,     0)[t]{$\bullet$}}
\put(60,      1){\makebox(0,      0)[t]{$\bullet$}}
\put(93,     1){\makebox(0,     0)[t]{$\bullet$}}

\put(126,      1){\makebox(0,    0)[t]{$\bullet$}}

\put(159,      1){\makebox(0,    0)[t]{$\bullet$}}

\put(126,      1){\makebox(0,    0)[t]{$\bullet$}}
\put(28,      -1){\line(1,      0){33}}
\put(61,      -1){\line(1,      0){30}}
\put(94,     -1){\line(1,      0){30}}
\put(126,      -1){\line(1,      0){33}}

\put(22,     10){$-1 ^{}$}
\put(40,      5){$-q^{}$}
\put(58,      -20){$q ^{}$}
\put(70,      -15){$-q ^{}$}

\put(91,      -20){$q ^{}$}

\put(102,    -15){$-q ^{}$}

\put(128,      10){$q ^{}$}

\put(135,     5){$-q ^{}$}

\put(161,      10){$q$}

\put(91,    38){\makebox(0,     0)[t]{$\bullet$}}

\put(91,      -1){\line(0,     1){40}}

\put(99,      33){${-1}$}

\put(99,      16){$-q ^{}$}

\put(166,        -1)  {   $q \in R_4$.It is empty by Row 14. }
\end{picture}$\\ \\
{\rm (c)} There are not any other cases  by Row 14 or  Lemma \ref {2.13}.

{\rm (ii)}    Adding on  Vertex 2.It is not quasi-affine since it contains a proper subGDD which is not an arithmetic GDD.

{\rm (iii)}    Omitting  Vertex 1 and  adding on  Vertex 4.
 There are not any other cases by Row 14.

 {\rm (iv)} Omitting Vertex 1 and adding on  Vertex  5.

 There are not any other cases by Row 14 or by Lemma \ref {2.77}.
   \subsection  {Quasi-affine  GDDs  over
GDD $1$ of Row $15$ in Table C }\label {sub3.2}
{\rm (i)}   Omitting Vertex 1 and adding on   Vertex 5.
 It    is empty   by Lemma \ref {2.6}.
    \subsection  {Quasi-affine  GDDs  over
GDD $2$ of Row $15$ in Table C }\label {sub3.2}
 {\rm (i)}   Omitting Vertex 1 and adding on   Vertex 5.
 It    is empty  by Row 15.

    {\rm (ii)}   Adding on  Vertex 2, 4. It is not quasi-affine since it contains a proper subGDD which is not an arithmetic GDD.

  {\rm (iii)}   Omitting Vertex 1 and adding on   Vertex 3.   It    is empty  by Row 15.
    \subsection  {Quasi-affine  GDDs  over
GDD $3$ of Row $15$ in Table C }\label {sub3.2}
 {\rm (i)}    Adding on  Vertex 2, 3, 4.
It is not quasi-affine since it contains a proper subGDD which is not an arithmetic GDD   by Row 15.
 {\rm (ii)}  Omitting Vertex 1 and adding on  Vertex  5.
It is empty by Lemma \ref {2.63}.
 \subsection  {Quasi-affine  GDDs  over
GDD $4$ of Row $15$ in Table C }\label {sub3.2}
   {\rm (i)}   Omitting  Vertex 5 and adding on   Vertex 3\\ \\

 $\ \ \ \ \  \ \   
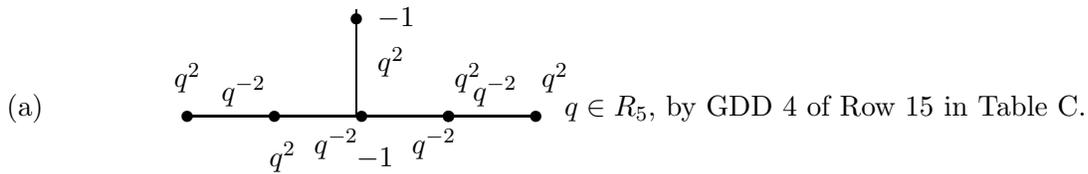
\begin{picture}(100,      15)  \put(-45,      -1){}

\put(-45,      -1){ {\rm (a)}}
\put(27,      1){\makebox(0,     0)[t]{$\bullet$}}
\put(60,      1){\makebox(0,      0)[t]{$\bullet$}}
\put(93,     1){\makebox(0,     0)[t]{$\bullet$}}

\put(126,      1){\makebox(0,    0)[t]{$\bullet$}}

\put(159,      1){\makebox(0,    0)[t]{$\bullet$}}

\put(126,      1){\makebox(0,    0)[t]{$\bullet$}}
\put(28,      -1){\line(1,      0){33}}
\put(61,      -1){\line(1,      0){30}}
\put(94,     -1){\line(1,      0){30}}
\put(126,      -1){\line(1,      0){33}}

\put(22,     10){$q ^{2}$}
\put(40,      5){$q^{-2}$}
\put(58,      -20){$q ^{2}$}
\put(75,      -15){$q ^{-2}$}

\put(91,      -20){${-1}$}

\put(112,    -15){$q ^{-2}$}

\put(128,      10){$q ^{2}$}

\put(135,     5){$q ^{-2}$}

\put(161,      10){$q^2$}

\put(91,    38){\makebox(0,     0)[t]{$\bullet$}}

\put(91,      -1){\line(0,     1){40}}

\put(99,      33){${-1}$}

\put(99,      16){$q ^{2}$}

\put(166,        -1)  {   $q \in R_5,$ by GDD $4$ of Row $15$  in Table C. }
\end{picture}$\\ \\
Omitting  Vertex 4 in whole  is not  arithmetic  by Row 15.\\ \\

 $\ \ \ \ \  \ \   
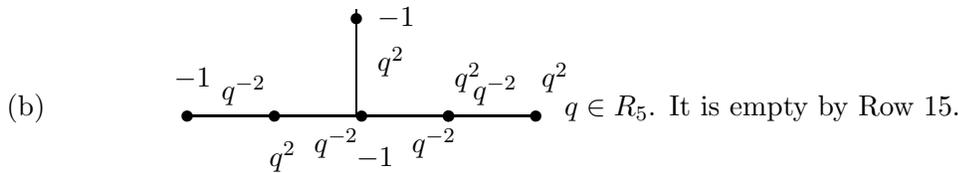
\begin{picture}(100,      15)  \put(-45,      -1){}

\put(-45,      -1){ {\rm (b)}}
\put(27,      1){\makebox(0,     0)[t]{$\bullet$}}
\put(60,      1){\makebox(0,      0)[t]{$\bullet$}}
\put(93,     1){\makebox(0,     0)[t]{$\bullet$}}

\put(126,      1){\makebox(0,    0)[t]{$\bullet$}}

\put(159,      1){\makebox(0,    0)[t]{$\bullet$}}

\put(126,      1){\makebox(0,    0)[t]{$\bullet$}}
\put(28,      -1){\line(1,      0){33}}
\put(61,      -1){\line(1,      0){30}}
\put(94,     -1){\line(1,      0){30}}
\put(126,      -1){\line(1,      0){33}}

\put(22,     10){$-1 ^{}$}
\put(40,      5){$q^{-2}$}
\put(58,      -20){$q ^{2}$}
\put(75,      -15){$q ^{-2}$}

\put(91,      -20){${-1}$}

\put(112,    -15){$q ^{-2}$}

\put(128,      10){$q ^{2}$}

\put(135,     5){$q ^{-2}$}

\put(161,      10){$q^2$}

\put(91,    38){\makebox(0,     0)[t]{$\bullet$}}

\put(91,      -1){\line(0,     1){40}}

\put(99,      33){${-1}$}

\put(99,      16){$q ^{2}$}

\put(166,        -1)  {  $q \in R_5$. It    is empty  by Row 15.}
\end{picture}$\\ \\
{\rm (c)}  There are not any other cases by Row 15 or by Lemma \ref {2.13}.

  {\rm (ii)} Omitting Vertex 1 and adding on  Vertex  5.  It    is empty  by Row 15.

 {\rm (iii)}   Adding on  Vertex 2 ,4, It is not quasi-affine since it contains a proper subGDD which is not an arithmetic GDD   by Row 15,
     \subsection  {Quasi-affine  GDDs  over
GDD $5$ of Row $15$ in Table C }\label {sub3.2}
 {\rm (i)}  Omitting Vertex  5 and adding on  Vertex  3.\\ \\

 $\ \ \ \ \  \ \   
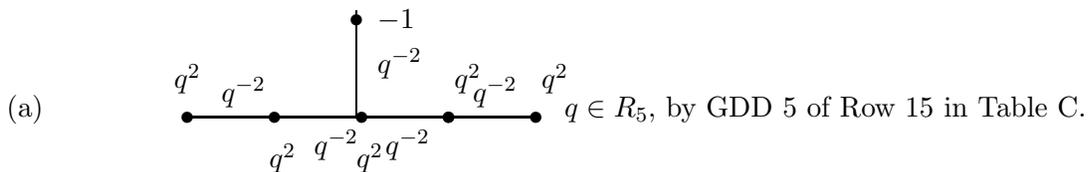
\begin{picture}(100,      15)  \put(-45,      -1){}

\put(-45,      -1){ {\rm (a)}}
\put(27,      1){\makebox(0,     0)[t]{$\bullet$}}
\put(60,      1){\makebox(0,      0)[t]{$\bullet$}}
\put(93,     1){\makebox(0,     0)[t]{$\bullet$}}

\put(126,      1){\makebox(0,    0)[t]{$\bullet$}}

\put(159,      1){\makebox(0,    0)[t]{$\bullet$}}

\put(126,      1){\makebox(0,    0)[t]{$\bullet$}}
\put(28,      -1){\line(1,      0){33}}
\put(61,      -1){\line(1,      0){30}}
\put(94,     -1){\line(1,      0){30}}
\put(126,      -1){\line(1,      0){33}}

\put(22,     10){$q ^{2}$}
\put(40,      5){$q^{-2}$}
\put(58,      -20){$q ^{2}$}
\put(75,      -15){$q ^{-2}$}

\put(91,      -20){$q ^{2}$}

\put(102,    -15){$q ^{-2}$}

\put(128,      10){$q ^{2}$}

\put(135,     5){$q ^{-2}$}

\put(161,      10){$q^2$}

\put(91,    38){\makebox(0,     0)[t]{$\bullet$}}

\put(91,      -1){\line(0,     1){40}}

\put(99,      33){${-1}$}

\put(99,      16){$q ^{-2}$}

\put(166,        -1)  {  $q \in R_5,$ by GDD $5$ of Row $15$  in Table C. }
\end{picture}$\\ \\
 It is quasi-affine    by Row 16-19 or by  Lemma \ref {2.2}
 {\rm (ii)}.\\ \\

 $\ \ \ \ \  \ \   
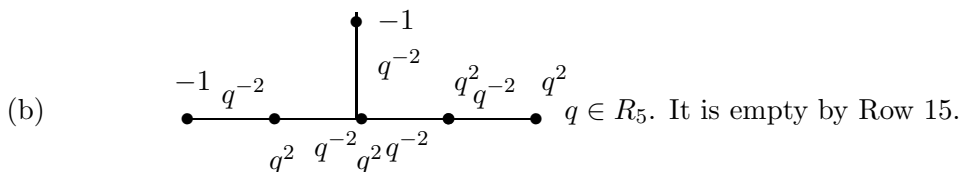
\begin{picture}(100,      15)  \put(-45,      -1){}

\put(-45,      -1){ {\rm (b)}}
\put(27,      1){\makebox(0,     0)[t]{$\bullet$}}
\put(60,      1){\makebox(0,      0)[t]{$\bullet$}}
\put(93,     1){\makebox(0,     0)[t]{$\bullet$}}

\put(126,      1){\makebox(0,    0)[t]{$\bullet$}}

\put(159,      1){\makebox(0,    0)[t]{$\bullet$}}

\put(126,      1){\makebox(0,    0)[t]{$\bullet$}}
\put(28,      -1){\line(1,      0){33}}
\put(61,      -1){\line(1,      0){30}}
\put(94,     -1){\line(1,      0){30}}
\put(126,      -1){\line(1,      0){33}}

\put(22,     10){$-1 ^{}$}
\put(40,      5){$q^{-2}$}
\put(58,      -20){$q ^{2}$}
\put(75,      -15){$q ^{-2}$}

\put(91,      -20){$q ^{2}$}

\put(102,    -15){$q ^{-2}$}

\put(128,      10){$q ^{2}$}

\put(135,     5){$q ^{-2}$}

\put(161,      10){$q^2$}

\put(91,    38){\makebox(0,     0)[t]{$\bullet$}}

\put(91,      -1){\line(0,     1){40}}

\put(99,      33){${-1}$}

\put(99,      16){$q ^{-2}$}

\put(166,        -1)  {   $q \in R_5$. It    is empty  by Row 15.}
\end{picture}$\\ \\
{\rm (c)} There are not any other cases by Row 15 or  by Lemma \ref {2.13}.

 {\rm (ii)}   Adding on  Vertex 2. It is not quasi-affine since it contains a proper subGDD which is not an arithmetic GDD   by Row 15.

   {\rm (iii)}  Omitting  Vertex 3 and   adding on    Vertex 4.
  There are not any  cases by Row 15.

  {\rm (iv)} Omitting  Vertex 3 and adding on  Vertex  5.
   \subsection  {Quasi-affine  GDDs  over
GDD $6$ of Row $15$ in Table C }\label {sub3.2}
  {\rm (i)}   Omitting Vertex 1 and adding on   Vertex 5.
   It    is empty  by Row 15.
    \subsection  {Quasi-affine  GDDs  over
GDD $7$ of Row $15$ in Table C }\label {sub3.2}
 {\rm (i)}   Omitting Vertex 1 and adding on   Vertex 5.
  It    is empty  by Row 15.
     \subsection  {Quasi-affine  GDDs  over
GDD $1$ of Row $16$ in Table C }\label {sub3.2}
 {\rm (i)} Omitting Vertex 1 and adding on   Vertex 6.

 {\rm (ii)}    Adding on  Vertex 2, 3, 5, It is not  It is quasi-affine    since it contains a proper subGDD which is not an arithmetic GDD.

 {\rm (iii)} Omitting Vertex 1 and adding on   Vertex 4.\\ \\ \\ \\

 $\ \ \ \ \  \ \
$\\ \\
 Omitting  Vertex 6  in whole is   GDD $1$ of Row $16$,  It is quasi-affine    by Row 17- 19 or Lemma \ref {2.13}.\\ \\ \\ \\

 $\ \ \ \ \  \ \   %
$\\ \\
It is repeated.
{\rm (c)} There are not any other cases by Lemma \ref {2.13}.
    \subsection  {Quasi-affine  GDDs  over
GDD $1$ of Row $17$ in Table C }\label {sub3.2}
 {\rm (i)} Omitting Vertex 1 and adding on   Vertex 6.\\

 $\ \ \ \ \  \ \   %
$
\\
 It is quasi-affine    by Row 21.\\

 $\ \ \ \ \  \ \   %
$\\
 It is repeated.
{\rm (c)} There are not any other cases by Lemma \ref {2.6}.
 {\rm (ii)} cycle.\\ \\ \\ \\ \\

 $ \ \ \ \ \  \ \  %
$\\
is Type 1. Omitting  Vertex 3, Vertex  4,  Vertex 5  in whole are simple chain.  It is quasi-affine.
{\rm (a)} to {\rm (b)}  is empty.{\rm (b)} to {\rm (a)}  is empty.
\\
\\ \\ \\ \\

 $ \ \ \ \ \  \ \  %
$\\
 It is repeated.
   \subsection  {Quasi-affine  GDDs  over
GDD $2$ of Row $17$ in Table C in Table C }\label {sub3.2}
 {\rm (i)} Adding on  Vertex 2 , 3, 4, 5.
   It is not quasi-affine since it contains a proper subGDD which is not an arithmetic GDD.

 {\rm (ii)} Omitting Vertex 1 and adding on   Vertex 6.\\

 $\ \ \ \ \  \ \   %
$\\ \\
 It is quasi-affine     by Lemma \ref {2.63}.\\

 $\ \ \ \ \  \ \   %
$\\  \\
 It is repeated.
{\rm (c)} There are not any other cases by Lemma \ref {2.13}.
    \subsection  {Quasi-affine  GDDs  over
GDD $3$ of Row $17$ in Table C }\label {sub3.2}
  {\rm (i)} Omitting Vertex 1 and adding on   Vertex 5.
It    is empty  by Lemma \ref {2.2} {\rm (ii)}.

 {\rm (ii)} Omitting Vertex 1 and adding on   Vertex 6.\\ \\

 $\ \ \ \ \  \ \   %
$\\ \\
 It is repeated.\\ \\

 $\ \ \ \ \  \ \   %
$\\ \\
 It is repeated.
{\rm (c)} There are not any other cases by Lemma \ref {2.13}.

 {\rm (iii)}  Adding on  Vertex  3, 4. It is not quasi-affine since it contains a proper subGDD which is not an arithmetic GDD.

 {\rm (iv)} Omitting  Vertex 6 and adding on  Vertex 2. It    is empty   by Lemma \ref {2.13}.
    \subsection  {Quasi-affine  GDDs  over
GDD $4$ of Row $17$ in Table C }\label {sub3.2}
 {\rm (i)} Omitting Vertex 1 and adding on   Vertex 6.\\ \\

 $\ \ \ \ \  \ \   %
$\\  \\
 It is quasi-affine    by Lemma \ref {2.63}.\\ \\

 $\ \ \ \ \  \ \   %
$\\ \\
 It is repeated.
{\rm (c)} There are not any other cases by Lemma \ref {2.13}.

 {\rm (ii)}   Adding on  Vertex 2, 3, 5. It is not quasi-affine since it contains a proper subGDD which is not an arithmetic GDD.

 {\rm (iii)}  Omitting Vertex 1 and  adding on   Vertex 4.  It    is empty   by Lemma \ref {2.2} {\rm (ii)}.
    \subsection  {Quasi-affine  GDDs  over
GDD $5$ of Row $17$ in Table C }\label {sub3.2}
 {\rm (i)}   Omitting Vertex 1 and adding on   Vertex 5.
It    is empty  by Lemma \ref {2.2} {\rm (ii)}.

   {\rm (ii)}   Omitting Vertex 1 and adding on  Vertex  6.\\ \\

 $\ \ \ \ \  \ \   %
$\\ \\
 It is repeated.
\\ \\

 $\ \ \ \ \  \ \   %
$\\ \\
It is repeated .
{\rm (c)} There are not any other cases by Lemma \ref {2.13}.

 {\rm (iii)}   Adding on  Vertex 3, 4. It is not quasi-affine since it contains a proper subGDD which is not an arithmetic GDD.

  {\rm (iv)}   Omitting Vertex 1 and adding on  Vertex 2.\\ \\

 $\ \ \ \ \  \ \   %
$\\ \\
It is repeated.
{\rm (b)} There are not any other cases by Lemma \ref {2.13}.
   \subsection  {Quasi-affine  GDDs  over
GDD $6$ of Row $17$ in Table C }\label {sub3.2}
 {\rm (i)}   Omitting Vertex 1 and adding on   Vertex 6.
{\ }

{\ }

{\ }

 $\ \ \ \ \  \ \   %
$\\ \\
 It is quasi-affine     by Lemma \ref {2.63}. \\ \\

 $\ \ \ \ \  \ \   %
$\\ \\
 It is repeated.
{\rm (c)} There are not any other cases by Lemma \ref {2.13}.

 {\rm (ii)}    Adding on  Vertex 2, 3, 5, It is not quasi-affine since it contains a proper subGDD which is not an arithmetic GDD.

 {\rm (iii)}   Omitting Vertex 1 and adding on   Vertex 4.  It    is empty  by  Lemma \ref {2.2} {\rm (ii)}.

    \subsection  {Quasi-affine  GDDs  over
GDD $7$ of Row $17$ in Table C }\label {sub3.2}
 {\rm (i)}   Omitting Vertex 1 and adding on   Vertex 6.

 {\rm (ii)}   Adding on  Vertex 2, 3, 5. It is not quasi-affine since it contains a proper subGDD which is not an arithmetic GDD.

  {\rm (iii)}   Omitting Vertex 1 and adding on  4.\\ \\ \\ \\

 $\ \ \ \ \  \ \   
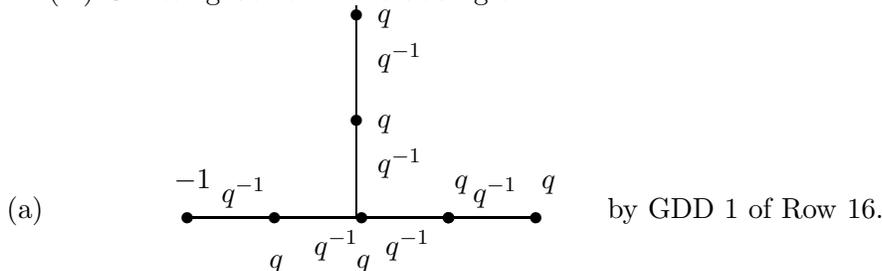
\begin{picture}(100,      15)  \put(-45,      -1){}

\put(-45,      -1){ {\rm (a)}}
\put(27,      1){\makebox(0,     0)[t]{$\bullet$}}
\put(60,      1){\makebox(0,      0)[t]{$\bullet$}}
\put(93,     1){\makebox(0,     0)[t]{$\bullet$}}

\put(126,      1){\makebox(0,    0)[t]{$\bullet$}}

\put(159,      1){\makebox(0,    0)[t]{$\bullet$}}

\put(126,      1){\makebox(0,    0)[t]{$\bullet$}}
\put(28,      -1){\line(1,      0){33}}
\put(61,      -1){\line(1,      0){30}}
\put(94,     -1){\line(1,      0){30}}
\put(126,      -1){\line(1,      0){33}}

\put(22,     10){$-1 ^{}$}
\put(40,      5){$q^{-1}$}
\put(58,      -20){${q}$}
\put(75,      -15){$q ^{-1}$}

\put(91,      -20){${q}$}

\put(102,    -15){$q ^{-1}$}

\put(128,      10){$q ^{}$}

\put(135,     5){$q ^{-1}$}

\put(161,      10){$q^{}$}

\put(91,    38){\makebox(0,     0)[t]{$\bullet$}}

\put(91,    78){\makebox(0,     0)[t]{$\bullet$}}

\put(91,      -1){\line(0,     1){80}}

\put(99,      16){$q ^{-1}$}

\put(99,      33){$q ^{}$}

\put(99,      56){$q ^{-1}$}

\put(99,      73){$q ^{}$}

\put(186,        -1)  {by  GDD $1$ of Row $16$. }
\end{picture}$\\ \\
 Omitting  Vertex 7 in whole is arithmetic by  GDD $7$ of Row $17$.
 It is quasi-affine    by Lemma \ref {2.63}.\\ \\ \\ \\

{\ }

{\ }

 $\ \ \ \ \  \ \   
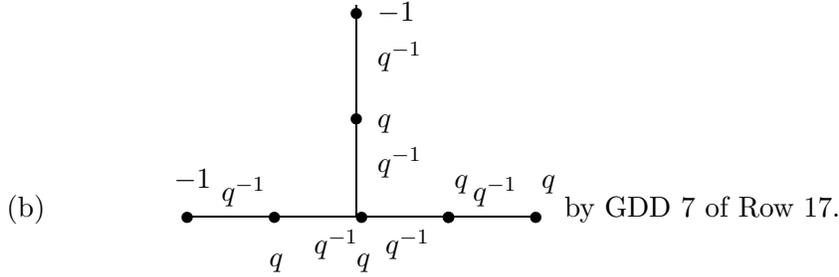
\begin{picture}(100,      15)  \put(-45,      -1){}

\put(-45,      -1){ {\rm (b)}}
\put(27,      1){\makebox(0,     0)[t]{$\bullet$}}
\put(60,      1){\makebox(0,      0)[t]{$\bullet$}}
\put(93,     1){\makebox(0,     0)[t]{$\bullet$}}

\put(126,      1){\makebox(0,    0)[t]{$\bullet$}}

\put(159,      1){\makebox(0,    0)[t]{$\bullet$}}

\put(126,      1){\makebox(0,    0)[t]{$\bullet$}}
\put(28,      -1){\line(1,      0){33}}
\put(61,      -1){\line(1,      0){30}}
\put(94,     -1){\line(1,      0){30}}
\put(126,      -1){\line(1,      0){33}}

\put(22,     10){$-1 ^{}$}
\put(40,      5){$q^{-1}$}
\put(58,      -20){${q}$}
\put(75,      -15){$q ^{-1}$}

\put(91,      -20){${q}$}

\put(102,    -15){$q ^{-1}$}

\put(128,      10){$q ^{}$}

\put(135,     5){$q ^{-1}$}

\put(161,      10){$q^{}$}

\put(91,    38){\makebox(0,     0)[t]{$\bullet$}}

\put(91,    78){\makebox(0,     0)[t]{$\bullet$}}

\put(91,      -1){\line(0,     1){80}}

\put(99,      16){$q ^{-1}$}

\put(99,      33){$q ^{}$}

\put(99,      56){$q ^{-1}$}

\put(99,      73){$-1 ^{}$}

\put(166,        -1)  { by GDD $7$ of Row $17$. }
\end{picture}$\\ \\
 Omitting  Vertex 7 in whole is not arithmetic by Row 17-18.
{\rm (c)} There are not any other cases by Lemma \ref {2.13}.
     \subsection  {Quasi-affine  GDDs  over
GDD $1$ of Row $18$ in Table C }\label {sub3.2}
 {\rm (i)}   Omitting Vertex 1 and adding on   Vertex 6.
   There are not any  cases by  Row 17-Row 18  or Lemma \ref {2.63}.
   \subsection  {Quasi-affine  GDDs  over
GDD $2$ of Row $18$ in Table C }\label {sub3.2}
 {\rm (i)}   Omitting Vertex 1 and adding on   Vertex 6.
It is empty by Lemma \ref {2.63}.
    \subsection  {Quasi-affine  GDDs  over
GDD $3$ of Row $18$ in Table C }\label {sub3.2}
 {\rm (i)}   Omitting Vertex  1 and adding on   Vertex 6.
There are not any  cases by  Lemma \ref {2.13}.
 There are not any  cases by Row 17-18.

   {\rm (iii)}   adding on   Vertex 3, 4. It is not quasi-affine since it contains a proper subGDD which is not an arithmetic GDD.

  {\rm (iv)}   Omitting Vertex  1 and adding on   Vertex 5.
It  is not  arithmetic   by Lemma \ref {2.63}.
  \subsection  {Quasi-affine  GDDs  over
GDD $4$ of Row $18$ in Table C }\label {sub3.2}
 {\rm (i)}   Omitting  Vertex 6 and adding on   Vertex 5.
Omitting Vertex 1  in whole is not an arithmetic GDD   by Lemma \ref {2.2}
 {\rm (ii)}.

   {\rm (ii)}   Adding on   Vertex 3,4. It is not quasi-affine since it contains a proper subGDD which is not an arithmetic GDD.

   {\rm (iii)}   omitting  Vertex 6 and adding on  Vertex 2.\\ \\

{\ }

 $\ \ \ \ \  \ \   
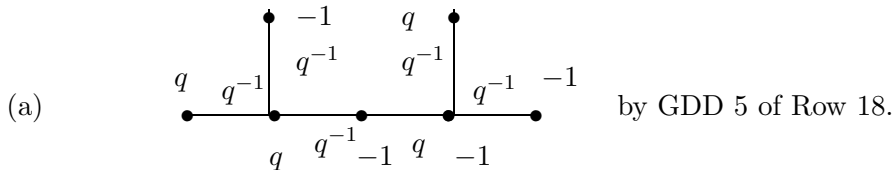
\begin{picture}(100,      15)  \put(-45,      -1){}

\put(-45,      -1){ {\rm (a)}}
\put(27,      1){\makebox(0,     0)[t]{$\bullet$}}
\put(60,      1){\makebox(0,      0)[t]{$\bullet$}}
\put(93,     1){\makebox(0,     0)[t]{$\bullet$}}

\put(126,      1){\makebox(0,    0)[t]{$\bullet$}}

\put(159,      1){\makebox(0,    0)[t]{$\bullet$}}

\put(126,      1){\makebox(0,    0)[t]{$\bullet$}}
\put(28,      -1){\line(1,      0){33}}
\put(61,      -1){\line(1,      0){30}}
\put(94,     -1){\line(1,      0){30}}
\put(126,      -1){\line(1,      0){33}}

\put(22,     10){$q$}
\put(40,      5){$q ^{-1}$}
\put(58,      -20){$q$}
\put(75,      -15){$q ^{-1}$}

\put(91,      -20){$-1$}

\put(112,     -15){$q ^{}$}

\put(128,      -20){$-1$}

\put(135,     5){$q ^{-1}$}

\put(161,      10){$-1$}

\put(58,    38){\makebox(0,     0)[t]{$\bullet$}}

\put(58,      -1){\line(0,     1){40}}

\put(68,      33){$-1$}

\put(68,      16){$q ^{-1}$}

\put(128,    38){\makebox(0,     0)[t]{$\bullet$}}

\put(128,      -1){\line(0,     1){40}}

\put(108,      33){$q$}

\put(108,      16){$q ^{-1}$}

\put(186,        -1)  { by  GDD $5$ of Row $18$. }
\end{picture}$\\
\\
Omitting  Vertex 6  in whole is not an arithmetic GDD   by Row 17-18.
{\rm (b)} There are not any other cases   by Row 17-Row 18.

   {\rm (v)}   omitting  Vertex 5 and adding on   Vertex 6.
Omitting Vertex 1  in whole is not an arithmetic GDD    Lemma \ref {2.63}.
     \subsection  {Quasi-affine  GDDs  over
GDD $5$ of Row $18$ in Table C }\label {sub3.2}
    {\rm (i)}   Omitting Vertex  1 and adding on   Vertex 5. There are not any  cases    by  Lemma \ref {2.2}  {\rm (ii)}.

 {\rm (ii)}   Omitting  Vertex 6 and adding on  Vertex 2.\\ \\

 $\ \ \ \ \  \ \   
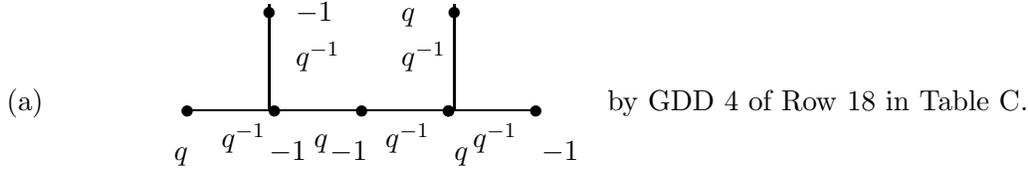
\begin{picture}(100,      15)  \put(-45,      -1){} \put(-45,      -1){ {\rm (a)}}
\put(27,      1){\makebox(0,     0)[t]{$\bullet$}}
\put(60,      1){\makebox(0,      0)[t]{$\bullet$}}
\put(93,     1){\makebox(0,     0)[t]{$\bullet$}}

\put(126,      1){\makebox(0,    0)[t]{$\bullet$}}

\put(159,      1){\makebox(0,    0)[t]{$\bullet$}}

\put(126,      1){\makebox(0,    0)[t]{$\bullet$}}
\put(28,      -1){\line(1,      0){33}}
\put(61,      -1){\line(1,      0){30}}
\put(94,     -1){\line(1,      0){30}}
\put(126,      -1){\line(1,      0){33}}

\put(22,    - 20){$q$}
\put(40,    -  15){$q ^{-1}$}
\put(58,      -20){$-1$}
\put(75,      -15){$q ^{}$}

\put(81,      -20){$-1$}

\put(102,     -15){$q ^{-1}$}

\put(128,      -20){$q$}

\put(135,     -15){$q ^{-1}$}

\put(161,      -20){$-1$}

\put(58,    38){\makebox(0,     0)[t]{$\bullet$}}

\put(58,      -1){\line(0,     1){40}}

\put(68,      33){$-1$}

\put(68,      16){$q ^{-1}$}

\put(128,    38){\makebox(0,     0)[t]{$\bullet$}}

\put(128,      -1){\line(0,     1){40}}

\put(108,      33){$q$}

\put(108,      16){$q ^{-1}$}

\put(186,        -1)  {by GDD $4$ of Row $18$  in Table C. }
\end{picture}$\\ \\
 Omitting  Vertex 6 in whole  is not an arithmetic GDD    by Row 17-19.
{\rm (b)} There are not any other cases by Row 17-18.

     {\rm (iii)}    Adding on   Vertex 3, 4. It is not quasi-affine since it contains a proper subGDD which is not an arithmetic GDD.

     {\rm (iv)}   Omitting Vertex 1 and adding on   Vertex 6.
 It    is empty   by Row 17-18.
     \subsection  {Quasi-affine  GDDs  over
GDD $6$ of Row $18$ in Table C }\label {sub3.2}
 {\rm (i)}   Omitting Vertex 1 and adding on   Vertex 5.
It is empty by  Lemma \ref {2.2} {\rm (ii)}.

   {\rm (ii)}   Omitting  Vertex 6 and adding on  Vertex 2.
 It    is empty  by Row 17-Row 18.

   {\rm (iii)}  Adding on   Vertex 3,4. It is not quasi-affine since it contains a proper subGDD which is not an arithmetic GDD.

   {\rm (iv)}   Omitting Vertex 1 and adding on   Vertex 6.
There are not any  cases by Row 17-18.
   \subsection  {Quasi-affine  GDDs  over
GDD $7$ of Row $18$ in Table C }\label {sub3.2}
   {\rm (i)}   Omitting  Vertex 6 and adding on   Vertex 5. It    is empty   Lemma \ref {2.2} {\rm (ii)}.\\ \\

 $\ \ \ \ \  \ \   
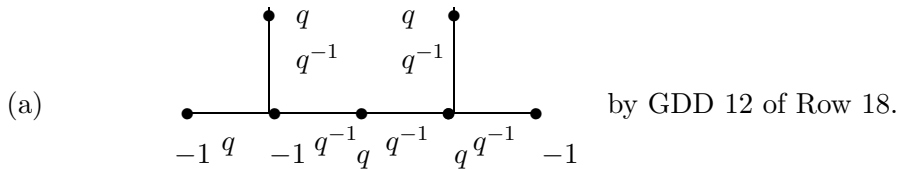
\begin{picture}(100,      15)  \put(-45,      -1){} \put(-45,      -1){ {\rm (a)}}
\put(27,      1){\makebox(0,     0)[t]{$\bullet$}}
\put(60,      1){\makebox(0,      0)[t]{$\bullet$}}
\put(93,     1){\makebox(0,     0)[t]{$\bullet$}}

\put(126,      1){\makebox(0,    0)[t]{$\bullet$}}

\put(159,      1){\makebox(0,    0)[t]{$\bullet$}}

\put(126,      1){\makebox(0,    0)[t]{$\bullet$}}
\put(28,      -1){\line(1,      0){33}}
\put(61,      -1){\line(1,      0){30}}
\put(94,     -1){\line(1,      0){30}}
\put(126,      -1){\line(1,      0){33}}

\put(22,    - 20){$-1$}
\put(40,    -  15){$q$}
\put(58,      -20){$-1$}
\put(75,      -15){$q ^{-1}$}

\put(91,      -20){$q$}

\put(102,     -15){$q ^{-1}$}

\put(128,      -20){$q$}

\put(135,     -15){$q ^{-1}$}

\put(161,      -20){$-1$}

\put(58,    38){\makebox(0,     0)[t]{$\bullet$}}

\put(58,      -1){\line(0,     1){40}}

\put(68,      33){$q$}

\put(68,      16){$q ^{-1}$}

\put(128,    38){\makebox(0,     0)[t]{$\bullet$}}

\put(128,      -1){\line(0,     1){40}}

\put(108,      33){$q$}

\put(108,      16){$q ^{-1}$}

\put(186,        -1)  {by GDD $12$ of Row $18.$   }
\end{picture}$\\ \\
Omitting Vertex 1  in whole is not an arithmetic GDD    Lemma \ref {2.13}.\\
 {\rm (b)} There are not any other cases by Row 17-18.

 {\rm (ii)}   Adding on   Vertex 3, 4. It is not quasi-affine since it contains a proper subGDD which is not an arithmetic GDD.

   {\rm (iii)}   Omitting Vertex 1 and  adding on   Vertex 6.  It    is empty  by Row 17-18.
    \subsection  {Quasi-affine  GDDs  over
GDD $8$ of Row $18$ in Table C }\label {sub3.2}
 {\rm (i)}   Omitting Vertex 1 and adding on  5.
 There are not any  cases by Lemma \ref {2.2} {\rm (ii)}.

     {\rm (ii)}   Omitting Vertex 1 and adding on   Vertex 6.
 There are not any  cases by Lemma \ref {2.63}.

         {\rm (iii)}   Adding on  Vertex 2.
Omitting Vertex  7 in whole  is not an arithmetic GDD   by Lemma \ref {2.13}.

 {\rm (iv)}   Adding on   Vertex 3, 4, It is not quasi-affine since it contains a proper subGDD which is not an arithmetic GDD.
    \subsection  {Quasi-affine  GDDs  over
GDD $9$ of Row $18$ in Table C }\label {sub3.2}
{\rm (i)}   Omitting Vertex 1 and adding on   Vertex 6.\\ \\

 $\ \ \ \ \  \ \   
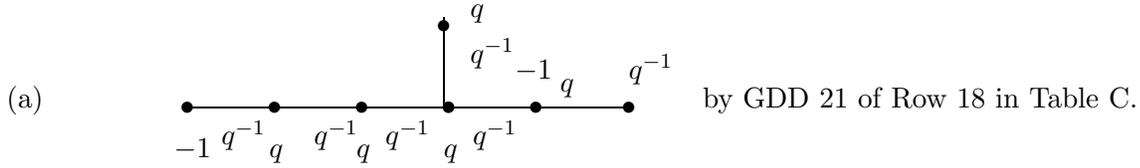
\begin{picture}(100,      15)  \put(-45,      -1){}

\put(-45,      -1){ {\rm (a)}}
\put(27,      1){\makebox(0,     0)[t]{$\bullet$}}
\put(60,      1){\makebox(0,      0)[t]{$\bullet$}}
\put(93,     1){\makebox(0,     0)[t]{$\bullet$}}

\put(126,      1){\makebox(0,    0)[t]{$\bullet$}}

\put(159,      1){\makebox(0,    0)[t]{$\bullet$}}
\put(194,      1){\makebox(0,    0)[t]{$\bullet$}}

\put(126,      1){\makebox(0,    0)[t]{$\bullet$}}
\put(28,      -1){\line(1,      0){33}}
\put(61,      -1){\line(1,      0){30}}
\put(94,     -1){\line(1,      0){30}}
\put(126,      -1){\line(1,      0){33}}
\put(159,      -1){\line(1,      0){33}}

\put(22,     -20){${-1}$}
\put(40,      -15){$q ^{-1}$}
\put(58,      -20){$q ^{}$}
\put(75,      -15){$q ^{-1}$}

\put(91,      -20){$q ^{}$}

\put(102,    -15){$q ^{-1}$}

\put(124,      -20){$q ^{}$}

\put(135,     -15){$q ^{-1}$}

\put(151,      10){$-1$}

\put(168,     5){$q$}

\put(194,      10){$q^{-1}$}

\put(124,    32){\makebox(0,     0)[t]{$\bullet$}}

\put(124,      -1){\line(0,     1){34}}

\put(134,      33){$q ^{}$}

\put(134,      16){$q ^{-1}$}

\put(222,        -1)  {by GDD $21$ of Row $18$  in Table C. }
\end{picture}$\\ \\
 It is repeated.\\ \\

 $\ \ \ \ \  \ \   
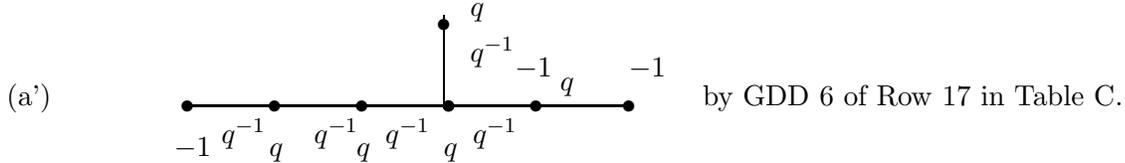
\begin{picture}(100,      15)  \put(-45,      -1){}

\put(-45,      -1){ {\rm (a')}}
\put(27,      1){\makebox(0,     0)[t]{$\bullet$}}
\put(60,      1){\makebox(0,      0)[t]{$\bullet$}}
\put(93,     1){\makebox(0,     0)[t]{$\bullet$}}

\put(126,      1){\makebox(0,    0)[t]{$\bullet$}}

\put(159,      1){\makebox(0,    0)[t]{$\bullet$}}
\put(194,      1){\makebox(0,    0)[t]{$\bullet$}}

\put(126,      1){\makebox(0,    0)[t]{$\bullet$}}
\put(28,      -1){\line(1,      0){33}}
\put(61,      -1){\line(1,      0){30}}
\put(94,     -1){\line(1,      0){30}}
\put(126,      -1){\line(1,      0){33}}
\put(159,      -1){\line(1,      0){33}}

\put(22,     -20){${-1}$}
\put(40,      -15){$q ^{-1}$}
\put(58,      -20){$q ^{}$}
\put(75,      -15){$q ^{-1}$}

\put(91,      -20){$q ^{}$}

\put(102,    -15){$q ^{-1}$}

\put(124,      -20){$q ^{}$}

\put(135,     -15){$q ^{-1}$}

\put(151,      10){$-1$}

\put(168,     5){$q ^{}$}

\put(194,      10){$-1$}

\put(124,    32){\makebox(0,     0)[t]{$\bullet$}}

\put(124,      -1){\line(0,     1){34}}

\put(134,      33){$q ^{}$}

\put(134,      16){$q ^{-1}$}

\put(222,        -1)  {by GDD $6$ of Row $17$  in Table C. }
\end{picture}$\\ \\
 It is repeated.
{\rm (c)} There are not any other cases by Lemma \ref {2.13}.

     {\rm (ii)}   Omitting Vertex 1 and adding on   Vertex 5.
 It    is empty  by  Lemma \ref {2.2}  {\rm (ii)}.

       {\rm (iii)}  Adding on  Vertex 2.\\ \\

 $\ \ \ \ \  \ \   
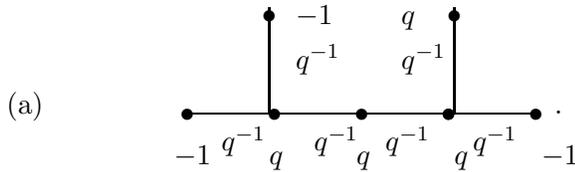
\begin{picture}(100,      15)  \put(-45,      -1){} \put(-45,      -1){ {\rm (a)}}
\put(27,      1){\makebox(0,     0)[t]{$\bullet$}}
\put(60,      1){\makebox(0,      0)[t]{$\bullet$}}
\put(93,     1){\makebox(0,     0)[t]{$\bullet$}}

\put(126,      1){\makebox(0,    0)[t]{$\bullet$}}

\put(159,      1){\makebox(0,    0)[t]{$\bullet$}}

\put(126,      1){\makebox(0,    0)[t]{$\bullet$}}
\put(28,      -1){\line(1,      0){33}}
\put(61,      -1){\line(1,      0){30}}
\put(94,     -1){\line(1,      0){30}}
\put(126,      -1){\line(1,      0){33}}

\put(22,    - 20){$-1$}
\put(40,    -  15){$q ^{-1}$}
\put(58,      -20){$q$}
\put(75,      -15){$q ^{-1}$}

\put(91,      -20){$q$}

\put(102,     -15){$q ^{-1}$}

\put(128,      -20){$q$}

\put(135,     -15){$q ^{-1}$}

\put(161,      -20){$-1$}

\put(58,    38){\makebox(0,     0)[t]{$\bullet$}}

\put(58,      -1){\line(0,     1){40}}

\put(68,      33){$-1$}

\put(68,      16){$q ^{-1}$}

\put(128,    38){\makebox(0,     0)[t]{$\bullet$}}

\put(128,      -1){\line(0,     1){40}}

\put(108,      33){$q$}

\put(108,      16){$q ^{-1}$}

\put(166,        -1)  {. }
\end{picture}$\\ \\  Omitting 1 in whole is  GDD $9$ of Row $18$  in Table C.
Omitting  Vertex 6 in whole  is not an arithmetic GDD   by  Lemma \ref {2.2}
 {\rm (i)}.\\  \\

 $\ \ \ \ \  \ \   
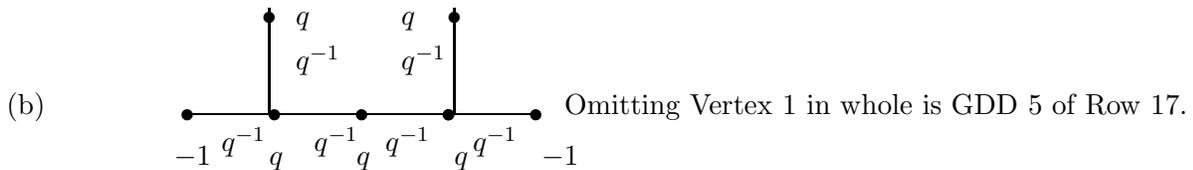
\begin{picture}(100,      15)  \put(-45,      -1){} \put(-45,      -1){ {\rm (b)}}
\put(27,      1){\makebox(0,     0)[t]{$\bullet$}}
\put(60,      1){\makebox(0,      0)[t]{$\bullet$}}
\put(93,     1){\makebox(0,     0)[t]{$\bullet$}}

\put(126,      1){\makebox(0,    0)[t]{$\bullet$}}

\put(159,      1){\makebox(0,    0)[t]{$\bullet$}}

\put(126,      1){\makebox(0,    0)[t]{$\bullet$}}
\put(28,      -1){\line(1,      0){33}}
\put(61,      -1){\line(1,      0){30}}
\put(94,     -1){\line(1,      0){30}}
\put(126,      -1){\line(1,      0){33}}

\put(22,    - 20){$-1$}
\put(40,    -  15){$q ^{-1}$}
\put(58,      -20){$q$}
\put(75,      -15){$q ^{-1}$}

\put(91,      -20){$q$}

\put(102,     -15){$q ^{-1}$}

\put(128,      -20){$q$}

\put(135,     -15){$q ^{-1}$}

\put(161,      -20){$-1$}

\put(58,    38){\makebox(0,     0)[t]{$\bullet$}}

\put(58,      -1){\line(0,     1){40}}

\put(68,      33){$q$}

\put(68,      16){$q ^{-1}$}

\put(128,    38){\makebox(0,     0)[t]{$\bullet$}}

\put(128,      -1){\line(0,     1){40}}

\put(108,      33){$q$}

\put(108,      16){$q ^{-1}$}

\put(166,        -1)  {  Omitting  Vertex 1 in whole is GDD $5$ of Row $17$. }
\end{picture}$\\ \\
Omitting  Vertex 7 in whole  is GDD $5$ of Row $17$. Omitting  Vertex 6  in whole is GDD $9$ of Row $18$.
 It is quasi-affine.\\ \\

 $\ \ \ \ \  \ \   
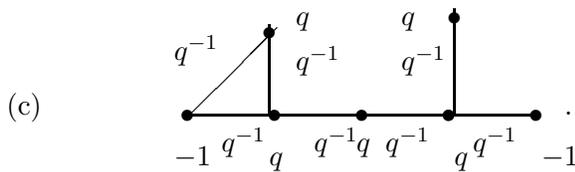
\begin{picture}(100,      15)  \put(-45,      -1){} \put(-45,      -1){ \rm (c)}
\put(27,      1){\makebox(0,     0)[t]{$\bullet$}}
\put(60,      1){\makebox(0,      0)[t]{$\bullet$}}
\put(93,     1){\makebox(0,     0)[t]{$\bullet$}}

\put(126,      1){\makebox(0,    0)[t]{$\bullet$}}

\put(159,      1){\makebox(0,    0)[t]{$\bullet$}}

\put(126,      1){\makebox(0,    0)[t]{$\bullet$}}
\put(28,      -1){\line(1,      0){33}}
\put(61,      -1){\line(1,      0){30}}
\put(94,     -1){\line(1,      0){30}}
\put(126,      -1){\line(1,      0){33}}

\put(22,     -20){${-1}$}
\put(40,      -15){$q ^{-1}$}
\put(58,      -20){$q$}
\put(75,      -15){$q ^{-1}$}

\put(91,      -15){$q$}

\put(102,     -15){$q ^{-1}$}

\put(128,      -20){$q$}

\put(135,     -15){$q ^{-1}$}

\put(161,      -20){$-1$}

\put(58,    32){\makebox(0,     0)[t]{$\bullet$}}

\put(58,      -1){\line(0,     1){34}}

\put(22,    20){$q ^{-1}$}

\put(68,      33){$q$}

\put(68,      16){$q ^{-1}$}

\put(28,      -1){\line(1,      1){33}}

\put(128,    38){\makebox(0,     0)[t]{$\bullet$}}

\put(128,      -1){\line(0,     1){40}}

\put(108,      33){$q$}

\put(108,      16){$q ^{-1}$}

\put(166,        -1)  {  . }
\end{picture}$\\ \\
Omitting  Vertex 6 in whole  is not  arithmetic  by Row 17-Row 18.
  {\rm (d)} There are not any other cases by Lemma \ref {2.63}.

 {\rm (iv)}   Adding on   Vertex 3, 4.  It is not quasi-affine since it contains a proper subGDD which is not an arithmetic GDD.
     \subsection  {Quasi-affine  GDDs  over
GDD $10$ of Row $18$ in Table C }\label {sub3.2}
 {\rm (i)}   Omitting Vertex 1 and adding on   Vertex 5.
There are not any  cases by Lemma \ref {2.2} {\rm (ii)}.

     {\rm (ii)}   Omitting Vertex 1 and adding on  Vertex  6.
 There are not any  cases by  Row 17-Row 18.

       {\rm (iii)}   Omitting Vertex 1 and adding on  Vertex 2.
Omitting  Vertex 7 in whole  is not an arithmetic GDD   by Row 17-18.

 {\rm (iv)} Adding on   Vertex 3, 4.  It is   not quasi-affine since it contains a proper subGDD which is not arithmetic.
             \subsection  {Quasi-affine  GDDs  over
GDD $11$ of Row $18$ in Table C }\label {sub3.2}
 {\rm (i)}   Omitting Vertex 1 and adding on   Vertex 6.
 There are not any  cases by Row 17-18 or  Lemma \ref {2.63}.

     {\rm (ii)}   Adding on  Vertex 2 ,3, 4, 5. It is not quasi-affine since it contains a proper subGDD which is not an arithmetic GDD.
\subsection  {Quasi-affine  GDDs  over
GDD $12$ of Row $18$ in Table C }\label {sub3.2}
 {\rm (i)}   Omitting Vertex 1 and adding on   Vertex 5.
 There are not any  cases by Lemma \ref {2.2} {\rm (ii)}.

     {\rm (ii)}   Omitting Vertex 1 and adding on   Vertex 6.\\ \\

 $\ \ \ \ \  \ \   
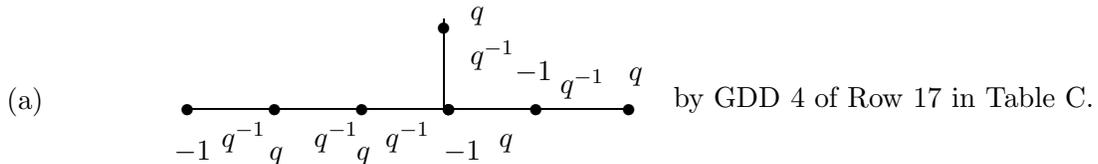
\begin{picture}(100,      15)  \put(-45,      -1){}

\put(-45,      -1){ {\rm (a)}}
\put(27,      1){\makebox(0,     0)[t]{$\bullet$}}
\put(60,      1){\makebox(0,      0)[t]{$\bullet$}}
\put(93,     1){\makebox(0,     0)[t]{$\bullet$}}

\put(126,      1){\makebox(0,    0)[t]{$\bullet$}}

\put(159,      1){\makebox(0,    0)[t]{$\bullet$}}
\put(194,      1){\makebox(0,    0)[t]{$\bullet$}}

\put(126,      1){\makebox(0,    0)[t]{$\bullet$}}
\put(28,      -1){\line(1,      0){33}}
\put(61,      -1){\line(1,      0){30}}
\put(94,     -1){\line(1,      0){30}}
\put(126,      -1){\line(1,      0){33}}
\put(159,      -1){\line(1,      0){33}}

\put(22,     -20){${-1}$}
\put(40,      -15){$q ^{-1}$}
\put(58,      -20){$q ^{}$}
\put(75,      -15){$q ^{-1}$}

\put(91,      -20){$q ^{}$}

\put(102,    -15){$q ^{-1}$}

\put(124,      -20){${-1}$}

\put(145,     -15){$q ^{}$}

\put(151,      10){$-1$}

\put(168,     5){$q ^{-1}$}

\put(194,      10){$q$}

\put(124,    32){\makebox(0,     0)[t]{$\bullet$}}

\put(124,      -1){\line(0,     1){34}}

\put(134,      33){$q ^{}$}

\put(134,      16){$q ^{-1}$}

\put(211,       -0)  {by GDD $4$ of Row $17$  in Table C. }
\end{picture}$\\ \\
 It is quasi-affine    by Row 21 or  Lemma \ref {2.63}.\\ \\

 $\ \ \ \ \  \ \   
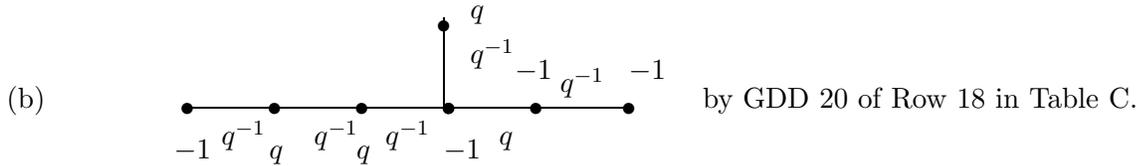
\begin{picture}(100,      15)  \put(-45,      -1){}

\put(-45,      -1){ {\rm (b)}}
\put(27,      1){\makebox(0,     0)[t]{$\bullet$}}
\put(60,      1){\makebox(0,      0)[t]{$\bullet$}}
\put(93,     1){\makebox(0,     0)[t]{$\bullet$}}

\put(126,      1){\makebox(0,    0)[t]{$\bullet$}}

\put(159,      1){\makebox(0,    0)[t]{$\bullet$}}
\put(194,      1){\makebox(0,    0)[t]{$\bullet$}}

\put(126,      1){\makebox(0,    0)[t]{$\bullet$}}
\put(28,      -1){\line(1,      0){33}}
\put(61,      -1){\line(1,      0){30}}
\put(94,     -1){\line(1,      0){30}}
\put(126,      -1){\line(1,      0){33}}
\put(159,      -1){\line(1,      0){33}}

\put(22,     -20){${-1}$}
\put(40,      -15){$q ^{-1}$}
\put(58,      -20){$q ^{}$}
\put(75,      -15){$q ^{-1}$}

\put(91,      -20){$q ^{}$}

\put(102,    -15){$q ^{-1}$}

\put(124,      -20){${-1}$}

\put(145,     -15){$q ^{}$}

\put(151,      10){$-1$}

\put(168,     5){$q ^{-1}$}

\put(194,      10){$-1$}

\put(124,    32){\makebox(0,     0)[t]{$\bullet$}}

\put(124,      -1){\line(0,     1){34}}

\put(134,      33){$q ^{}$}

\put(134,      16){$q ^{-1}$}

\put(222,        -1)  {by GDD $20$ of Row $18$  in Table C. }
\end{picture}$\\ \\
 It is repeated.
{\rm (c)} There are not any other cases by Lemma \ref {2.13}.

       {\rm (iii)}   Omitting Vertex 1 and adding on  Vertex 2.
       Omitting  Vertex 7 in whole  is not an arithmetic GDD   by Lemma \ref {2.13}.

 {\rm (iv)}   Adding on   Vertex 3,4. It is not quasi-affine since it contains a proper subGDD which is not an arithmetic GDD.
             \subsection  {Quasi-affine  GDDs  over
GDD $13$ of Row $18$ in Table C }\label {sub3.2}
 {\rm (i)}   Omitting Vertex 1 and adding on   Vertex 6.

 There are not any other cases  by Row 17-Row 18.

          (ii)   Adding on  Vertex 2, 3, 4, 5. They are  not quasi-affine since it contains a proper subGDD which is not an arithmetic GDD.
            \subsection  {Quasi-affine  GDDs  over
GDD $14$ of Row $18$ in Table C }\label {sub3.2}
 {\rm (i)}   Omitting Vertex 1 and adding on   Vertex 6.
 There are not any  cases  by Row 17-Row 18 or  Lemma \ref {2.63}.
             \subsection  {Quasi-affine  GDDs  over
GDD $15$ of Row $18$ in Table C } \label {sub3.2}
{\rm (i)}   Omitting Vertex 1 and adding on   Vertex 6.\\

 $\ \ \ \ \  \ \   
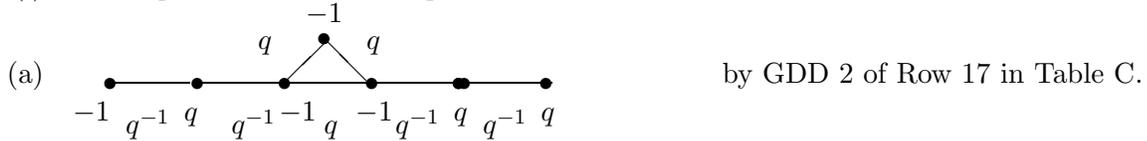
\begin{picture}(100,      15)  \put(-45,      -1){} \put(-45,      -1){{\rm (a)} }

\put(-6,      1){\makebox(0,     0)[t]{$\bullet$}}
\put(27,      1){\makebox(0,     0)[t]{$\bullet$}}
\put(60,      1){\makebox(0,      0)[t]{$\bullet$}}
\put(93,     1){\makebox(0,     0)[t]{$\bullet$}}

\put(126,      1){\makebox(0,    0)[t]{$\bullet$}}

\put(128,      1){\makebox(0,    0)[t]{$\bullet$}}
\put(159,      1){\makebox(0,    0)[t]{$\bullet$}}

\put(126,      1){\makebox(0,    0)[t]{$\bullet$}}
\put(28,      -1){\line(1,      0){33}}
\put(61,      -1){\line(1,      0){30}}
\put(94,     -1){\line(1,      0){30}}

\put(128,      -1){\line(1,      0){33}}
\put(-20,     -15){$-1$}
\put(0,     -20){$q ^{-1}$}

\put(22,     -15){$q$}
\put(40,     -20){$q ^{-1}$}
\put(58,      -15){$-1$}
\put(75,      -20){$q ^{}$}

\put(87,      -15){${-1}$}

\put(102,     -20){$q ^{-1}$}

\put(124,      -15){$q ^{}$}
\put(135,     -20){$q ^{-1}$}

\put(157,      -15){$q ^{}$}

\put(75,    18){\makebox(0,     0)[t]{$\bullet$}}

\put(-6,     -1){\line(1,      0){30}}
\put(60,     -1){\line(1,      0){30}}

\put(91,      0){\line(-1,     1){17}}
\put(60,     -1){\line(1,      1){17}}

\put(50,    12){$q$}

\put(68,    22){$-1$}
\put(91,    12){$q$}

\put(222,        -1)  { by GDD $2$ of Row $17$  in Table C. }
\end{picture}$\\ \\   It is quasi-affine     by Row 21 or by Lemma \ref {2.63}\\

 $\ \ \ \ \  \ \   
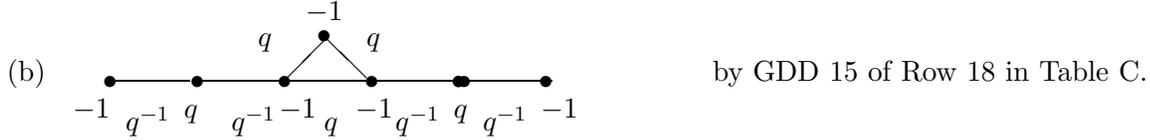
\begin{picture}(100,      15)  \put(-45,      -1){}
 \put(-45,      -1){{\rm (b)} }

\put(-6,      1){\makebox(0,     0)[t]{$\bullet$}}
\put(27,      1){\makebox(0,     0)[t]{$\bullet$}}
\put(60,      1){\makebox(0,      0)[t]{$\bullet$}}
\put(93,     1){\makebox(0,     0)[t]{$\bullet$}}

\put(126,      1){\makebox(0,    0)[t]{$\bullet$}}

\put(128,      1){\makebox(0,    0)[t]{$\bullet$}}
\put(159,      1){\makebox(0,    0)[t]{$\bullet$}}

\put(126,      1){\makebox(0,    0)[t]{$\bullet$}}
\put(28,      -1){\line(1,      0){33}}
\put(61,      -1){\line(1,      0){30}}
\put(94,     -1){\line(1,      0){30}}

\put(128,      -1){\line(1,      0){33}}
\put(-20,     -15){$-1$}
\put(0,     -20){$q ^{-1}$}

\put(22,     -15){$q$}
\put(40,     -20){$q ^{-1}$}
\put(58,      -15){$-1$}
\put(75,      -20){$q ^{}$}

\put(87,      -15){${-1}$}

\put(102,     -20){$q ^{-1}$}

\put(124,      -15){$q ^{}$}
\put(135,     -20){$q ^{-1}$}

\put(157,      -15){$-1 ^{}$}

\put(75,    18){\makebox(0,     0)[t]{$\bullet$}}

\put(-6,     -1){\line(1,      0){30}}
\put(60,     -1){\line(1,      0){30}}

\put(91,      0){\line(-1,     1){17}}
\put(60,     -1){\line(1,      1){17}}

\put(50,    12){$q$}

\put(68,    22){$-1$}
\put(91,    12){$q$}

\put(222,        -1)  {by GDD $15$ of Row $18$  in Table C. }
\end{picture}$\\ \\
  It is quasi-affine     by Row 21 or  Lemma \ref {2.63}.
{\rm (c)} There are not any other cases by Lemma \ref {2.13}.

 {\rm (ii)}  Adding on  Vertex 2, Vertex 3, Vertex  4, Vertex 5. It is not quasi-affine since it contains a proper subGDD which is not arithmetic.
            \subsection  {Quasi-affine  GDDs  over
GDD $16$ of Row $18$ in Table C } \label {sub3.2}
 {\rm (i)} Omitting Vertex 1 and adding on   Vertex 6.
\\
 There are not any other cases by Row $17$ and  Row $18$.
            \subsection  {Quasi-affine  GDDs  over
GDD $17$ of Row $18$ in Table C }\label {sub3.2}
 {\rm (i)} Omitting Vertex 1 and  adding on  Vertex  6.
 There are not any cases  by Row 17-Row 18.

 {\rm (ii)} Omitting Vertex 1 and  adding on   Vertex 4.

 There are not any other cases by Lemma \ref {2.13}.

 {\rm (iii)}   Adding on  Vertex 2 ,3, 5. It is not quasi-affine since it contains a proper subGDD which is not an arithmetic GDD.
          \subsection  {Quasi-affine  GDDs  over
GDD $18$ of Row $18$ in Table C } \label {sub3.2}
{\rm (i)} Omitting Vertex 1 and adding on   Vertex 6.
 There are not any  cases by Row 17-Row 18 or Lemma \ref {2.63}.
             \subsection  {Quasi-affine  GDDs  over
GDD $19$ of Row $18$ in Table C } \label {sub3.2}
 {\rm (i)} Omitting Vertex 1 and  adding on  Vertex  6.\\

 $\ \ \ \ \  \ \   
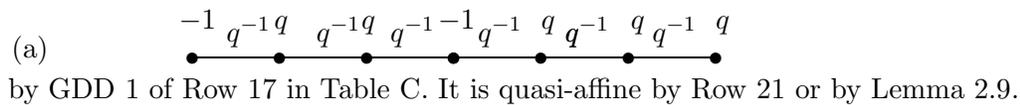
\begin{picture}(100,      15)  \put(-45,      -1){}

\put(-45,      -1){ {\rm (a)}}
\put(27,      1){\makebox(0,     0)[t]{$\bullet$}}
\put(60,      1){\makebox(0,      0)[t]{$\bullet$}}
\put(93,     1){\makebox(0,     0)[t]{$\bullet$}}

\put(126,      1){\makebox(0,    0)[t]{$\bullet$}}

\put(159,      1){\makebox(0,    0)[t]{$\bullet$}}
\put(192,      1){\makebox(0,    0)[t]{$\bullet$}}
\put(225,      1){\makebox(0,    0)[t]{$\bullet$}}

\put(126,      1){\makebox(0,    0)[t]{$\bullet$}}
\put(28,      -1){\line(1,      0){33}}
\put(61,      -1){\line(1,      0){30}}
\put(94,     -1){\line(1,      0){30}}
\put(126,      -1){\line(1,      0){33}}
\put(159,      -1){\line(1,      0){33}}
\put(192,      -1){\line(1,      0){33}}

\put(22,     10){${-1}$}
\put(40,      5){$q ^{-1}$}
\put(58,      10){$q$}
\put(74,      5){$q ^{-1}$}

\put(91,      10){$q ^{}$}

\put(102,     5){$q ^{-1}$}

\put(120,      10){${-1}$}
\put(135,     5){$q ^{-1}$}
\put(159,      10){$q ^{}$}

\put(168,     5){$q ^{-1}$}

\put(168,     5){$q ^{}$}

\put(193,      10){$q$}

\put(201,     5){$q ^{-1}$}
\put(225,      10){$q$}

\put(222,        -1)  { }
\end{picture}$\\ by GDD $1$ of Row $17$  in Table C.
  It is quasi-affine    by Row 21 or  by Lemma \ref {2.63}.\\

 $\ \ \ \ \  \ \   
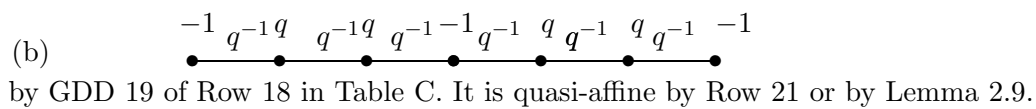
\begin{picture}(100,      15)  \put(-45,      -1){}

\put(-45,      -1){ {\rm (b)}}
\put(27,      1){\makebox(0,     0)[t]{$\bullet$}}
\put(60,      1){\makebox(0,      0)[t]{$\bullet$}}
\put(93,     1){\makebox(0,     0)[t]{$\bullet$}}

\put(126,      1){\makebox(0,    0)[t]{$\bullet$}}

\put(159,      1){\makebox(0,    0)[t]{$\bullet$}}
\put(192,      1){\makebox(0,    0)[t]{$\bullet$}}
\put(225,      1){\makebox(0,    0)[t]{$\bullet$}}

\put(126,      1){\makebox(0,    0)[t]{$\bullet$}}
\put(28,      -1){\line(1,      0){33}}
\put(61,      -1){\line(1,      0){30}}
\put(94,     -1){\line(1,      0){30}}
\put(126,      -1){\line(1,      0){33}}
\put(159,      -1){\line(1,      0){33}}
\put(192,      -1){\line(1,      0){33}}

\put(22,     10){${-1}$}
\put(40,      5){$q ^{-1}$}
\put(58,      10){$q$}
\put(74,      5){$q ^{-1}$}

\put(91,      10){$q ^{}$}

\put(102,     5){$q ^{-1}$}

\put(120,      10){${-1}$}
\put(135,     5){$q ^{-1}$}
\put(159,      10){$q ^{}$}

\put(168,     5){$q ^{-1}$}

\put(168,     5){$q ^{}$}

\put(193,      10){$q$}

\put(201,     5){$q ^{-1}$}
\put(225,      10){$-1$}

\put(222,        -1)  { }
\end{picture}$\\ by GDD $19$ of Row $18$  in Table C.
  It is quasi-affine    by Row 21 or by Lemma \ref {2.63}.\\
{\rm (c)} There are not any other cases by Lemma \ref {2.6}.
 {\rm (ii)} cycle.
{\rm (a)} to {\rm (a)}  is empty.{\rm (a)} to {\rm (b)}  is empty.{\rm (b)} to {\rm (a)}  is empty.\\ \\ \\ \\ \\

 $\ \ \ \ \  \ \   \begin{picture}(100,      15)  \put(-45,      -1){ {\rm (b)} to {\rm (b)}}

\put(-45,      -1){ }
\put(27,      1){\makebox(0,     0)[t]{$\bullet$}}
\put(60,      1){\makebox(0,      0)[t]{$\bullet$}}
\put(93,     1){\makebox(0,     0)[t]{$\bullet$}}

\put(126,      1){\makebox(0,    0)[t]{$\bullet$}}

\put(159,      1){\makebox(0,    0)[t]{$\bullet$}}
\put(192,      1){\makebox(0,    0)[t]{$\bullet$}}

\put(126,      1){\makebox(0,    0)[t]{$\bullet$}}
\put(28,      -1){\line(1,      0){33}}
\put(61,      -1){\line(1,      0){30}}
\put(94,     -1){\line(1,      0){30}}
\put(126,      -1){\line(1,      0){33}}
\put(159,      -1){\line(1,      0){33}}

\put(22,     10){${-1}$}
\put(40,      5){$q ^{-1}$}
\put(58,      10){$q$}
\put(74,      5){$q ^{-1}$}

\put(91,      10){$q ^{}$}

\put(102,     5){$q ^{-1}$}

\put(120,      10){${-1}$}
\put(135,     5){$q ^{-1}$}
\put(159,      10){$q ^{}$}

\put(168,     5){$q ^{-1}$}
\put(192,      10){$q$}

\put(112,     85){\makebox(0,     0)[t]{$\bullet$}}

\put(28,      -1){\line(1,      1){85}}

\put(192,      -1){\line(-1,      1){85}}

\put(100,      65){${-1}$}
\put(50,      35){$q$}
\put(155,      35){$q ^{-1}$}

\put(200,        -1)  {.  }
\end{picture}$\\
\\
Omitting  Vertex 2
 $\ \ \ \ \  \ \   \begin{picture}(100,      15)  \put(-45,      -1){}

\put(-45,      -1){}
\put(27,      1){\makebox(0,     0)[t]{$\bullet$}}
\put(60,      1){\makebox(0,      0)[t]{$\bullet$}}
\put(93,     1){\makebox(0,     0)[t]{$\bullet$}}

\put(126,      1){\makebox(0,    0)[t]{$\bullet$}}

\put(159,      1){\makebox(0,    0)[t]{$\bullet$}}
\put(192,      1){\makebox(0,    0)[t]{$\bullet$}}

\put(126,      1){\makebox(0,    0)[t]{$\bullet$}}
\put(28,      -1){\line(1,      0){33}}
\put(61,      -1){\line(1,      0){30}}
\put(94,     -1){\line(1,      0){30}}
\put(126,      -1){\line(1,      0){33}}
\put(159,      -1){\line(1,      0){33}}

\put(22,     10){$q ^{}$}
\put(40,      5){$q ^{-1}$}
\put(58,      10){$-1$}
\put(74,      5){$q ^{-1}$}

\put(91,      10){$q ^{}$}

\put(102,     5){$q ^{-1}$}

\put(124,      10){$q ^{}$}
\put(135,     5){$q ^{-1}$}
\put(153,      10){${-1}$}

\put(168,     5){$q ^{}$}
\put(192,      10){$-1$}

\put(200,        -1)  { is Type 1.}
\end{picture}$\\
Omitting  Vertex 3,  4, 5  in whole are simple chains.  It is quasi-affine.
             \subsection  {Quasi-affine  GDDs  over
GDD $20$ of Row $18$ in Table C }\label {sub3.2}
 {\rm (i)} Omitting Vertex 1 and  adding on   Vertex 6.\\ \\

 $\ \ \ \ \  \ \   
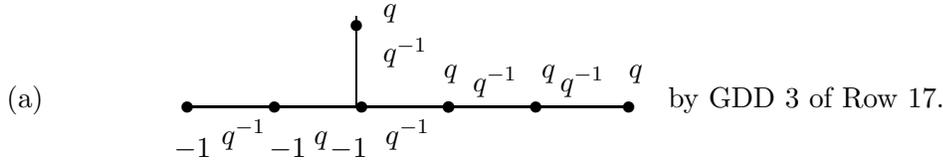
\begin{picture}(100,      15)  \put(-45,      -1){}

\put(-45,      -1){ {\rm (a)}}
\put(27,      1){\makebox(0,     0)[t]{$\bullet$}}
\put(60,      1){\makebox(0,      0)[t]{$\bullet$}}
\put(93,     1){\makebox(0,     0)[t]{$\bullet$}}

\put(126,      1){\makebox(0,    0)[t]{$\bullet$}}

\put(159,      1){\makebox(0,    0)[t]{$\bullet$}}
\put(194,      1){\makebox(0,    0)[t]{$\bullet$}}

\put(126,      1){\makebox(0,    0)[t]{$\bullet$}}
\put(28,      -1){\line(1,      0){33}}
\put(61,      -1){\line(1,      0){30}}
\put(94,     -1){\line(1,      0){30}}
\put(126,      -1){\line(1,      0){33}}
\put(159,      -1){\line(1,      0){33}}

\put(22,     -20){${-1}$}
\put(40,      -15){$q ^{-1}$}
\put(58,      -20){$-1$}
\put(75,      -15){$q ^{}$}

\put(81,      -20){$-1$}

\put(102,    -15){$q ^{-1}$}

\put(124,      10){$q ^{}$}

\put(135,     5){$q ^{-1}$}

\put(161,      10){$q$}

\put(168,     5){$q ^{-1}$}

\put(194,      10){$q$}

\put(91,    32){\makebox(0,     0)[t]{$\bullet$}}

\put(91,      -1){\line(0,     1){34}}

\put(101,      33){$q ^{}$}

\put(101,      16){$q ^{-1}$}

\put(209,        -1)  {by GDD $3$ of Row $17$. }
\end{picture}$\\ \\
 It is quasi-affine    by Row 21 or by Lemma \ref {2.63}.\\ \\

 $\ \ \ \ \  \ \   
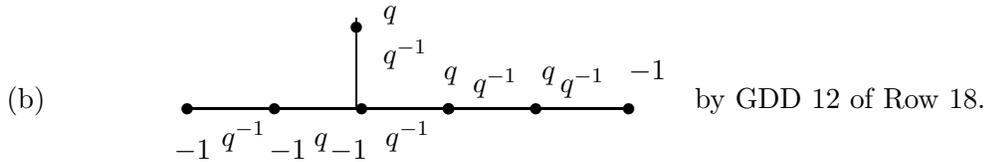
\begin{picture}(100,      15)  \put(-45,      -1){}

\put(-45,      -1){ {\rm (b)}}
\put(27,      1){\makebox(0,     0)[t]{$\bullet$}}
\put(60,      1){\makebox(0,      0)[t]{$\bullet$}}
\put(93,     1){\makebox(0,     0)[t]{$\bullet$}}

\put(126,      1){\makebox(0,    0)[t]{$\bullet$}}

\put(159,      1){\makebox(0,    0)[t]{$\bullet$}}
\put(194,      1){\makebox(0,    0)[t]{$\bullet$}}

\put(126,      1){\makebox(0,    0)[t]{$\bullet$}}
\put(28,      -1){\line(1,      0){33}}
\put(61,      -1){\line(1,      0){30}}
\put(94,     -1){\line(1,      0){30}}
\put(126,      -1){\line(1,      0){33}}
\put(159,      -1){\line(1,      0){33}}

\put(22,     -20){${-1}$}
\put(40,      -15){$q ^{-1}$}
\put(58,      -20){$-1$}
\put(75,      -15){$q ^{}$}

\put(81,      -20){$-1$}

\put(102,    -15){$q ^{-1}$}

\put(124,      10){$q ^{}$}

\put(135,     5){$q ^{-1}$}

\put(161,      10){$q$}

\put(168,     5){$q ^{-1}$}

\put(194,      10){$-1$}

\put(91,    32){\makebox(0,     0)[t]{$\bullet$}}

\put(91,      -1){\line(0,     1){34}}

\put(101,      33){$q ^{}$}

\put(101,      16){$q ^{-1}$}

\put(219,        -1)  {by GDD $12$ of Row $18$. }
\end{picture}$\\ \\
 It is quasi-affine    by Row 21 or  Lemma \ref {2.63}.
{\rm (c)} There are not any other cases by Lemma \ref {2.13}.

 {\rm (ii)}   Adding on  Vertex 2, 3, 5. It is not  It is quasi-affine    since it contains a proper subGDD which is not an arithmetic GDD.

 {\rm (iii)} Omitting Vertex 1 and adding on   Vertex 4.  It    is empty  by  Lemma \ref {2.2} {\rm (ii)}.
            \subsection  {Quasi-affine  GDDs  over
GDD $21$ of Row $18$ in Table C }\label {sub3.2}
{\rm (i)} Omitting Vertex 1 and  adding on   Vertex 6.
 \\ \\

 $\ \ \ \ \  \ \   \begin{picture}(100,      15)  \put(-45,      -1){}

\put(-45,      -1){ {\rm (a)}}
\put(27,      1){\makebox(0,     0)[t]{$\bullet$}}
\put(60,      1){\makebox(0,      0)[t]{$\bullet$}}
\put(93,     1){\makebox(0,     0)[t]{$\bullet$}}

\put(126,      1){\makebox(0,    0)[t]{$\bullet$}}

\put(159,      1){\makebox(0,    0)[t]{$\bullet$}}
\put(194,      1){\makebox(0,    0)[t]{$\bullet$}}

\put(126,      1){\makebox(0,    0)[t]{$\bullet$}}
\put(28,      -1){\line(1,      0){33}}
\put(61,      -1){\line(1,      0){30}}
\put(94,     -1){\line(1,      0){30}}
\put(126,      -1){\line(1,      0){33}}
\put(159,      -1){\line(1,      0){33}}

\put(22,     -20){$q ^{-1}$}
\put(40,      -15){$q ^{}$}
\put(58,      -20){${-1}$}
\put(75,      -15){$q ^{-1}$}

\put(91,      -20){$q ^{}$}

\put(102,    -15){$q ^{-1}$}

\put(124,      10){$q ^{}$}

\put(135,     5){$q ^{-1}$}

\put(161,      10){$q$}

\put(168,     5){$q ^{-1}$}

\put(194,      10){$q$}

\put(91,    32){\makebox(0,     0)[t]{$\bullet$}}

\put(91,      -1){\line(0,     1){34}}

\put(101,      33){$q ^{}$}

\put(101,      16){$q ^{-1}$}

\put(200,        -1)  { GDD $5$ of Row $17$  in Table C. }
\end{picture}$\\ \\
Omitting  Vertex 4  in whole is a simple chain.
  It is quasi-affine    by Row 21 or  Lemma \ref {2.63}.\\ \\

 $\ \ \ \ \  \ \   \begin{picture}(100,      15)  \put(-45,      -1){}

\put(-45,      -1){ {\rm (b)}}
\put(27,      1){\makebox(0,     0)[t]{$\bullet$}}
\put(60,      1){\makebox(0,      0)[t]{$\bullet$}}
\put(93,     1){\makebox(0,     0)[t]{$\bullet$}}

\put(126,      1){\makebox(0,    0)[t]{$\bullet$}}

\put(159,      1){\makebox(0,    0)[t]{$\bullet$}}
\put(194,      1){\makebox(0,    0)[t]{$\bullet$}}

\put(126,      1){\makebox(0,    0)[t]{$\bullet$}}
\put(28,      -1){\line(1,      0){33}}
\put(61,      -1){\line(1,      0){30}}
\put(94,     -1){\line(1,      0){30}}
\put(126,      -1){\line(1,      0){33}}
\put(159,      -1){\line(1,      0){33}}

\put(22,     -20){$q ^{-1}$}
\put(40,      -15){$q ^{}$}
\put(58,      -20){${-1}$}
\put(75,      -15){$q ^{-1}$}

\put(91,      -20){$q ^{}$}

\put(102,    -15){$q ^{-1}$}

\put(124,      10){$q ^{}$}

\put(135,     5){$q ^{-1}$}

\put(161,      10){$q$}

\put(168,     5){$q ^{-1}$}

\put(194,      10){$-1$}

\put(91,    32){\makebox(0,     0)[t]{$\bullet$}}

\put(91,      -1){\line(0,     1){34}}

\put(101,      33){$q ^{}$}

\put(101,      16){$q ^{-1}$}

\put(210,        -1)  { GDD $9$ of Row $18$  in Table C. }
\end{picture}$\\ \\
Omitting  Vertex 4  in whole is a simple chain.
  It is quasi-affine    by Row 21 or  Lemma \ref {2.63}.
{\rm (c)} There are not any other cases by Lemma \ref {2.13}.

 {\rm (ii)}    Adding on  Vertex 2, 3, 5. It is not quasi-affine since it contains a proper subGDD which is not an arithmetic GDD.

   {\rm (iii)}   Omitting Vertex 1 and adding on   Vertex 4.  It    is empty  by Lemma \ref {2.2} {\rm (ii)}.
             \subsection  {Quasi-affine  GDDs  over
GDD $1$ of Row $19$ in Table C } \label {sub3.2}
     {\rm (i)}   Omitting Vertex 1 and adding on   Vertex 6.  It    is empty  by Lemma \ref {2.63}.
                 \subsection  {Quasi-affine  GDDs  over
GDD $2$ of Row $19$ in Table C }\label {sub3.2}
    {\rm (i)}   Omitting Vertex 1 and adding on   Vertex 6.  It    is empty  by Lemma \ref {2.63}.

     {\rm (ii)} Adding on  Vertex 2,  3, 4, 5.  It is not quasi-affine since it contains a proper subGDD which is not an arithmetic GDD.
             \subsection  {Quasi-affine  GDDs  over
GDD $3$ of Row $19$ in Table C } \label {sub3.2}
     {\rm (i)}   Adding on   Vertex 3, 4, 5, 6. It is not quasi-affine since it contains a proper subGDD which is not an arithmetic GDD.

       {\rm (ii)}   Omitting  Vertex 6 and adding on  Vertex 2.  It    is empty  by Lemma \ref {2.13}.
            \subsection  {Quasi-affine  GDDs  over
GDD $4$ of Row $19$ in Table C } \label {sub3.2}
 {\rm (i)}   Omitting Vertex 1 and adding on  Vertex  6.
 It    is empty  by Lemma \ref {2.63}.
             \subsection  {Quasi-affine  GDDs  over
GDD $5$ of Row $19$ in Table C }\label {sub3.2}
 {\rm (i)}   Omitting Vertex 1 and adding on  Vertex  6.
    It    is empty  by   Row 19.

     {\rm (ii)}  Adding on  Vertex 2, 3, 4, 5. It is not quasi-affine since it contains a proper subGDD which is not an arithmetic GDD.
         \subsection  {Quasi-affine  GDDs  over
GDD $6$ of Row $19$ in Table C } \label {sub3.2}
 {\rm (i)}   Omitting Vertex 1 and adding on   Vertex 6.
    It    is empty  by  Row 19.

     {\rm (ii)}  Adding on  Vertex 2, 3, 4, 5. It is not quasi-affine since it contains a proper subGDD which is not an arithmetic GDD.
          \subsection  {Quasi-affine  GDDs  over
GDD $7$ of Row $19$ in Table C } \label {sub3.2}
 {\rm (i)}   Omitting Vertex 1 and adding on   Vertex 6.

     {\rm (ii)}  Adding on  Vertex 2, 3, 5. It is not quasi-affine since it contains a proper subGDD which is not an arithmetic GDD.

 {\rm (iii)}  Omitting Vertex 1 and adding on   Vertex 4.\\ \\ \\ \\

{\ }

{\ }

{\ }

{\ }

 $\ \ \ \ \  \ \   
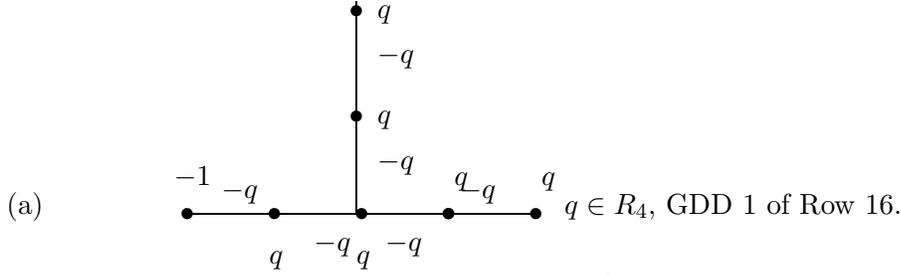
\begin{picture}(100,      15)  \put(-45,      -1){}

\put(-45,      -1){ {\rm (a)}}
\put(27,      1){\makebox(0,     0)[t]{$\bullet$}}
\put(60,      1){\makebox(0,      0)[t]{$\bullet$}}
\put(93,     1){\makebox(0,     0)[t]{$\bullet$}}

\put(126,      1){\makebox(0,    0)[t]{$\bullet$}}

\put(159,      1){\makebox(0,    0)[t]{$\bullet$}}

\put(126,      1){\makebox(0,    0)[t]{$\bullet$}}
\put(28,      -1){\line(1,      0){33}}
\put(61,      -1){\line(1,      0){30}}
\put(94,     -1){\line(1,      0){30}}
\put(126,      -1){\line(1,      0){33}}

\put(22,     10){$-1 ^{}$}
\put(40,      5){$-q^{}$}
\put(58,      -20){${q}$}
\put(75,      -15){$-q ^{}$}

\put(91,      -20){${q}$}

\put(102,    -15){$-q ^{}$}

\put(128,      10){$q ^{}$}

\put(130,     5){$-q ^{}$}

\put(161,      10){$q^{}$}

\put(91,    38){\makebox(0,     0)[t]{$\bullet$}}

\put(91,    78){\makebox(0,     0)[t]{$\bullet$}}

\put(91,      -1){\line(0,     1){80}}

\put(99,      16){$-q ^{}$}

\put(99,      33){$q ^{}$}

\put(99,      56){$-q ^{}$}

\put(99,      73){$q ^{}$}

\put(166,        -1)  {  $q \in R_4,$   GDD $1$ of Row $16$. }
\end{picture}$\\ \\
Omitting  Vertex 7  in whole is an arithmetic GDD   by   GDD $7$ of Row $19$.
 It is quasi-affine    by   Lemma \ref {2.13}.
\\ \\ \\ \\

 $\ \ \ \ \  \ \   
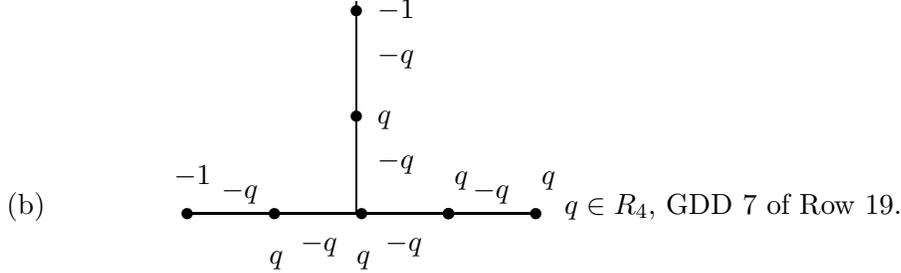
\begin{picture}(100,      15)  \put(-45,      -1){}

\put(-45,      -1){ {\rm (b)}}
\put(27,      1){\makebox(0,     0)[t]{$\bullet$}}
\put(60,      1){\makebox(0,      0)[t]{$\bullet$}}
\put(93,     1){\makebox(0,     0)[t]{$\bullet$}}

\put(126,      1){\makebox(0,    0)[t]{$\bullet$}}

\put(159,      1){\makebox(0,    0)[t]{$\bullet$}}

\put(126,      1){\makebox(0,    0)[t]{$\bullet$}}
\put(28,      -1){\line(1,      0){33}}
\put(61,      -1){\line(1,      0){30}}
\put(94,     -1){\line(1,      0){30}}
\put(126,      -1){\line(1,      0){33}}

\put(22,     10){$-1 ^{}$}
\put(40,      5){$-q^{}$}
\put(58,      -20){${q}$}
\put(70,      -15){$-q ^{}$}

\put(91,      -20){${q}$}

\put(102,    -15){$-q ^{}$}

\put(128,      10){$q ^{}$}

\put(135,     5){$-q ^{}$}

\put(161,      10){$q^{}$}

\put(91,    38){\makebox(0,     0)[t]{$\bullet$}}

\put(91,    78){\makebox(0,     0)[t]{$\bullet$}}

\put(91,      -1){\line(0,     1){80}}

\put(99,      16){$-q ^{}$}

\put(99,      33){$q ^{}$}

\put(99,      56){$-q ^{}$}

\put(99,      73){$-1 ^{}$}

\put(166,        -1)  {   $q \in R_4,$   GDD $7$ of Row $19$. }
\end{picture}$\\ \\
Omitting  Vertex 7  in whole is not  an arithmetic GDD   by   Row $19$.
{\rm (c)} There are not any other cases by Lemma \ref {2.13}.
           \subsection  {Quasi-affine  GDDs  over
GDD $1$ of Row $20$ in Table C }\label {sub3.2}
 {\rm (i)}   Omitting Vertex 1 and adding on  Vertex  7.\\ \\

 $\ \ \ \ \  \ \   
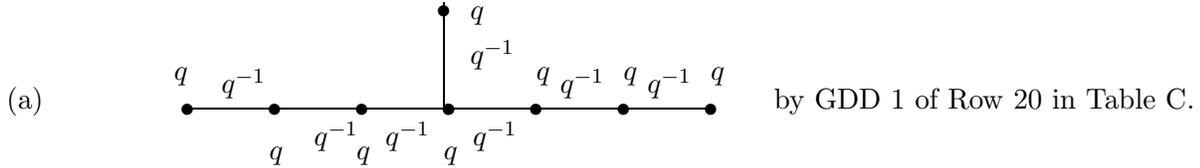
\begin{picture}(100,      15)  \put(-45,      -1){}
\put(-45,      -1){ {\rm (a)}}
\put(27,      1){\makebox(0,     0)[t]{$\bullet$}}
\put(60,      1){\makebox(0,      0)[t]{$\bullet$}}
\put(93,     1){\makebox(0,     0)[t]{$\bullet$}}

\put(126,      1){\makebox(0,    0)[t]{$\bullet$}}

\put(159,      1){\makebox(0,    0)[t]{$\bullet$}}
\put(192,      1){\makebox(0,    0)[t]{$\bullet$}}
\put(225,      1){\makebox(0,    0)[t]{$\bullet$}}

\put(126,      1){\makebox(0,    0)[t]{$\bullet$}}
\put(28,      -1){\line(1,      0){33}}
\put(61,      -1){\line(1,      0){30}}
\put(94,     -1){\line(1,      0){30}}
\put(126,      -1){\line(1,      0){33}}
\put(159,      -1){\line(1,      0){33}}
\put(192,      -1){\line(1,      0){33}}

\put(22,     10){$q$}
\put(40,      5){$q^{-1}$}
\put(58,      -20){$q^{}$}
\put(75,      -15){$q^{-1}$}

\put(91,      -20){$q^{}$}

\put(102,     -15){$q^{-1}$}

\put(124,      -20){$q^{}$}
\put(135,     -15){$q^{-1}$}
\put(159,      10){$q^{}$}

\put(168,     5){$q^{-1}$}
\put(192,      10){$q^{}$}

\put(201,     5){$q ^{-1}$}
\put(225,      10){$q ^{}$}

\put(124,    38){\makebox(0,     0)[t]{$\bullet$}}

\put(124,      -1){\line(0,     1){40}}

\put(134,      33){$q^{}$}

\put(134,      16){$q^{-1}$}

\put(249,        -1)  {by GDD $1$ of Row $20$  in Table C. }
\end{picture}$\\ \\
 It is quasi-affine.\\ \\\

 $\ \ \ \ \  \ \   
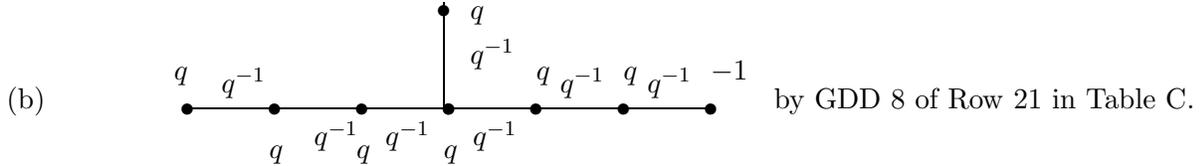
\begin{picture}(100,      15)  \put(-45,      -1){}
\put(-45,      -1){ {\rm (b)}}
\put(27,      1){\makebox(0,     0)[t]{$\bullet$}}
\put(60,      1){\makebox(0,      0)[t]{$\bullet$}}
\put(93,     1){\makebox(0,     0)[t]{$\bullet$}}

\put(126,      1){\makebox(0,    0)[t]{$\bullet$}}

\put(159,      1){\makebox(0,    0)[t]{$\bullet$}}
\put(192,      1){\makebox(0,    0)[t]{$\bullet$}}
\put(225,      1){\makebox(0,    0)[t]{$\bullet$}}

\put(126,      1){\makebox(0,    0)[t]{$\bullet$}}
\put(28,      -1){\line(1,      0){33}}
\put(61,      -1){\line(1,      0){30}}
\put(94,     -1){\line(1,      0){30}}
\put(126,      -1){\line(1,      0){33}}
\put(159,      -1){\line(1,      0){33}}
\put(192,      -1){\line(1,      0){33}}

\put(22,     10){$q$}
\put(40,      5){$q^{-1}$}
\put(58,      -20){$q^{}$}
\put(75,      -15){$q^{-1}$}

\put(91,      -20){$q^{}$}

\put(102,     -15){$q^{-1}$}

\put(124,      -20){$q^{}$}
\put(135,     -15){$q^{-1}$}
\put(159,      10){$q^{}$}

\put(168,     5){$q^{-1}$}
\put(192,      10){$q^{}$}

\put(201,     5){$q ^{-1}$}
\put(225,      10){$-1 ^{}$}

\put(124,    38){\makebox(0,     0)[t]{$\bullet$}}

\put(124,      -1){\line(0,     1){40}}

\put(134,      33){$q^{}$}

\put(134,      16){$q^{-1}$}

\put(249,        -1)  {by GDD $8$ of Row $21$  in Table C. }
\end{picture}$\\ \\
 It is repeated.
{\rm (c)} There are not any other cases by Lemma \ref {2.13}.

 {\rm (ii)}   Adding on  Vertex 2, 3, 4, 5, 6. It is not quasi-affine since it contains a proper subGDD which is not an arithmetic GDD.
            \subsection  {Quasi-affine  GDDs  over
GDD $1$ of Row $21$ in Table C } \label {sub3.2}
{\rm (i)}   Omitting Vertex 1 and adding on   Vertex 7.
It    is empty  by Row 21 or Lemma \ref {2.63}.
            \subsection  {Quasi-affine  GDDs  over
GDD $2$ of Row $21$ in Table C }\label {sub3.2}
 {\rm (i)}   Omitting Vertex 1 and adding on   Vertex 7.
It    is empty  by Row 21 or  Lemma \ref {2.63}.

 {\rm (ii)}   Adding on  Vertex 2, 3, 4, 5, 6. It is not  It is quasi-affine    since it contains a proper subGDD which is not an arithmetic GDD.
 \subsection  {Quasi-affine  GDDs  over
GDD $3$ of Row $21$ in Table C }\label {sub3.2}
  {\rm (i)}    Adding on   Vertex 3, 4, 5.
It is not quasi-affine since it contains a proper subGDD which is not an arithmetic GDD.

 {\rm (ii)}   Omitting Vertex 1 and adding on   Vertex 7.   It    is empty  by Row 21 or Lemma \ref {2.63}.

 {\rm (iii)}   Omitting  Vertex 6 and adding on  Vertex 1.

  {\rm (iv)}   Omitting  Vertex 6 and adding on  Vertex 2.\\ \\

 $\ \ \ \ \  \ \   
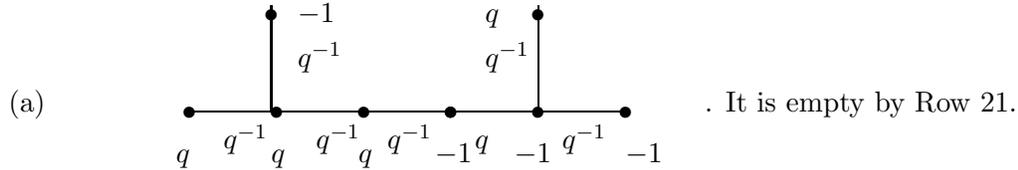
\begin{picture}(100,      15)  \put(-45,      -1){} \put(-45,      -1){ {\rm (a)}}
\put(27,      1){\makebox(0,     0)[t]{$\bullet$}}
\put(60,      1){\makebox(0,      0)[t]{$\bullet$}}
\put(93,     1){\makebox(0,     0)[t]{$\bullet$}}

\put(126,      1){\makebox(0,    0)[t]{$\bullet$}}

\put(159,      1){\makebox(0,    0)[t]{$\bullet$}}
\put(192,      1){\makebox(0,    0)[t]{$\bullet$}}

\put(126,      1){\makebox(0,    0)[t]{$\bullet$}}
\put(28,      -1){\line(1,      0){33}}
\put(61,      -1){\line(1,      0){30}}
\put(94,     -1){\line(1,      0){30}}
\put(126,      -1){\line(1,      0){33}}
\put(159,      -1){\line(1,      0){33}}

\put(22,     -20){$q$}
\put(40,      -15){$q ^{-1}$}
\put(58,      -20){$q$}
\put(75,      -15){$q ^{-1}$}

\put(91,      -20){$q$}

\put(102,     -15){$q ^{-1}$}

\put(120,      -20){$-1$}
\put(135,     -15){$q ^{}$}

\put(150,     -20){${-1}$}

\put(168,     -15){$q ^{-1}$}
\put(192,      -20){$-1 ^{}$}

\put(58,    38){\makebox(0,     0)[t]{$\bullet$}}

\put(58,      -1){\line(0,     1){40}}

\put(68,      33){$-1$}

\put(68,      16){$q ^{-1}$}

\put(159,    38){\makebox(0,     0)[t]{$\bullet$}}

\put(159,      -1){\line(0,     1){40}}

\put(139,      33){$q$}

\put(139,      16){$q ^{-1}$}

\put(222,        -1)  {. It    is empty  by Row 21. }
\end{picture}$\\ \\
{\rm (b)} There are not any other cases by Lemma \ref {2.63}.
            \subsection  {Quasi-affine  GDDs  over
GDD $4$ of Row $21$ in Table C }\label {sub3.2}
{\rm (i)}  Adding on   Vertex 3, 4, 5.
It is not quasi-affine since it contains a proper subGDD which is not an arithmetic GDD.

   {\rm (ii)}   Omitting Vertex 1 and adding on   Vertex 6.  It    is empty  by Row 21 or  Lemma \ref {2.2} {\rm (ii)}.

   {\rm (iii)}   Omitting Vertex 1 and adding on   Vertex 7.  It    is empty  by Row 21 or  Lemma \ref {2.63}.

  {\rm (iv)}   Omitting Vertex 1 and adding on  Vertex 2.\\ \\

 $\ \ \ \ \  \ \   
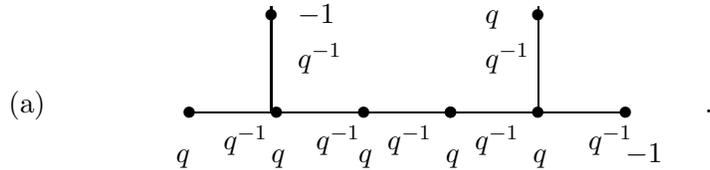
\begin{picture}(100,      15)  \put(-45,      -1){} \put(-45,      -1){ {\rm (a)}}
\put(27,      1){\makebox(0,     0)[t]{$\bullet$}}
\put(60,      1){\makebox(0,      0)[t]{$\bullet$}}
\put(93,     1){\makebox(0,     0)[t]{$\bullet$}}

\put(126,      1){\makebox(0,    0)[t]{$\bullet$}}

\put(159,      1){\makebox(0,    0)[t]{$\bullet$}}
\put(192,      1){\makebox(0,    0)[t]{$\bullet$}}

\put(126,      1){\makebox(0,    0)[t]{$\bullet$}}
\put(28,      -1){\line(1,      0){33}}
\put(61,      -1){\line(1,      0){30}}
\put(94,     -1){\line(1,      0){30}}
\put(126,      -1){\line(1,      0){33}}
\put(159,      -1){\line(1,      0){33}}

\put(22,     -20){$q$}
\put(40,      -15){$q ^{-1}$}
\put(58,      -20){$q$}
\put(75,      -15){$q ^{-1}$}

\put(91,      -20){$q$}

\put(102,     -15){$q ^{-1}$}

\put(124,      -20){$q$}
\put(135,     -15){$q ^{-1}$}

\put(157,     -20){$q ^{}$}

\put(178,     -15){$q ^{-1}$}
\put(192,      -20){$-1 ^{}$}

\put(58,    38){\makebox(0,     0)[t]{$\bullet$}}

\put(58,      -1){\line(0,     1){40}}

\put(68,      33){$-1$}

\put(68,      16){$ q ^{-1}$}

\put(159,    38){\makebox(0,     0)[t]{$\bullet$}}

\put(159,      -1){\line(0,     1){40}}

\put(139,      33){$q$}

\put(139,      16){$q ^{-1}$}

\put(222,        -1)  {. }
\end{picture}$\\ \\  Omitting  Vertex 7 in whole is not arithmetic by Row 21.
{\rm (b)} There are not any other cases by Lemma \ref {2.13}.
             \subsection  {Quasi-affine  GDDs  over
GDD $5$ of Row $21$ in Table C }\label {sub3.2}
 {\rm (i)}   Adding on  Vertex 2, 3, 4, 5, 6. It is not quasi-affine since it contains a proper subGDD which is not an arithmetic GDD.

 {\rm (ii)} Omitting Vertex 1 and adding on   Vertex 7.  It    is empty  by Row 21.
             \subsection  {Quasi-affine  GDDs  over
GDD $6$ of Row $21$ in Table C }\label {sub3.2}
 {\rm (i)}  Adding on  Vertex 2, 3, 4, 5, 6. It is not quasi-affine since it contains a proper subGDD which is not an arithmetic GDD.

 {\rm (ii)} Omitting Vertex 1 and adding on   Vertex 7.  It    is empty  by Row 21.
             \subsection  {Quasi-affine  GDDs  over
GDD $7$ of Row $21$ in Table C }\label {sub3.2}
 {\rm (i)}  Adding on  Vertex 2, 3, 4, 5, 6. It is not quasi-affine since it contains a proper subGDD which is not an arithmetic GDD.

 {\rm (ii)} Omitting Vertex 1 and adding on  Vertex  7.  It    is empty  by Row 21.
             \subsection  {Quasi-affine  GDDs  over
GDD $8$ of Row $21$ in Table C }\label {sub3.2}
 {\rm (i)}  Adding on  Vertex 2, 3, 4, 5, 6. It is not quasi-affine since it contains a proper subGDD which is not an arithmetic GDD.

 {\rm (ii)} Omitting Vertex 1 and adding on   Vertex 7. \\ \\

 $\ \ \ \ \  \ \   
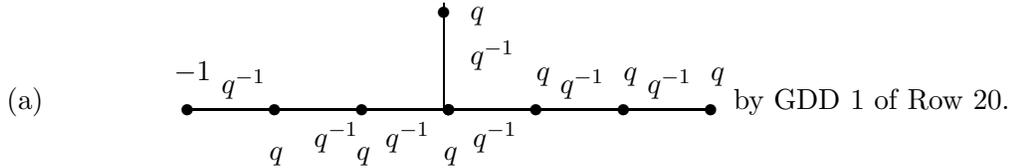
\begin{picture}(100,      15)  \put(-45,      -1){}
\put(-45,      -1){ {\rm (a)}}
\put(27,      1){\makebox(0,     0)[t]{$\bullet$}}
\put(60,      1){\makebox(0,      0)[t]{$\bullet$}}
\put(93,     1){\makebox(0,     0)[t]{$\bullet$}}

\put(126,      1){\makebox(0,    0)[t]{$\bullet$}}

\put(159,      1){\makebox(0,    0)[t]{$\bullet$}}
\put(192,      1){\makebox(0,    0)[t]{$\bullet$}}
\put(225,      1){\makebox(0,    0)[t]{$\bullet$}}

\put(126,      1){\makebox(0,    0)[t]{$\bullet$}}
\put(28,      -1){\line(1,      0){33}}
\put(61,      -1){\line(1,      0){30}}
\put(94,     -1){\line(1,      0){30}}
\put(126,      -1){\line(1,      0){33}}
\put(159,      -1){\line(1,      0){33}}
\put(192,      -1){\line(1,      0){33}}

\put(22,     10){$-1$}
\put(40,      5){$q^{-1}$}
\put(58,      -20){$q^{}$}
\put(75,      -15){$q^{-1}$}

\put(91,      -20){$q^{}$}

\put(102,     -15){$q^{-1}$}

\put(124,      -20){$q^{}$}
\put(135,     -15){$q^{-1}$}
\put(159,      10){$q^{}$}

\put(168,     5){$q^{-1}$}
\put(192,      10){$q^{}$}

\put(201,     5){$q ^{-1}$}
\put(225,      10){$q ^{}$}

\put(124,    38){\makebox(0,     0)[t]{$\bullet$}}

\put(124,      -1){\line(0,     1){40}}

\put(134,      33){$q^{}$}

\put(134,      16){$q^{-1}$}

\put(230,  -1)  { by  GDD $1$ of Row $20$. }
\end{picture}$\\ \\
 It is quasi-affine    by Row 22.
\\ \\

 $\ \ \ \ \  \ \   
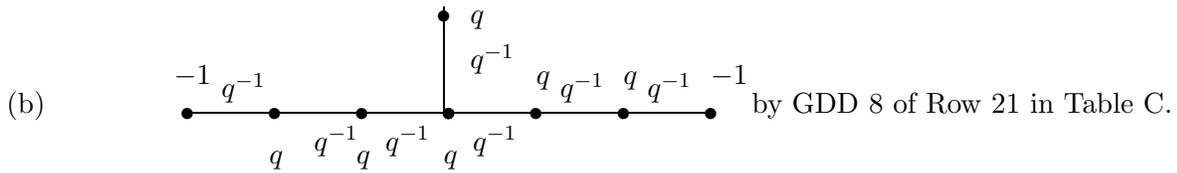
\begin{picture}(100,      15)  \put(-45,      -1){}
\put(-45,      -1){ {\rm (b)}}
\put(27,      1){\makebox(0,     0)[t]{$\bullet$}}
\put(60,      1){\makebox(0,      0)[t]{$\bullet$}}
\put(93,     1){\makebox(0,     0)[t]{$\bullet$}}

\put(126,      1){\makebox(0,    0)[t]{$\bullet$}}

\put(159,      1){\makebox(0,    0)[t]{$\bullet$}}
\put(192,      1){\makebox(0,    0)[t]{$\bullet$}}
\put(225,      1){\makebox(0,    0)[t]{$\bullet$}}

\put(126,      1){\makebox(0,    0)[t]{$\bullet$}}
\put(28,      -1){\line(1,      0){33}}
\put(61,      -1){\line(1,      0){30}}
\put(94,     -1){\line(1,      0){30}}
\put(126,      -1){\line(1,      0){33}}
\put(159,      -1){\line(1,      0){33}}
\put(192,      -1){\line(1,      0){33}}

\put(22,     10){$-1$}
\put(40,      5){$q^{-1}$}
\put(58,      -20){$q^{}$}
\put(75,      -15){$q^{-1}$}

\put(91,      -20){$q^{}$}

\put(102,     -15){$q^{-1}$}

\put(124,      -20){$q^{}$}
\put(135,     -15){$q^{-1}$}
\put(159,      10){$q^{}$}

\put(168,     5){$q^{-1}$}
\put(192,      10){$q^{}$}

\put(201,     5){$q ^{-1}$}
\put(225,      10){$-1 ^{}$}

\put(124,    38){\makebox(0,     0)[t]{$\bullet$}}

\put(124,      -1){\line(0,     1){40}}

\put(134,      33){$q^{}$}

\put(134,      16){$q^{-1}$}

\put(237,        -1)  { by GDD $8$ of Row $21$  in Table C. }
\end{picture}$\\ \\
Omitting  Vertex 5  in whole is a simple chain.  It is quasi-affine.
{\rm (c)} There are not any other cases by Lemma \ref {2.13}.
          \subsection  {Quasi-affine  GDDs  over
GDD $1$ of Row $22$ in Table C }\label {ss99}
 {\rm (i)}    Adding on  Vertex 2, 3, 4, 5, 6, 7. It is not quasi-affine since it contains a proper subGDD which is not an arithmetic GDD.

  {\rm (ii)} Omitting Vertex 1 and adding on   Vertex 8.
It is empty by Row 22.

\section {Appendix}

Given a generalized Cartan matrix $A$ we can obtain a Kac-Moody Lie algebra $ \mathfrak g (A)$. If $A$ is a Cartan matrix, then $L$ is a semi-simple Lie algebra;   If $A$ is a quasi-affine  generalized Cartan matrix, then  $ \mathfrak g (A)$ is a quasi-affine  Lie algebra, which is infinite dimensional; Given a braided matrix $Q$,  its braided vector $V$  we can obtain a Nichols  algebra $\mathfrak B(V)$ and  Nichols Lie braided algebra $\mathfrak L (V)$. If the GDD is an arithmetic GDD  which   fixed
parameter is of finite order, then $\mathfrak B(V)$ and $\mathfrak L (V)$  are  finite dimensional;   If the GDD is  quasi-affine   then $\mathfrak B(V)$ and  $\mathfrak L (V)$  are infinite dimensional (see \cite {WZZ15a, WZZ15b, WWZZ}).

If a GDD is obtained  by relabelling   other GDD, then we view the two GDD are the same since their YD modules are isomorphic.
For example, \\

 $\begin{picture}(100,      15)

\put(60,      1){\makebox(0,      0)[t]{$\bullet$}}

\put(28,      -1){\line(1,      0){33}}
\put(27,      1){\makebox(0,     0)[t]{$\bullet$}}

\put(-14,      1){\makebox(0,     0)[t]{$\bullet$}}

\put(-14,     -1){\line(1,      0){50}}

\put(-18,     10){$q$}
\put(0,      5){$q ^{-1}$}
\put(22,     10){$q$}
\put(40,      5){$q ^{-1}$}

\put(58,      10){$q ^{-3}$}

\put(122,      -1){and }

  \ \ \ \ \ \ \ \ \ \ \ \ \ \ \ \ \ \ \ {$ q \in R_9$ }
\end{picture}$\ \ \ \ \ \ \ \ \ \ \ \
  \ \ \ \ \ \
 $\begin{picture}(100,      15)

\put(60,      1){\makebox(0,      0)[t]{$\bullet$}}

\put(28,      -1){\line(1,      0){33}}
\put(27,      1){\makebox(0,     0)[t]{$\bullet$}}

\put(-14,      1){\makebox(0,     0)[t]{$\bullet$}}

\put(-14,     -1){\line(1,      0){50}}

\put(-18,     10){ $q ^{-3}$}
\put(0,      5){$q ^{-1}$}
\put(22,     10){$q$}
\put(40,      5){$q ^{-1}$}

\put(58,      10){$q$}

  \ \ \ \ \ \ \ \ \ \ \ \ \ \ \ \ \ \ \ {$ q \in R_9$\ \ \  are the same. }
\end{picture}$
\subsection {Relation between quasi-affine  GDDs and affine  matrices }
If for any $1\le i\not =j \le n$, there exist $b_{ij} \in - \mathbb N_0$ such that $\widetilde{q}_{ij} = q_{ii} ^{b_{ij}}$, then GDD is said to be of Cartan type.
Let $a_{ij}$ be the maximal number in $\{ b_{ij} \mid $ $\widetilde{q}_{ij} = q_{ii} ^{b_{ij}}$, $b_{ij} \in - \mathbb N_0\}$. $a_{ii} =: 2.$ $A=:(a_{ij})_{n \times n}$ is called a braiding exponential matrix(see \cite {AS00}).
It is clear that $A$ is a generalized Cartan matrix. If $A$ is affine, then GDD is called an affine GDD.

\begin {Proposition} \label {a5.1} (\cite [Th. 4]{He06b}) If GDD is of  Cartan type with braiding exponential matrix $A$,
then GDD is arithmetic if and only if $A$ is a Cartan matrix.

\end {Proposition}

{\bf Proof.} It follows from \cite [Pro. 4.7]{Ka90} and \cite [Th. 4]{He06b}.
\hfill $\Box$

\begin {Proposition} \label {a5.2} All affine GDDs are listed:\\

$\ \ \ \ \  \ \ \begin{picture}(100,      15)  \put(-45,      -1){$A^{(1)}_1$ }

\put(60,      1){\makebox(0,      0)[t]{$\bullet$}}

\put(28,      -1){\line(1,      0){33}}
\put(27,      1){\makebox(0,     0)[t]{$\bullet$}}

%\put(-14,      1){\makebox(0,     0)[t]{$\bullet$}}

%\put(-14,     -1){\line(1,      0){50}}

\put(22,     10){$q ^{}$}
\put(38,      5){$q^{-2}$}

\put(60,      10){$q ^{}$}

  \ \ \ \ \ \ \ \ \ \ \ \ \ \ \ \ \ \ \ {$q^2 \not=1$. }
\end{picture}$\\ \\ \\ \\ \\

$\ \ \ \ \  \ \ \begin{picture}(130,      15)

\put(-45,      -1){$A^{(1)}_N, N\ge 2$ }

\put(27,      1){\makebox(0,     0)[t]{$\bullet$}}
\put(60,      1){\makebox(0,      0)[t]{$\bullet$}}
\put(93,     1){\makebox(0,     0)[t]{$\bullet$}}

%\put(126,      1){\makebox(0,    0)[t]{$\bullet$}}

\put(159,      1){\makebox(0,    0)[t]{$\bullet$}}
\put(192,      1){\makebox(0,    0)[t]{$\bullet$}}

%\put(126,      1){\makebox(0,    0)[t]{$\bullet$}}
\put(28,      -1){\line(1,      0){33}}
\put(61,      -1){\line(1,      0){30}}
% \put(94,     -1){\line(1,      0){30}}
%% \put(126,      -1){\line(1,      0){33}}
\put(159,      -1){\line(1,      0){33}}

\put(130,     1){\makebox(0,     0)[t]{$\cdots\cdots\cdots\cdots$}}

\put(22,     10){$q ^{}$}
\put(40,      5){$q ^{-1}$}
\put(58,      10){$q$}
\put(74,      5){$q ^{-1}$}

\put(91,      10){$q ^{}$}

%\put(102,     5){$q ^{-1}$}

%\put(124,      10){$q ^{-1}$}
%\put(135,     5){$q ^{-1}$}
\put(159,      10){$q ^{}$}

\put(168,     5){$q ^{-1}$}
\put(192,      10){$q$}

\put(112,     85){\makebox(0,     0)[t]{$\bullet$}}

\put(28,      -1){\line(1,      1){85}}

\put(192,      -1){\line(-1,      1){85}}

\put(100,      65){$q$}
\put(50,      35){$q^{-1}$}
\put(155,      35){$q ^{-1}$}

\put(210,        -1)  {$q \not=1$. }
\end{picture}$\\ \\

$\ \ \ \ \  \ \ \begin{picture}(130,      15)

\put(-45,      -1){$B^{(1)}_N, N\ge 3$  }
\put(27,      1){\makebox(0,     0)[t]{$\bullet$}}
\put(60,      1){\makebox(0,      0)[t]{$\bullet$}}
\put(93,     1){\makebox(0,     0)[t]{$\bullet$}}

\put(126,      1){\makebox(0,    0)[t]{$\bullet$}}

\put(159,      1){\makebox(0,    0)[t]{$\bullet$}}

\put(126,      1){\makebox(0,    0)[t]{$\bullet$}}
\put(28,      -1){\line(1,      0){33}}
\put(61,      -1){\line(1,      0){30}}
%\put(94,     -1){\line(1,      0){30}}
\put(126,      -1){\line(1,      0){33}}

\put(104,     1){\makebox(0,     0)[t]{$\cdots\cdots$}}

\put(22,     10){$q$}
\put(40,      5){$q^{-1}$}
\put(58,      -20){$q$}
\put(75,      -15){$q ^{-1}$}

\put(91,      10){$q ^{}$}

%\put(102,     5){$q ^{-1}$}

\put(128,      10){$q ^{}$}

\put(135,     5){$q ^{-2}$}

\put(161,      10){$q ^{2}$}

\put(58,    38){\makebox(0,     0)[t]{$\bullet$}}

\put(58,      -1){\line(0,     1){40}}

\put(68,      33){$q ^{}$}

\put(68,      16){$q ^{-1}$}

\put(176,        -1)  { $q ^{2} \not= 1$. }
\end{picture}$\\ \\

$\ \ \ \ \  \ \ \begin{picture}(100,      15)

\put(-45,      -1){ $C^{(1)}_N, N\ge 2$  }

\put(27,      1){\makebox(0,     0)[t]{$\bullet$}}
\put(60,      1){\makebox(0,      0)[t]{$\bullet$}}
\put(93,     1){\makebox(0,     0)[t]{$\bullet$}}
\put(159,      1){\makebox(0,      0)[t]{$\bullet$}}
\put(192,     1){\makebox(0,      0)[t]{$\bullet$}}
\put(28,      -1){\line(1,      0){30}}
\put(61,      -1){\line(1,      0){30}}
\put(130,     1){\makebox(0,     0)[t]{$\cdots\cdots\cdots\cdots$}}
\put(160,     -1){\line(1,      0){30}}
\put(22,     -15){}
\put(58,      -15){}
\put(91,      -15){}
\put(157,      -15){$$}
\put(191,      -15){$$}
\put(22,     10){$q$}

\put(40,      5){$q^{-2}$}

\put(58,      10){$q^{2}$}
\put(73,      5){$q^{-2}$}

\put(91,      10){$q^{2}$}

\put(157,      10){$q^{2}$}

\put(172,     5){$q^{-2}$}
\put(191,      10){$q$}

\put(210,        -1) { $q ^{2} \not= 1$.}
\end{picture}$\\
\\

$\ \ \ \ \  \ \ \begin{picture}(130,      15)\put(-45,      -1){  $D^{(1)}_N, N\ge 4$   }
\put(27,      1){\makebox(0,     0)[t]{$\bullet$}}
\put(60,      1){\makebox(0,      0)[t]{$\bullet$}}
\put(93,     1){\makebox(0,     0)[t]{$\bullet$}}

\put(126,      1){\makebox(0,    0)[t]{$\bullet$}}

\put(159,      1){\makebox(0,    0)[t]{$\bullet$}}
\put(192,      1){\makebox(0,    0)[t]{$\bullet$}}

\put(126,      1){\makebox(0,    0)[t]{$\bullet$}}
\put(28,      -1){\line(1,      0){33}}
\put(61,      -1){\line(1,      0){30}}
%\put(94,     -1){\line(1,      0){30}}
\put(126,      -1){\line(1,      0){33}}
\put(159,      -1){\line(1,      0){33}}

\put(112,     1){\makebox(0,     0)[t]{$\cdots\cdots$}}

\put(22,     -20){$q$}
\put(40,      -15){$q ^{-1}$}
\put(58,      -20){$q$}
\put(75,      -15){$q ^{-1}$}

\put(91,      -20){$q$}

%\put(102,     -15){$q ^{-1}$}

\put(120,      -20){$q$}
\put(135,     -15){$q ^{-1}$}

\put(150,     -20){$q$}

\put(168,     -15){$q ^{-1}$}
\put(192,      -20){$q ^{}$}

\put(58,    38){\makebox(0,     0)[t]{$\bullet$}}

\put(58,      -1){\line(0,     1){40}}

\put(68,      33){$q$}

\put(68,      16){$q ^{-1}$}

\put(159,    38){\makebox(0,     0)[t]{$\bullet$}}

\put(159,      -1){\line(0,     1){40}}

\put(139,      33){$q$}

\put(139,      16){$q ^{-1}$}

\put(200,  -1)  { $q^{} \not=1$. }
\end{picture}$\\ \\ \\ \\ \\ \\

$\ \ \ \ \  \ \ \begin{picture}(130,      15)

\put(-45,      -1){$E^{(1)}_6$ }
\put(27,      1){\makebox(0,     0)[t]{$\bullet$}}
\put(60,      1){\makebox(0,      0)[t]{$\bullet$}}
\put(93,     1){\makebox(0,     0)[t]{$\bullet$}}

\put(126,      1){\makebox(0,    0)[t]{$\bullet$}}

\put(159,      1){\makebox(0,    0)[t]{$\bullet$}}

\put(126,      1){\makebox(0,    0)[t]{$\bullet$}}
\put(28,      -1){\line(1,      0){33}}
\put(61,      -1){\line(1,      0){30}}
\put(94,     -1){\line(1,      0){30}}
\put(126,      -1){\line(1,      0){33}}

\put(22,     10){$q ^{}$}
\put(40,      5){$q^{-1}$}
\put(58,      -20){${q}$}
\put(75,      -15){$q ^{-1}$}

\put(91,      -20){${q}$}

\put(102,    -15){$q ^{-1}$}

\put(128,      10){$q ^{}$}

\put(135,     5){$q ^{-1}$}

\put(161,      10){$q^{}$}

\put(91,    38){\makebox(0,     0)[t]{$\bullet$}}

\put(91,    78){\makebox(0,     0)[t]{$\bullet$}}

\put(91,      -1){\line(0,     1){80}}

\put(99,      16){$q ^{-1}$}

\put(99,      33){$q ^{}$}

\put(99,      56){$q ^{-1}$}

\put(99,      73){$q ^{}$}

\put(166,        -1)  { $q^{} \not=1$.. }
\end{picture}$\\ \\  \\

$\ \ \ \ \  \ \ \begin{picture}(130,      15)
\put(-45,      -1){$E^{(1)}_7$}
\put(27,      1){\makebox(0,     0)[t]{$\bullet$}}
\put(60,      1){\makebox(0,      0)[t]{$\bullet$}}
\put(93,     1){\makebox(0,     0)[t]{$\bullet$}}

\put(126,      1){\makebox(0,    0)[t]{$\bullet$}}

\put(159,      1){\makebox(0,    0)[t]{$\bullet$}}
\put(192,      1){\makebox(0,    0)[t]{$\bullet$}}
\put(225,      1){\makebox(0,    0)[t]{$\bullet$}}

\put(126,      1){\makebox(0,    0)[t]{$\bullet$}}
\put(28,      -1){\line(1,      0){33}}
\put(61,      -1){\line(1,      0){30}}
\put(94,     -1){\line(1,      0){30}}
\put(126,      -1){\line(1,      0){33}}
\put(159,      -1){\line(1,      0){33}}
\put(192,      -1){\line(1,      0){33}}

\put(22,     10){$q$}
\put(40,      5){$q^{-1}$}
\put(58,      -20){$q^{}$}
\put(75,      -15){$q^{-1}$}

\put(91,      -20){$q^{}$}

\put(102,     -15){$q^{-1}$}

\put(124,      -20){$q^{}$}
\put(135,     -15){$q^{-1}$}
\put(159,      10){$q^{}$}

\put(168,     5){$q^{-1}$}
\put(192,      10){$q^{}$}

\put(201,     5){$q ^{-1}$}
\put(225,      10){$q ^{}$}

\put(124,    38){\makebox(0,     0)[t]{$\bullet$}}

\put(124,      -1){\line(0,     1){40}}

\put(134,      33){$q^{}$}

\put(134,      16){$q^{-1}$}

\put(262,        -1)  {$q^{} \not=1$. }
\end{picture}$\\ \\ \\

$\ \ \ \ \  \ \
\begin{picture}(130,      15)\put(-45,      -1){$E^{(1)}_8$}
\put(-9,      1){\makebox(0,     0)[t]{$\bullet$}}

\put(27,      1){\makebox(0,     0)[t]{$\bullet$}}
\put(60,      1){\makebox(0,      0)[t]{$\bullet$}}
\put(93,     1){\makebox(0,     0)[t]{$\bullet$}}

\put(126,      1){\makebox(0,    0)[t]{$\bullet$}}

\put(159,      1){\makebox(0,    0)[t]{$\bullet$}}
\put(194,      1){\makebox(0,    0)[t]{$\bullet$}}

\put(126,      1){\makebox(0,    0)[t]{$\bullet$}}

\put(-10,      -1){\line(1,      0){33}}

\put(28,      -1){\line(1,      0){33}}
\put(61,      -1){\line(1,      0){30}}
\put(94,     -1){\line(1,      0){30}}
\put(126,      -1){\line(1,      0){33}}
\put(159,      -1){\line(1,      0){33}}

\put(-9,     -20){$q ^{}$}
\put(4,      -15){$q ^{-1}$}

\put(22,     -20){$q ^{}$}
\put(40,      -15){$q ^{-1}$}
\put(58,      -20){$q ^{}$}
\put(70,      -15){$q ^{-1}$}

\put(85,      -20){$q^{}$}

\put(102,    -15){$q ^{-1}$}

\put(122,      -20){$q$}

\put(135,     -15){$q ^{-1}$}

\put(161,      10){$q$}

\put(168,     5){$q ^{-1}$}

\put(194,      10){$q$}

\put(124,    32){\makebox(0,     0)[t]{$\bullet$}}

\put(124,      -1){\line(0,     1){34}}

\put(134,      33){$q ^{}$}

\put(134,      16){$q ^{-1}$}

\put(218,        -1)  {$q^{} \not=1$.
}
\end{picture}$\\ \\

$\ \ \ \ \  \ \ \begin{picture}(130,      15)

\put(-45,      -1){$F^{(1)}_4$}
\put(27,      1){\makebox(0,     0)[t]{$\bullet$}}
\put(60,      1){\makebox(0,      0)[t]{$\bullet$}}
\put(93,     1){\makebox(0,     0)[t]{$\bullet$}}

\put(126,      1){\makebox(0,    0)[t]{$\bullet$}}

\put(159,      1){\makebox(0,    0)[t]{$\bullet$}}
\put(28,      -1){\line(1,      0){33}}
\put(61,      -1){\line(1,      0){30}}
\put(94,     -1){\line(1,      0){30}}
\put(126,      -1){\line(1,      0){33}}

\put(22,     10){$q ^{}$}
\put(40,      5){$q ^{-1}$}
\put(58,      10){$q ^{}$}
\put(74,      5){$q ^{-1}$}

\put(91,      10){$q ^{}$}

\put(102,     5){$q ^{-2}$}

\put(124,      10){$q ^{2}$}
\put(135,     5){$q ^{-2}$}
\put(159,      10){$q ^{2}$}

\put(222,        -1)  { $q^{2} \not=1$. }
\end{picture}$\\

 $\ \ \ \ \  \ \ \  \begin{picture}(100,      15)  \put(-45,      -1){$G^{(1)}_2$ }

\put(60,      1){\makebox(0,      0)[t]{$\bullet$}}

\put(28,      -1){\line(1,      0){33}}
\put(27,      1){\makebox(0,     0)[t]{$\bullet$}}

\put(-14,      1){\makebox(0,     0)[t]{$\bullet$}}

\put(-14,     -1){\line(1,      0){50}}

\put(-18,     10){$q$}
\put(-5,      5){$q^{-1}$}
\put(22,     10){$q^{}$}
\put(38,      5){$q ^{-3}$}

\put(63,      10){$q^{3}$}

  \ \ \ \ \ \ \ \ \ \ \ \ \ \ \ \ \ \ \ {  ord ($q ) >3$. }
\end{picture}$\\

$\ \ \ \ \  \ \ \ \begin{picture}(100,      15)  \put(-45,      -1){  $A^{(2)}_2$   }

\put(60,      1){\makebox(0,      0)[t]{$\bullet$}}

\put(28,      -1){\line(1,      0){33}}
\put(27,      1){\makebox(0,     0)[t]{$\bullet$}}

%\put(-14,      1){\makebox(0,     0)[t]{$\bullet$}}

%\put(-14,     -1){\line(1,      0){50}}

\put(22,     10){$q ^{4}$}
\put(38,      5){$q^{-4}$}

\put(60,      10){$q ^{}$}

  \ \ \ \ \ \ \ \ \ \ \ \ \ \ \ \ \ \ \ { ord ($q ) >4,$ }
\end{picture}$\\

$\ \ \ \ \  \ \ \ \begin{picture}(100,      15)

\put(-45,      -1){$A^{(2)}_{2N}, N\ge 2$  }

\put(27,      1){\makebox(0,     0)[t]{$\bullet$}}
\put(60,      1){\makebox(0,      0)[t]{$\bullet$}}
\put(93,     1){\makebox(0,     0)[t]{$\bullet$}}
\put(159,      1){\makebox(0,      0)[t]{$\bullet$}}
\put(192,     1){\makebox(0,      0)[t]{$\bullet$}}
\put(28,      -1){\line(1,      0){30}}
\put(61,      -1){\line(1,      0){30}}
\put(130,     1){\makebox(0,     0)[t]{$\cdots\cdots\cdots\cdots$}}
\put(160,     -1){\line(1,      0){30}}
\put(22,     -15){}
\put(58,      -15){}
\put(91,      -15){}
\put(157,      -15){$$}
\put(191,      -15){$$}

\put(22,     10){$q^{4}$}

\put(40,      5){$q^{-4}$}

\put(58,      10){$q^{2}$}
\put(73,      5){$q^{-2}$}

\put(91,      10){$q^{2}$}

\put(157,      10){$q^{2}$}

\put(172,     5){$q^{-2}$}
\put(191,      10){$q$}

\put(210,        -1) {  ord ($q_{}) >4$.}
\end{picture}$\\
\\

$\ \ \ \ \  \ \ \ \ \ \ \begin{picture}(130,      15)

\put(-65,      -1){ $A^{(2)}_{2N-1}, N\ge 3$ }
\put(27,      1){\makebox(0,     0)[t]{$\bullet$}}
\put(60,      1){\makebox(0,      0)[t]{$\bullet$}}
\put(93,     1){\makebox(0,     0)[t]{$\bullet$}}

\put(126,      1){\makebox(0,    0)[t]{$\bullet$}}

\put(159,      1){\makebox(0,    0)[t]{$\bullet$}}

\put(126,      1){\makebox(0,    0)[t]{$\bullet$}}
\put(28,      -1){\line(1,      0){33}}
\put(61,      -1){\line(1,      0){30}}
%\put(94,     -1){\line(1,      0){30}}
\put(126,      -1){\line(1,      0){33}}

\put(104,     1){\makebox(0,     0)[t]{$\cdots\cdots$}}

\put(22,     10){$q^{2}$}
\put(40,      5){$q^{-2}$}
\put(58,      -20){$q^{2}$}
\put(75,      -15){$q ^{-2}$}

\put(91,      10){$q ^{2}$}

%\put(102,     5){$q ^{-2}$}

\put(128,      10){$q ^{2}$}

\put(135,     5){$q ^{-2}$}

\put(161,      10){$q ^{}$}

\put(58,    38){\makebox(0,     0)[t]{$\bullet$}}

\put(58,      -1){\line(0,     1){40}}

\put(68,      33){$q ^{2}$}

\put(68,      16){$q ^{-2}$}

\put(176,        -1)  { $q ^{2} \not= 1$. }
\end{picture}$\\ \\

$\ \ \ \ \  \ \ \ \ \ \ \begin{picture}(100,      15)

\put(-65,      -1){ $D^{(2)}_{N+1}, N\ge 2$ }

\put(27,      1){\makebox(0,     0)[t]{$\bullet$}}
\put(60,      1){\makebox(0,      0)[t]{$\bullet$}}
\put(93,     1){\makebox(0,     0)[t]{$\bullet$}}
\put(159,      1){\makebox(0,      0)[t]{$\bullet$}}
\put(192,     1){\makebox(0,      0)[t]{$\bullet$}}
\put(28,      -1){\line(1,      0){30}}
\put(61,      -1){\line(1,      0){30}}
\put(130,     1){\makebox(0,     0)[t]{$\cdots\cdots\cdots\cdots$}}
\put(160,     -1){\line(1,      0){30}}
\put(22,     -15){}
\put(58,      -15){}
\put(91,      -15){}
\put(157,      -15){$$}
\put(191,      -15){$$}

\put(22,     10){$q^{2}$}

\put(40,      5){$q^{-2}$}

\put(58,      10){$q^{}$}
\put(73,      5){$q^{-1}$}

\put(91,      10){$q^{}$}

\put(157,      10){$q^{}$}

\put(172,     5){$q^{-2}$}
\put(191,      10){$q^{2}$}

\put(210,        -1) { $q ^{2} \not= 1$.}
\end{picture}$\\

$\ \ \ \ \  \ \  \ \ \ \ \  \begin{picture}(130,      15)

\put(-65,      -1){$E^{(2)}_{6}$ }
\put(27,      1){\makebox(0,     0)[t]{$\bullet$}}
\put(60,      1){\makebox(0,      0)[t]{$\bullet$}}
\put(93,     1){\makebox(0,     0)[t]{$\bullet$}}

\put(126,      1){\makebox(0,    0)[t]{$\bullet$}}

\put(159,      1){\makebox(0,    0)[t]{$\bullet$}}
\put(28,      -1){\line(1,      0){33}}
\put(61,      -1){\line(1,      0){30}}
\put(94,     -1){\line(1,      0){30}}
\put(126,      -1){\line(1,      0){33}}

\put(22,     10){$q ^{2}$}
\put(40,      5){$q ^{-2}$}
\put(58,      10){$q ^{2}$}
\put(74,      5){$q ^{-2}$}

\put(91,      10){$q ^{2}$}

\put(102,     5){$q ^{-2}$}

\put(124,      10){$q ^{}$}
\put(135,     5){$q ^{-1}$}
\put(159,      10){$q ^{}$}

\put(222,        -1)  { $q^{2} \not=1$. }
\end{picture}$\\

 $\ \ \ \ \  \ \  \ \ \ \ \ \ \ \ \ \begin{picture}(100,      15)  \put(-85,      -1){$D^{(3)}_{4}, N\ge 2$ }

\put(60,      1){\makebox(0,      0)[t]{$\bullet$}}

\put(28,      -1){\line(1,      0){33}}
\put(27,      1){\makebox(0,     0)[t]{$\bullet$}}

\put(-14,      1){\makebox(0,     0)[t]{$\bullet$}}

\put(-14,     -1){\line(1,      0){50}}

\put(-18,     10){$q^{3}$}
\put(-5,      5){$q^{-3}$}
\put(22,     10){$q^{3}$}
\put(38,      5){$q ^{-3}$}

\put(63,      10){$q^{}$}

  \ \ \ \ \ \ \ \ \ \ \ \ \ \ \ \ \ \ \ {  ord ($q ) >3$. }
\end{picture}$

\end {Proposition}

All finite Cartan GDDs can be  obtained by omitting the first vertexes of
$A^{(1)}_N,$  $B^{(1)}_N,$ $C^{(1)}_N,$ $D^{(1)}_N,$ $ E^{(1)}_6,$$E^{(1)}_7,$ $E^{(1)}_8,$ $F^{(1)}_4,$ $G^{(1)}_2.$
%\subsection {About bi-classical, classical+ semi-classical and bi-semi-classical }
\subsection {About semi-classical, continual  and tail  }
\begin {Lemma}\label {2.84}   Assume rank $n >4.$  Tails of all  GDDs of Table C are listed.
  1 is a tail  of all  arithmetic  GDDs in Table C.
 5 is a tail  of  GDD $1$ of Row $11$.    5 is a tail  of  GDD $2$ of Row $11$.  5 is a tail  of    GDD $3$ of Row $11$.  5 is a tail  of    GDD $4$ of Row $11$.   5 is a tail  of   GDD $7$ of Row $12$.   5 is a tail  of   GDD $10$ of Row $12$.  5 is a tail  of   GDD $11$ of Row $12$.
 5 is a tail  of   GDD $13$ of Row $12$.  5 is a tail  of   GDD $14$ of Row $12$.
    5 is a tail   of  GDD $6$ of Row $14$.   5 is a tail   of  GDD $5$ of Row $15$.
    6 is a tail   of GDD $1$ of Row $16$.
    6 is a tail   of GDD $7$ of Row $17$. 6 is a tail   of  GDD $7$ of Row $19$.

\end {Lemma}

\begin {Lemma}\label {2.87}  Assume rank $n >4.$
 {\rm (i)}  All strict  semi-classical GDDs which Tails are  1 are listed:
 GDD $4$ of Row $11.$
          GDD $2$ of Row $12$.  GDD $3$ of Row $12$.  GDD $4$ of Row $12$.
                     GDD $8$ of Row $12$. GDD $10$ of Row $12$.
                              GDD $3$ of Row $13$.  GDD $5$ of Row $13$.
                                       GDD $6$ of Row $13$.
                                        GDD $8$ of Row $13$.
 GDD $9$ of Row $13$.
 GDD $10$ of Row $13$.  GDD $11$ of Row $13$.  GDD $12$ of Row $13$.  GDD $13$ of Row $13$.  GDD $14$ of Row $13$.
    GDD $4$ of Row $14$.
   GDD $6$ of Row $15$.
   GDD $5$ of Row $17$.  GDD $2$ of Row $18$.  GDD $3$ of Row $18$. GDD $4$ of Row $18$.
 GDD $5$ of Row $18$. GDD $7$ of Row $18$.  GDD $9$ of Row $18$. GDD $4$ of Row $19$.  GDD $3$ of Row $21$. GDD $4$ of Row $21$.

 {\rm (ii)} All strict semi-classical GDDs which Tails are not 1 are listed:
 5 is a tail   of  GDD $6$ of Row $14.$ \  5 is a tail   of  GDD $5$ of Row $15.$
\end {Lemma}

\begin {Lemma}\label {2.85}   Assume rank $n >4.$  Continuation of all non-classical GDDs of Table C are listed.
GDD $1$ of Row $11$ is continual  on 1:
GDD $1$ of Row $17$ in  T5 and  GDD $19$ of Row $18$ in  T6.

GDD $2$ of Row $11$ is continual  on 1:
 GDD $15$ of Row $18$ via T6  and  GDD $2$ of Row $17$ via T5.

GDD $3$ of Row $11$ is continual :
 GDD $4$ of Row $17$ via T5  on 5.   GDD $20$ of Row $18$ via T6  on 5.   GDD $3$ of Row $17$ via T5  on 1.  GDD $12$ of Row $18$ via T6  on 1.

 GDD $4$ of Row $11$ is continual :
 GDD $21$ of Row $18$ via T5  on 5.  GDD $6$ of Row $17$ via T6  on 5.  GDD $9$ of Row $18$ via T6  on 1. GDD $5$ of Row $17$ via T6  on 1.

 GDD $1$ of Row $12$ is continual  on 1:
 GDD $1$ of Row $18$ via T5.

GDD $2$ of Row $12$ is continual  on 1:
 GDD $2$ of Row $18$ via T5.

GDD $3$ of Row $12$ is continual  on 1:
 GDD $3$ of Row $18$ via T5

 GDD $4$ of Row $12$ is continual  on 1:
 GDD $4$ of Row $18$ via T5.

  GDD $5$ of Row $12$ is continual  on 1:
 GDD $6$ of Row $18$ via T5.

GDD $6$ of Row $12$ is continual  on 1:
 GDD $8$ of Row $18$ via T5

 GDD $7$ of Row $12$  is continual  on 1:
GDD $10$ of Row $18$ via T6

GDD $8$ of Row $12$ is continual  on 1:
 GDD $5$ of Row $18$ via T5

 GDD $9$ of Row $12$ is continual  on 1:
 GDD $11$ of Row $18$ via T5

 GDD $10$ of Row $12$ is continual  on 1:
 GDD $7$ of Row $18$ via T6

 GDD $11$ of Row $12$ is continual  on 1:
 GDD $13$ of Row $18$ via T6

GDD $12$ of Row $12$ is continual  on 1:
 GDD $14$ of Row $18$ via T5.

GDD $13$ of Row $12$ is continual  on 1:
 GDD $16$ of Row $18$ via T6

 GDD $14$ of Row $12$ is continual  on 1:
 GDD $17$ of Row $18$ via T6.

GDD $15$ of Row $12$ is continual  on 1:
 GDD $18$ of Row $18$ via T5.

GDD $1$ of Row $14$ is continual  on 1:
GDD $1$ of Row $19$ via T5

GDD $2$ of Row $14$ is continual  on 1:
 GDD $2$ of Row $19$ via T5

GDD $3$ of Row $14$ is continual  on 1:
GDD $3$ of Row $19$ via T5.

GDD $4$ of Row $14$ is continual  on 1:
GDD $4$ of Row $19$ via T5

GDD $5$ of Row $14$ is continual  on 1:
 GDD $5$ of Row $19$ via T5.

 GDD $6$ of Row $14$ is continual  on 1:
GDD $6$ of Row $19$ via T6

GDD $1$ of Row $16$ is continual  on 1:
 GDD $1$ of Row $20$ via T5. GDD $8$ of Row $21$ via T6

GDD $1$ of Row $17$ is continual  on 1:
 GDD $1$ of Row $21$ via T5.

GDD $2$ of Row $17$ is continual  on 1:
 GDD $2$ of Row $21$ via T5.

GDD $4$ of Row $17$ is continual  on 1:
GDD $5$ of Row $21$ via T5

GDD $5$ of Row $17$ is continual  on 1:
 GDD $4$ of Row $21$ via T5.

GDD $6$ of Row $17$is continual  on 1:
 GDD $6$ of Row $21$ via T5.

 GDD $7$ of Row $17$ is continual  on 1:
  GDD $7$ of Row $21$ via T6.

  GDD $4$ of Row $18$ is continual  on 1:
GDD $3$ of Row $21$ via T5.

 GDD $1$ of Row $20$ is continual  on 1:
GDD $1$ of Row $22$ via T5.

\end {Lemma}
\subsection {About proof of main result}
We  find all  quasi-affine  GDDs over an  arithmetic GDD   A as follows. We first find all  quasi-affine  GDDs which are GDDs adding   a vertex on vertex $1, 2, \cdots, n$, respectively, of A. For example,    We find all  quasi-affine  GDDs which are GDDs adding   a vertex on vertex $1$. Choice $i$, $2\le i \le n.$
Let B the GDD omitting vertex i  of A. We find all  GDDs  by adding  a vertex on vertex 1 of B such that they are  arithmetic. Next we put vertex 1 into these GDDs. Finally we determine if they are quasi-affine.

\end {document}